\documentclass[twoside,a4paper,reqno]{amsart}
\usepackage[colorlinks=true, breaklinks=true, urlcolor=webbrown, linkcolor=RoyalBlue, citecolor=webgreen, backref=page]{hyperref}

\usepackage{style}

\begin{document}

\title[Spectral continuity]{Direct and inverse spectral continuity for Dirac operators}

\author{R.\,V.\,Bessonov}
\author{P.\,V.\,Gubkin}

\address{
\begin{flushleft}
R.V.Bessonov: roman.bessonov@fmf.uni-lj.si\\\vspace{0.1cm}
University of Ljubljana\\  
Jadranska 19, 1000 Ljubljana, SLOVENIA\\
\vspace{0.1cm}
St.\,Petersburg Department of Steklov Mathematical Institute\\ Russian Academy of Sciences\\
Fontanka 27, 191023 St.Petersburg,  RUSSIA\\
\end{flushleft}
}

\address{
\begin{flushleft}
P.V.Gubkin: 
gubkinpavel@pdmi.ras.ru\\\vspace{0.1cm}
St.\,Petersburg State University\\  
Universitetskaya nab. 7-9, 199034 St.\,Petersburg, RUSSIA\\
\vspace{0.1cm}
St.\,Petersburg Department of Steklov Mathematical Institute\\ Russian Academy of Sciences\\
Fontanka 27, 191023 St.Petersburg,  RUSSIA\\
\end{flushleft}
}

\thanks{The work of RB in Sections 2, 4 is supported by grants P1-0291 and N1-0237 from Slovenian Research Agency ARIS. The work of PG in Sections 1, 3 is supported by the Russian Science Foundation grant RScF 19-71-30002. PG  is a winner of the BASIS PhD Student competition and would like to thank its sponsors and jury.}

\begin{abstract}
The half-line Dirac operators with $L^2$-potentials can be characterized by their spectral data. It is known that the spectral correspondence is a homeomorphism: close potentials give rise to close spectral data and vice versa. We prove the first explicit two-sided uniform estimate related to this continuity in the general $L^2$-case. The proof is based on an exact solution of the inverse spectral problem for Dirac operators with $\delta$-interactions on a half-lattice  in terms of the Schur's algorithm for analytic functions.
\end{abstract} 

\vspace{1cm}


\maketitle

\setcounter{tocdepth}{2}

\tableofcontents

\section{Introduction}\label{intro}
The problem of stable reconstruction of differential operators from their spectral data is a classical subject of spectral theory. Here, {\it stability} means that the nonlinear map relating coefficients of the operators and their spectral data is  a {\it homeomorphism} between some topological spaces.  An overview of  stability results for operators with discrete spectrum can be found in Hryniv \cite{Hr11, Hr11b}, Horv\'ath and Kiss \cite{HK10}; see also Savchuk and Shkalikov \cite{SS10}, Chelkak, Kargaev, and Korotyaev \cite{KCH09, ChKK04}. Spectral stability for operators with nonempty continuous spectrum is much less studied. We give some historical remarks below in Section \ref{sHR}. This paper is devoted to the stability of the solutions of direct and inverse spectral problems for Dirac operators on the half-line $\R_+ = [0, +\infty)$.

\subsection{Dirac operators} The Dirac operator $\Di_{q}$ on $\R_+$ with a square summable potential $q \in L^2(\R_+)$ is defined by the differential expression
\begin{equation}\label{do}
\Di_{q}: X \mapsto JX' + QX, \qquad Q = \begin{pmatrix}	\Im q & \Re q \\	\Re q & -\Im q \end{pmatrix},
\end{equation}
where $J = \jm$, and $X$ belongs to the set of all locally absolutely continuous vector-valued functions in the Hilbert space 
$$
L^2(\R_+, \C^2) = \Bigl\{Y: \R_+ \to \C^2:\; \|Y\|_{L^2(\R, \C^2)}^{2} = \int_{\R_+}\|Y(x)\|_{\C^2}^{2}\,dx < \infty\Bigr\}
$$ 
satisfying a self-adjoint boundary condition $\left\langle X(0),\left(\!
\begin{smallmatrix}
\sin \alpha \\ \cos \alpha
\end{smallmatrix}
\!\right)\right\rangle_{\C^2} = 0$, and such that $\Di_{q}X \in L^2(\R_+, \C^2)$.
With this domain, $\Di_{q}$ is the densely defined self-adjoint operator on $L^2(\R_+, \C^2)$, see  \cite{Levitanb}, \cite{Den06}. Without loss of generality, we will work with the boundary condition  corresponding to $\alpha = 0$. To describe appropriate spectral data for $\Di_{q}$, consider the matrix-valued solution $N = \left(\begin{smallmatrix}
n_{11} & n_{12} \\
n_{21} & n_{22} 
\end{smallmatrix}\right)$ of the Dirac system $JN'(x,z) + Q(x)N(x, z) = zN(x, z)$, $N(0, z) = \idm$, $z \in \C$. The Weyl function of $\Di_{q}$ is defined by
$$
m_{q}(z) = \lim_{x \to +\infty} \frac{n_{22}(x,z)}{n_{21}(x,z)}, \qquad \Im z > 0.
$$
This function belongs to the Herglotz class in the upper half-plane $\C_+ = \{z \in \C:\; \Im z > 0\}$, i.e., it is analytic and takes $\C_+$ into itself. See \cite{Den06} for the Weyl theory of Dirac operators from the perspective of Krein systems. It is well-known that the Weyl function $m_q$ determines $\Di_{q}$ uniquely. From the point of view of spectral correspondence it is more convenient to work with the Cayley transform of $m_q$, i.e., with the Schur function $f_q$ of $\Di_{q}$. Recall that an analytic function $f$ on $\C_+$ is said to belong to the Schur class $S(\C_+)$ if $|f(z)| \le 1$ for all $z \in \C_+$. The Schur function $f_{q}$ of $\Di_{q}$ is determined by
\begin{equation}\label{sf}
m_{q} = i\frac{1+f_{q}}{1-f_{q}}.
\end{equation}
Each function $f \in S(\C_+)$ has the nontangential boundary values almost everywhere on $\R$ \cite{Garnett}. As usual, we use the same letter $f$ for the function in the unit ball of $L^\infty(\R)$ defined by these boundary values. In a moment we will see that Schur functions of $\Di_{q}$ with $q \in L^2(\R_+)$ belong to the set 
\begin{equation}\label{eq: spaceS}
S_{2}(\C_+) = \big\{f \in S(\C_+):\; \log(1-|f|^2) \in L^1(\R) \big\}.
\end{equation} 
The set $S_{2}(\C_+)$ is a complete metric space with respect to the metric
\begin{equation}\label{eq42}
\rho_{S_{2}}(f, g) = \sqrt{\int_{\R}-\log\Bigl(1 - \Bigl|\frac{f - g}{1 - \ov{f}g}\Bigr|^2\Bigr)\,dx}.
\end{equation} 
It can be shown that $f_n \to g$ in $S_{2}(\C_+)$ if and only if $\|\log(1-|f_n|^2)\|_{L^1(\R_+)} \to \|\log(1-|g|^2)\|_{L^1(\R_+)}$ and $f_n \to g$ in Lebesgue measure on $\R$, see Lemma \ref{lem: conv in S2 and conv in measure v2}.  

\subsection{Sylvester-Winebrenner theorem} Our starting point is the following fundamental result that stems from the paper \cite{Sylvester} by Sylvester and Winebrenner. For Dirac operators, it was first proved by Denisov \cite{Den06} in a somewhat restricted form, see details in Section \ref{section: tSWD} below. 
\begin{Thm}[Sylvester-Winebrenner theorem]\label{tSW}
The correspondence $\F: q \mapsto f_{q}$ is a homeo\-morphism from $L^2(\R_+)$ onto $S_{2}(\C_+)$. Moreover, we have
\begin{equation}\label{eqSW} 
\int_{\R_+}|q(x)|^2\,dx = \frac{1}{\pi}\int_{\R}-\log(1-|f_q(x)|^2)\,dx.	
\end{equation}
\end{Thm}
Theorem \ref{tSW} belongs to a general direction in spectral theory that relies on the usage of trace formulae (or {\it sum rules}, in the terminology of B.\,Simon \cite{SimonDes}). This direction often leads to the most general results when one is interested in complete characterization theorems
({\it ``spectral gems''} \cite{SimonDes}). See, e.g., \cite{KS03, KS09, DKS10, Yu18, BD20, BD21b, DEY21}. The proofs of such theorems, however, do not involve reconstruction procedures for potentials from the spectral data, and in particular, they do not imply any {\it continuity estimates} related to the spectral correspondence. To illustrate the situation, let us rewrite the sum rule \eqref{eqSW} in the form
$$
\|q-0\|_{L^2(\R_+)}^{2} = \rho_{S_2}(f_q, 0)^{2}/\pi.
$$
Having this identity, it is natural to expect that quantities $\|q - \tilde q\|_{L^2(\R_+)}$, $\rho_{S_2}(f_q, f_{\tilde q})$ control each other  for $q$, $\tilde q \in L^2(\R)$. Moreover, Theorem \ref{tSW} says that $\|q_n - \tilde q\|_{L^2(\R_+)} \to 0$ if and only if $\rho_{S_2}(f_{q_n}, f_{\tilde q}) \to 0$, making this expectation even more plausible. It turns out, however, that it is false. In fact, we have the following theorem. 
\begin{Thm}\label{t5}
There are potentials $u_n$, $\tilde u_n$, $q_n$, $\tilde q_n$ in the unit ball of $L^2(\R_+)$ such that
\begin{align} 
	&\lim_{n \to \infty}\|u_n - \tilde u_n\|_{L^2(\R_+)} = 0, \;\; \mbox{but} \;\; \lim_{n \to \infty} \rho_{S_{2}}(f_{u_n}, f_{\tilde u_n}) > 0, \label{eq12}\\
	&\lim_{n \to \infty}\|q_n - \tilde q_n\|_{L^2(\R_+)} \,> 0, \;\; \mbox{but} \;\;\, \lim_{n \to \infty} \rho_{S_{2}}(f_{q_n}, f_{\tilde q_n}) = 0. \label{eq11}
\end{align}
In other words, the homeomorphisms $\F$, $\F^{-1}$ in Theorem \ref{tSW} are not uniformly continuous on bounded subsets of $L^2(\R_+)$, $S_2(\C_+)$. 
\end{Thm}
As an ``explanation'' for \eqref{eq12}, \eqref{eq11}, let us mention that the continuous operators $\F$, $\F^{-1}$ are not linear and the closed unit balls in $L^2(\R_+)$, $S_{2}(\C_+)$ are not compact. In particular, the standard general arguments are not applicable here and the lack of uniform continuity of $\F$, $\F^{-1}$ is possible. To prove \eqref{eq12}, we construct some explicit sequences of potentials $u_n$, $\tilde u_n$. The proof of \eqref{eq11} is more delicate. Here we use a very important observation of Volberg and Yuditskii \cite{VYu02} on the non-injectivity of the scattering map for Jacobi matrices. This observation was transferred to the setting of the nonlinear Fourier transform by Tao and Thiele \cite{TT} and to Dirac operators by the first author and Denisov \cite{BD24}.

\medskip

\subsection{The main result} Let us now turn to positive results. For every $f \in S_{2}(\C_+)$, the function $|f|^2$ is comparable to $|\log(1-|f|^2)|$ on the set $E$ where $|f| \le 1/2$, and the complement $\R \setminus E$ has a finite measure as $\log(1-|f|^2) \in L^1(\R)$. It follows that $f \in L^2(\R)$. Thus, the Fourier transform of any element $f \in S_{2}(\C_+)$ is well defined and belongs to $L^2(\R)$. We will denote it by $\hat f$, so that
\begin{gather}
\label{eq: fourier on R def}
\hat f(\xi) = \frac{1}{\sqrt{2\pi}}\int_{\R}f(x)e^{-i\xi x}\,dx
\end{gather}
if $f$ is integrable. We will need the following Wiener-type norm and the weighted $L^1$-norm:
$$
\|f\|_{W^1_A(\R_+)} = \int_{\R_+}|\hat f(\xi)|e^{-A\xi}\,d\xi, \qquad \|q\|_{L^1_A(\R_+)} = \int_{\R_+}|q(\xi)|e^{-A\xi}\,d\xi.
$$
Our main result is the following theorem. 
\begin{Thm}\label{t1}
Let $q$, $\tilde q \in L^2(\R_+)$, and let $f_{q}, f_{\tilde q}$ be the Schur functions \eqref{sf} of the corresponding Dirac operators $\Di_{q}$, $\Di_{\tilde q}$ \eqref{do}. Then we have 
\begin{equation}\label{t1est}
	c_{1}\|q - \tilde q\|_{L^1_{2A}(\R_+)} \le \|f_{q} - f_{\tilde q}\|_{W^1_{A}(\R_+)}
	\le c_2\|q - \tilde q\|_{L^1_{2A}(\R_+)}
\end{equation}
for $c_1 = \sqrt{\pi/2}$, $c_2 = 2\sqrt{2\pi}$ and any $A \ge 12\max\big(\|q\|_{L^2(\R_+)}^2,  \|\tilde q\|^2_{L^2(\R_+)}\big)$. 
\end{Thm}
We would like to stress that \eqref{t1est} is a uniform estimate. This makes it much stronger than the continuity property in Theorem \ref{tSW}, cf. \eqref{eq12}, \eqref{eq11}. Note also that \eqref{t1est} is nontrivial and new even in the case $\tilde q = 0$.

\subsection{Historical remarks}\label{sHR} Before proceed with further results, let us give a few historical remarks. Perhaps, the most general stability result in the one-dimensional spectral theory is the Krein-de Branges spectral theorem (see Section 5.2 of \cite{Remlingb}). It gives stability of the solution of direct and inverse spectral problems for canonical Hamiltonian systems and implies spectral stability for various other classical operators (Schr\"odinger and Dirac operators, Krein strings, Jacobi matrices). However, the proof of Krein-de Branges theorem, at least in its current form, cannot give explicit stability estimates, because it uses the following  general topological argument to prove the fact that the solution map is a homeomorphism: 
\begin{equation}\label{continuity}
\mbox{\it a continuous bijection between two Hausdorff compacts is a homeomorphism.}
\end{equation}
Moreover, the usage of compactness arguments similar to \eqref{continuity} forces to deal with very weak variants of stability, because for this approach closed bounded subsets in the metric spaces under consideration need to be compact. Even implicit stability control with respect to the norms like $\|\cdot \|_{L^p(\R_+)}$ or $\|\cdot\|_{L^1_A(\R_+)}$ via general Krein-de Branges theory is not possible because these norms define topologies that are not locally compact.

\medskip

For classical operators, there are well-known constructive methods to solve inverse spectral problems developed by Gelfand-Levitan and Krein (see Marchenko \cite{March06} for an excellent historical overview) and more recent methods by Belishev-Mikhaylov \cite{BM12} and Makarov-Poltoratski \cite{MP23}. In principle, one can use these methods to prove spectral continuity by accurate estimation of all quantities appearing in the proofs. See, e.g.,\ Section 6 in \cite{Den06} for some spectral continuity results for $L^2_{\loc}(\R_+)$-potentials.  However, it seems difficult to get optimal estimates on this way, because the methods are rather involved and even ``just'' characterization theorems (isomorphisms between potentials and spectral data without proving continuity) are already very nontrivial.  This explains why we choose another road and first treat the case of the so-called {\it exactly solvable models}. In this case a detailed analysis is possible and gives explicit two-sided estimates. We then approximate a general Dirac operator with potential $q \in L^2(\R_+)$ by a sequence of exactly solvable models and arrive at Theorem \ref{t1}. Our main instrument in the analysis of exactly solvable models is the classical Schur's algorithm, whose definition we now recall.

\subsection{Schur's algorithm}\label{section: schur algorithm} The Schur class $S(\D)$ in the open unit disk $\D = \{z \in \C: |z|<1\}$ consists of analytic functions $F$ on $\D$ such that $|F(z)| \le 1$ for all $z \in \D$. We will deal with the subset of $S(\D)$ defined by
$$
S_*(\D) = \{F \in S(\D)\colon \; \text{$F$ is not a finite Blaschke product}\}.
$$ 
The Schur's algorithm for $F \in S_*(\D)$ is the following iterative procedure: 
\begin{equation}\label{sa}
F_0 = F, \qquad zF_{k+1} = \frac{F_k - F_k(0)}{1 - \ov{F_k(0)}F_k}, \qquad k\ge 0, \qquad z \in \D.
\end{equation}
All functions produced by Schur's algorithm belong to $S_*(\D)$, moreover, each step $F_{k} \mapsto F_{k+1}$ is the bijection from $S_*(\D)$ onto itself. The numbers $\{F_k(0)\}_{k \ge 0}$ in \eqref{sa} are called the recurrence coefficients of $F$.  We also will need a version of Schur's algorithm for periodic functions in the the upper half-plane $\C_+$. Fix $\ell > 0$ and define
\begin{equation}\label{eq103}
S_{\ell, *}(\C_+) = \bigl\{f:\C_+ \to \D: \;\; f(z) = F(e^{2i\ell z}),\; F \in S_*(\D)\bigr\}.
\end{equation}
A function $f \in S(\C_+)$ belongs to $S_{\ell, *}(\C_+)$ if and only if $f(z+\pi/\ell) = f(z)$ for every $z \in \C_+$, and there is no finite Blaschke product $B$ in $\D$ such that $f = B(e^{2i\ell z})$. Relation \eqref{sa} for $f(z) = F(e^{2i\ell z})$, $F \in S_{*}(\D)$, takes the form
\begin{equation}\label{sa-hp}
f_0 = f, \qquad e^{2i\ell z}f_{k+1} = \frac{f_k - f_k(\infty)}{1 - \ov{f_k(\infty)}f_k}, \qquad k \ge 0, \qquad z \in \C_+,
\end{equation}
where $f_k(\infty) = \lim_{y \to +\infty}f_k(iy) = F_{k}(0)$, see Lemma \ref{l10}. It is natural to call the numbers $\{f_{k}(\infty)\}_{k \in \Z_+}$ the recurrence coefficients of $f \in S_{\ell, *}(\C_+)$. They determine functions $f$, $F$ uniquely. In fact, the knowledge of first $n$ recurrence coefficients of $f$ allows to approximate it with accuracy $2e^{-2n\ell y}$ in the half-plane $\Im z \ge  y > 0$, cf. $(1.3.43)$ in \cite{Simonbook1}.

\medskip

\subsection{Kronig-Penney model} Let us consider the half-line Dirac operators $\Di_{q}$ \eqref{do} on $\R_+$ whose potentials 
\begin{equation}\label{dp}
q = \sum_{k \in \Z_+}c_k \delta_{\ell k}, \qquad c_k \in \C,
\end{equation} 
now are linear combinations of point masses (usually, they are called $\delta$-interactions) supported on the half-lattice $\ell\Z_+ = \{\ell k:\;  k \in \Z, \; k \ge 0\}$ of sparseness $\ell > 0$. The spectral theory of this class of Dirac operators can be considered in the framework of the relativistic Kronig-Penney model for massless particles (for comparison, in the original Kronig-Penney model \cite{KP31} for a non-relativistic electron in a one-dimensional crystal, Schr\"odinger operators with $\delta$-interactions on the lattice $\Z$ were used, and $c_k \equiv c$ in the classical case). The Kronig-Penney model and its various generalizations are called exactly solvable models, meaning that the resolvents of the operators under consideration (Schr\"odinger or Dirac) as well as associated quantities (spectra, resonances, eigenfunctions, etc) can be often found explicitly in terms of the potential $q$. Solvable models attracted an enormous attention in theoretical physics and mathematics. We refer the reader to the classical monograph \cite{AGHH88} (1988) by Albeverio, Gesztesy, H{\o}egh-Krohn, and Holden, to the review chapter by Exner in the second edition \cite{AGHH05} (2005) of this monograph, and to later survey by Kostenko and Malamud \cite{KM13} (2013). The most close mathematical references to our work are \cite{CMP, LS14, Hug98, GS87}. The main distinction of our setting from the previous considerations comes from the fact that exactly solvable relativistic models are usually studied for radial massive Dirac operators
$$
\Di_{\mathbf{m}, q}: X \mapsto c\begin{pmatrix}0 & -1 \\ 1 & 0\end{pmatrix} X'(x) + \mathbf{m}c^2 \begin{pmatrix}1 & 0 \\ 0 & -1\end{pmatrix} X(x) + Q(x) X(x), 
$$	
with some positive parameters $c$, $\mathbf{m}$ (corresponding to the velocity of light and the mass of the particle, see Section 4.6.6 in \cite{Thaller}). We consider the case where $\mathbf{m} = 0$ and choose physical units so that $c=1$. We also assume that $Q$ is real and written in the second canonical form, i.e., $Q = Q^*$, $\trace Q = 0$, see~\eqref{do}. These assumptions are standard for the spectral theory of Dirac operators \cite{LSb} and for its applications to the nonlinear Schr\"odinger equation \cite{FTbook} (the massless Dirac operator is the auxiliary operator for the inverse scattering transform method for NLSE) but less common in the area of exactly solvable models.

\medskip

Essential part of the literature devoted to exactly solvable models deals with direct problems: knowing potential $q$ (a measure on some discrete subset of $\R$), one determines some spectral characteristics of the corresponding Schr\"odinger or Dirac operator. On the other hand, the full spectral characterization (the Weyl function or the spectral measure and the corresponding Fourier transform) is not known even for potentials $q$ supported on a lattice in the simplest massless case $\mathbf{m} = 0$. Moreover, it is not immediate if it is possible to describe spectral measures in closed form in terms of $q$. Indeed, for general $q$ of the form \eqref{dp} the corresponding spectral measures could have a complex structure and arbitrary spectral type (e.g., singular continuous component is not excluded). Below we show that such a description indeed exists, and, moreover, it turns out to be very simple and explicit (modulo nonlinearity and generality of the problem). To state our second main result, we need the following bijection $\varkappa: \D \to \C$: 
\begin{equation}\label{eq44}
\varkappa(w) = \frac{\ov{w}}{2|w|}\log\frac{1+|w|}{1-|w|}, \qquad w \in \C.
\end{equation}
\begin{Thm}\label{t2}
Let $q$ be a discrete complex-valued measure on $\R_+$ such that $\supp q \subset \ell \Z_{+}$, $\ell > 0$. Then the Schur function $f_q$ of $\Di_{q}$ is a nondegenerate $\pi/\ell$-periodic function, i.e.,\ $f_q \in S_{\ell, *}(\C_+)$, where the set $S_{\ell, *}(\C_+)$ is defined in \eqref{eq103}. Any element of $S_{\ell, *}(\C_+)$ arises uniquely in this way.  Moreover, $q$ and $f_q$ determine each other via
\begin{equation}\label{eq39}
	q(\{\ell k\}) = \varkappa(f_{q,k}(\infty)), \qquad k \in \Z_+,
\end{equation}
where $\{f_{q,k}(\infty)\}_{k \in \Z_+}$ is the sequence of recurrence coefficients of $f_q$, see \eqref{sa-hp}.
\end{Thm}
It is astonishing that the Schur's algorithm \cite{Schur17} (1917) and the Kronig-Penney model \cite{KP31} (1931) met in Theorem \ref{t2} almost a century after their independent development and generalizations in various directions.  In fact, Theorem~\ref{t2} completely reduces the spectral theory of the relativistic Kronig-Penney model for massless particles to the theory of orthogonal polynomials on the unit circle \cite{Szbook, Simonbook1}. To illustrate this, we formulate some immediate spectral consequences. Below we assume that $q$ is an arbitrary discrete complex-valued measure on $\R_+$ such that $\supp q \subset \Z_{+}$, and $\Di_{q}$ is the corresponding Dirac operator. We will use the standard notation  $\sigma_{ac}(\Di_{q})$, $\sigma_{p}(\Di_{q})$, $\sigma_{sc}(\Di_{q})$ for the absolutely continuous, pure point, and singular continuous parts of its spectrum. 
\begin{Cor}[Rakhmanov]\label{c3}
If $\sigma_{ac}(\Di_{q}) = \R$, then $q(\{k\}) \to 0$ as $k \to +\infty$.
\end{Cor}
\begin{Cor}[Baxter]\label{c4}
If $\sum_{k \ge 0}|q(\{k\})| < \infty$, then $\sigma_{p}(\Di_{q}) = \sigma_{sc}(\Di_{q}) = \emptyset$, and the main spectral measure $\mu_q$ of $\Di_{q}$ has the form $\mu_q = w_q\,dx$, where $w_q = \sum_{k \in \Z}c_k e^{2ikx}$ is such that $\sum_{k \in \Z}|c_k| < \infty$, $\min_\R w_q > 0$. The converse is also true. 
\end{Cor}
\begin{Cor}[Szeg\H{o}-Golinskii-Ibragimov]\label{c5}
If $\sum_{k \ge 0}k|q(\{k\})|^2 < \infty$, then $\sigma_{p}(\Di_{q}) = \sigma_{sc}(\Di_{q}) = \emptyset$, and the main spectral measure $\mu_q$ of $\Di_{q}$ has the form $\mu_q = w_q\,dx$, where $\log w_q = \sum_{k \in \Z}c_k e^{2ikx}$ is such that $\sum_{k \in \Z}|k||c_k|^2 < \infty$. The converse is also true. 
\end{Cor}
Generalized versions of Schur's algorithm appeared in spectral theory previously. For example, in Section 10 of \cite{Den06} Denisov discusses an analogue of Schur's algorithm in the form of the differential equation
$$
\frac{df_\ell(z)}{d\ell} = -iz f_\ell(z) + A(\ell) - \overline{A(\ell)}f_\ell^2(z), \qquad \ell \in \R_+, \quad z \in \C_+, \quad A(\ell) = q(\ell/2)/2,
$$ 
for Schur functions $f_\ell$ of potentials $q_\ell = q(\ell + \cdot)$ generated by $q \in L^1_{\loc}(\R_+)$, see $(10.3)$ in \cite{Den06}. See also Poltoratski \cite{Polt24} for a similar Ricatti equation for certain meromorphic inner functions arising on the spectral side of Dirac system. These equations are more difficult to analyze in the perturbation regime than the recursive relation \eqref{sa-hp} in the classical Schur's algorithm.  

\medskip

\subsection{Plan of the paper} We start the next section with the introduction of canonical Hamiltonian systems and Dirac operators with  measures. We prove Theorem \ref{t2} in Section \ref{s32}. Theorem \ref{t1} is proved in Section \ref{s4}. We derive it from continuity estimates for Schur's algorithm (Section \ref{s41}) and Theorem \ref{t2}. In the last section, we show that the homeomorphism in Theorem \ref{tSW} is not uniformly continuous in both directions, thus proving Theorem \ref{t5}.

\section{Dirac operators with periodic Weyl functions}\label{s5}
\subsection{Canonical Hamiltonian systems}
A canonical Hamiltonian system on $\R_+ = [0, +\infty)$ is the differential equation of the form
\begin{equation}\label{csv}
JX'(x, z) = z\Hh(x) X(x, z), \qquad x \in \R_+, \quad z \in \C,
\end{equation}
where, as before, $J = \jm$, the derivative of the vector-valued function $X: \R_+ \times \C \to \C^2$ is taken with respect to $x \in \R_+$, and $\Hh$ is a Hamiltonian. Here and below by a Hamiltonian we mean a matrix-valued mapping on $\R_+$, 
$$
\Hh: x \mapsto \begin{pmatrix}
h_1 & h\\h&h_2
\end{pmatrix}, \qquad \Hh(x) \ge 0,
$$
whose entries $h_1$, $h_2$, $h$ are real-valued functions in $L_{\loc}^1(\R_+) = \{f:\; f \in L^1[0,\ell]\mbox{ for all }\ell >0\}$. We also assume that $\Hh$ is not identically zero on any subset of $\R_+$ of positive Lebesgue measure. A Hamiltonian $\Hh$ is called singular if $\Hh \notin L^1(\R_+)$, or, equivalently,
$$
\int_{\R_+}\trace{\Hh(x)}\,dx = +\infty.
$$ 
The set of singular Hamiltonians on $\R_+$ will be denoted by $\Hbbs$. Let $\Theta$, $\Phi$ denote the solutions of \eqref{csv} satisfying $\Theta(0, z) = \oz$, $\Phi(0,z) = \zo$. Take $\omega \in \C_+ \cup \R \cup\{\infty\}$ (here and below $\infty$ is regarded as the element of the Riemann sphere, and the linear fractional transformations involving $\infty$ are understood accordingly). To each Hamiltonian $\Hh\in \Hbbs$ one can associate the Weyl function, 
\begin{equation}\label{wf}
m_{\Hh}(z) = \lim_{x \to +\infty}\frac{\Phi^+(x,z) + \Phi^-(x,z)\omega}{\Theta^+(x,z) + \Theta^-(x,z)\omega},
\qquad
\Theta = \begin{pmatrix}
	\Theta^+\\
	\Theta^-
\end{pmatrix}, \quad \Phi = \begin{pmatrix}
	\Phi^+\\
	\Phi^-
\end{pmatrix}, \quad z \in \C_+.
\end{equation}
It is known that the limit above exists for every $z \in \C_+$ and does not depend on the choice of $\omega$.
Moreover, the Weyl function, $m_{\Hh}$, is analytic in $\C_+$ and takes $\C_+$ into $\C_+$ unless it coincides with a constant $c \in \R \cup \{\infty\}$. Weyl's theory for canonical Hamiltonian systems can be found in \cite{HSW2000}, \cite{Remlingb}, \cite{Romanov}. We also define the Schur's  function $f_\Hh$ by 
\begin{gather}\label{sf cs}
f_\Hh = \frac{m_\Hh - i}{m_\Hh + i},\qquad\qquad m_\Hh = i\frac{1+f_\Hh}{1-f_\Hh}.
\end{gather}

Analytic functions taking $\C_+$ into $\C_+$ form the Herglotz-Nevanlinna class $N(\C_+)$. The set  
$$
\ov{N}(\C_+) 
= N(\C_+) \cup \R \cup \{\infty\}
$$ is the compactification of $N(\C_+)$ when the latter is equipped with the topology of convergence on compact subsets in $\C_+$. This topology (we extend it to $\ov{N}(\C_+)$) is metrizable with the metric, e.g.,
$$
\rho_{c}(m, \tilde m) = \max_{|z-i|\le 1/2}\frac{2|m(z) - \tilde m(z)|}{\sqrt{1+|m(z)|^2}\sqrt{1+|\tilde m(z)|^2}}, \qquad \rho_{c}(m, \infty) = \max_{|z-i| \le 1/2}\frac{2}{\sqrt{1+|m(z)|^2}}.
$$
One can choose other metrics on $\ov{N}(\C_+)$ determining the same compact topological space, see discussion on page 109 in \cite{Remlingb}.

\medskip

Different singular Hamiltonians $\Hh, \tilde\Hh$ can have equal Weyl functions. For instance, it is not difficult to check that if  $\tilde\Hh(x) = \xi'(x)\Hh(\xi(x))$ almost everywhere on $\R_+$ for some locally absolutely continuous increasing bijection $\xi: \R_+ \to \R_+$, then $m_{\Hh} = m_{\tilde\Hh}$. It will be convenient to call such $\Hh$, $\tilde \Hh$ equivalent, so that $\Hbbs$ becomes the set of classes of equivalent Hamiltonians: 
$$
\Hbbs = \{\Hh\mbox{ is a singular Hamiltonian on }\R_+\}/\sim.
$$
One can check that each class of equivalence in $\Hbbs$ contains the unique (up to values on a set of measure zero) element $\Hh^{tr}$ such that $\trace \Hh^{tr} = 1$ on $\R_+$. One can turn $\Hbbs$ into a compact Hausdorff space by defining the topology via the metric, e.g., 
\begin{equation}\label{eq8}
d(\Hh, \tilde\Hh) = \sum_{n \ge 0}2^{-n}\frac{d_n(\Hh^{tr}, \tilde\Hh^{tr})}{1+d_n(\Hh^{tr}, \tilde\Hh^{tr})}, \qquad d_n(\Hh^{tr}, \tilde\Hh^{tr})= \sup_{0 \le t \le n}\left\|\int_{0}^{t}(\Hh^{tr}(s) - \tilde \Hh^{tr}(s)) \,ds\right\|.
\end{equation}
Compactness of $\Hbbs$ follows from Riesz representation theorem, see the proof of Theorem 5.4 in \cite{Remlingb}.

\medskip

The following theorem is a key result of Krein -- de Branges spectral theory of canonical Hamiltonian systems \cite{KK68}, \cite{dBbook}, \cite{DMbook}, \cite{Remlingb}, \cite{Romanov}. 
\begin{Thm}[Krein -- de Branges theorem]\label{KdB}
The correspondence $\Hh \mapsto m_{\Hh}$ is the homeomorphism of compact metric spaces $\Hbbs$,  $\ov{N}(\C_+)$.
\end{Thm}
For a discussion (and some surprisingly deep applications) of the continuity part in Theorem \ref{KdB}, see Section 5 in \cite{Remlingb}, Section 2 in \cite{EKT18}, or Section 3 in \cite{ELS21}. Currently, explicit estimates related to the continuity properties of the homeomorphism in Theorem \ref{KdB} are not known. However, in contrast to the Sylvester-Winebrenner theorem (Theorem \ref{tSW}), the homeomorphism in Krein-de Branges theorem are uniformly continuous in both directions by compactness and the Heine-Cantor theorem. 

\medskip

Let $\Hh \in \Hbbs$ be a singular Hamiltonian on $\R_+ = [0, +\infty)$. The fundamental matrix solution corresponding to $\Hh$ is the locally absolutely-continuous (with respect to $x \in \R_+$) matrix-valued mapping $M: \R_+ \times \C \to \sltc$ satisfying the differential equation
\begin{equation}\label{cs}
JM'(x, z) = z\Hh(x) M(x,z), \qquad M(0, z) = \idm,
\end{equation}
almost everywhere on $\R_+$, where the differentiation is taken with respect to $x \in \R_+$ and $z \in \C$ is a spectral parameter. Note that we have $M = \bigl(\Theta, \Phi \bigr)$ in terms of the solutions $\Theta$, $\Phi$ of \eqref{csv} satisfying $\Theta(0, z) = \oz$, $\Phi(0,z) = \zo$. Let us write $w_1 \doteq A w_2$ for two complex numbers $w_1$, $w_2$ and a matrix $A$ with $\det A \neq 0$ if 
$$
w_1 = \frac{a_{11}w_2 + a_{12}}{a_{21}w_2 + a_{22}}, \qquad A = \begin{pmatrix}
a_{11} & a_{12}\\ a_{21} & a_{22}
\end{pmatrix}.
$$
It is straightforward to generalize this definition to the case where $w_1$, $w_2$ can admit the value $\infty$. For a matrix $A \in \sltc$, we have $w_1 \doteq A w_2$ if and only if $w_2 \doteq A^{-1} w_1$, $w_{1,2} \in \C \cup \{\infty\}$. Note also that the definition \eqref{wf} for $m_{\Hh}$ can be rewritten in the following form:
\begin{equation}\label{wfbis}
m_{\Hh} \doteq \lim_{x \to \infty}\sigma_1 M(x,z)^{T} \sigma_1 \omega,
\end{equation}
where $\sigma_1 = \left(\begin{smallmatrix}
0 &1\\1 & 0\end{smallmatrix}\right)$ and $M^{T}$ stands for the transposed matrix $M$.  Relation \eqref{sf cs} between the Weyl and Schur functions reads as
\begin{gather}
\label{eq: LR matrices}
m_q \doteq 
Lf_q,\quad f_q \doteq 
Rm_q,\qquad
L = 
\begin{pmatrix}
	i & i\\
	-1 & 1
\end{pmatrix},
\quad 
R = 
\begin{pmatrix}
	1 & -i\\
	1 & i
\end{pmatrix}.
\end{gather}

\begin{Lem}\label{l16}
For $\Hh \in \Hbbs$  and  $A\in\sltr$ define $\Hh_{A} = A^*\Hh A$. Then we have
\begin{itemize}
	\item[$(a)$] $M_A = A^{-1}MA$, where $M_A$ and $M$ are solutions of \eqref{cs} with Hamiltonians $\Hh_A$ and $\Hh$\textup{;}
	\item[$(b)$] $m_{\Hh_A} \doteq \sigma_1 A^{*}\sigma_1 m_\Hh$, where $m_{\Hh_A}$ and $m_{\Hh}$ are the Weyl functions corresponding to $\Hh_A$ and $\Hh$.
\end{itemize}
\end{Lem}
\beginpf The proof is a calculation. We have $J^{-1} = -J$ and
\begin{gather*}
A = 
\begin{pmatrix}
	a_{11} & a_{12}
	\\
	a_{21} & a_{22}
\end{pmatrix},
\quad 
A^{-1} = 
\begin{pmatrix}
	a_{22} & -a_{12}
	\\
	-a_{21} & a_{11}
\end{pmatrix},
\quad 
JA^{-1}J = 
\begin{pmatrix}
	-a_{11} & -a_{21}
	\\
	-a_{12} & -a_{22}
\end{pmatrix}
=-A^*.
\end{gather*}
Therefore we can write
\begin{gather*}
J (A^{-1}MA)' = J A^{-1}M'A = (J A^{-1} J)( J^{-1}M'A) = A^*(JM')A = zA^*HMA = z(A^*HA)A^{-1}MA,
\end{gather*}
which proves assertion $(a)$ of the lemma. To prove $(b)$, take $\omega\in \R$. Formula \eqref{wfbis} for $m_{\Hh_A}$ takes the form
$$m_{\Hh_A} \doteq \lim_{x \to \infty}\sigma_1 M_A(x,z)^{T} \sigma_1 \omega \doteq \lim_{x \to \infty}\sigma_1 (A^{-1}M(x,z)A)^{T} \sigma_1 \omega \doteq \lim_{x \to \infty}\sigma_1 A^{*}M(x,z)^{T}(A^{-1})^{*} \sigma_1\omega.
$$
Note that  $\sigma_1(A^{-1})^{*} \sigma_1 \in \sltr$, hence the number $\tilde \omega \doteq \sigma_1(A^{-1})^{*} \sigma_1\omega$ is in $\R$. Then, we have 
$$m_{\Hh_A} \doteq \lim_{x \to \infty}\sigma_1 A^{*}M(x,z)^{T}\sigma_1\tilde\omega \doteq
\sigma_1 A^{*}\sigma_1\lim_{x \to \infty}\sigma_1 M(x,z)^{T}\sigma_1\tilde\omega
\doteq \sigma_1 A^{*}\sigma_1 m_\Hh.$$ 
This proves the lemma. \qed

\medskip

Given a Hamiltonian $\Hh \in \Hbbs$ and $\ell>0$, define the Hamiltonian $\Hh_\ell = \Hh(\ell + \cdot)$ on $\R_+$. Let $m_{\Hh_\ell}$ be the Weyl function of $\Hh_\ell$.  For the following lemma see, e.g., the proof of formula $(2.13)$ in \cite{BD20}.
\begin{Lem}\label{l15}
We have 	$m_{\Hh}(z) \doteq \sigma_1 M(\ell,z)^{T} \sigma_1m_{\Hh_\ell}(z)$ for all $z \in \C_+$, where $M$ is the solution of~\eqref{cs}.
\end{Lem}
\begin{Lem}\label{l110}
Let $H$ be a singular Hamiltonian on $\R_+$. Define
\begin{gather*}
	\Hh(x) = 
	\begin{cases}
		\idm, & x \in [0, \ell),\\
		H(x-\ell), & x \ge \ell,
	\end{cases}
	\qquad 
	E_{lz} = 
	\begin{pmatrix}
		\cos \ell z & \sin \ell z\\
		-\sin \ell z & \cos \ell z
	\end{pmatrix}.
\end{gather*}
Then we have $m_{\Hh} \doteq  E_{lz}m_H$ and $f_{\Hh} =  e^{2ilz}f_H$.
\end{Lem}
\beginpf Notice that $M(t,x) = E_{xz}$ solves \eqref{csv} with the Hamiltonian $\Hh$ for $x\le \ell$. Lemma  \ref{l15} then gives
\begin{align}
m_{\Hh} 
&\doteq \sigma_1 \begin{pmatrix}
	\cos \ell z & \sin \ell z\\
	-\sin \ell z & \cos \ell z\\
\end{pmatrix}^{T} \sigma_1 m_H \doteq \begin{pmatrix}
	\cos \ell z & \sin \ell z\\
	-\sin \ell z & \cos \ell z\\
\end{pmatrix} 
m_H
\end{align}
as claimed. Applying relation \eqref{eq: LR matrices} we get
\begin{align*}
f_{\Hh} \doteq R m_\Hh &\doteq RE_{\ell z}m_H = RE_{\ell z}Lf_H
\\ 
&\doteq
\begin{pmatrix}
	1 & -i\\
	1 & i
\end{pmatrix}
\begin{pmatrix}
	\cos \ell z & \sin \ell z\\
	-\sin \ell z & \cos \ell z
\end{pmatrix}
\begin{pmatrix}
	i & i\\
	-1 & 1
\end{pmatrix}
f_H \doteq     \begin{pmatrix}
	e^{2ilz} & 0\\
	0 & e^{2ilz}
\end{pmatrix}f_H.
\end{align*}
The proof of the lemma is concluded.
\qed

\medskip

\subsection{Dirac operators with measures}\label{section: measure-valued potentials} Let us denote by $\cM$ the set of all signed complex Borel measures on $\R$ with $\supp\mu \subset \R_+$ such that the total variation of $\mu$ is finite on all intervals $[0, L]$, $L \ge 0$. Each element $q \in \cM$ generates the matrix potential
\begin{equation}\label{eq45}
Q = 
\begin{pmatrix}
\mu_2 & \mu_1 \\
\mu_1 & -\mu_2 
\end{pmatrix}, 
\qquad 
\begin{cases}
\mu_1 = \Re q,\\
\mu_2 = \Im q.
\end{cases}
\end{equation}
Conversely, each matrix potential $Q$ of the form \eqref{eq45} with real entries $\mu_{1,2} \in \cM$ generates a complex-valued measure $q = \mu_1 + i\mu_2$ in $\cM$ that we also refer to as the potential. It is possible to associate a self-adjoint Dirac operator $\Di_{q}$ to each $q \in \cM$. That was done by Zeng in paper \cite{Zeng23} devoted to Dirac operators with measures. In particular, Zeng \cite{Zeng23} proved relation of these operators to canonical Hamiltonian systems, that we will use in the proof of Theorem \ref{t2}. For the reader's convenience, we give a summary of some results from \cite{Zeng23}. Since our measures $q \in \cM$ might have a nontrivial point mass at $x =0$ (the case formally excluded by Zeng, though he mentioned that it can be covered as well), we make necessary modifications of definitions.

\medskip

Take $q \in \cM$, construct the corresponding matrix potential $Q$ as in \eqref{eq45}, and define the solution of the differential equation 
\begin{equation}\label{eq26}
JN_q'(x) + Q(x)N_q(x) = 0, \qquad \lim_{x \to 0, \; x< 0}N_q(x) = \IDM,
\end{equation}
as the unique $2\times 2$ continuous from the right matrix function $N_{q}$ such that 
\begin{equation}\label{eq27}
JN_q(x) + \int_{(-\infty, x]}\mathbf{g}(Q(\{x_1\})J^*)Q(x_1)N_q(x_1)= J,
\qquad x \in \R,
\end{equation}
where 
$$
\mathbf{g}(T) = \IDM +  \sum_{n=2}^{\infty}\frac{T^{n-1}}{n!},
$$
for a $2\times 2$ matrix $T$. This definition of the solution $N_{q}$ of \eqref{eq26} via \eqref{eq27}, as well as the existence and uniqueness of such a solution are due to Persson \cite{Persson} (see Theorem 3.1 in \cite{Persson}). The reader might note appearance of $J^*$ in the expression $\mathbf{g}(Q(\{x_1\})J^*)$ in \eqref{eq27}. This multiplicative factor $J^*$ does not appear in the work \cite{Persson} of Persson. The explanation is simple: \eqref{eq26} is in fact the differential equation for $JN_q$ (namely, $JN_q' + QJ^* \cdot JN_q = 0$), and the coefficient in front of $JN_q$ in \eqref{eq26} is $QJ^*$, not $Q$. Then, we need to use same coefficient in \eqref{eq27}. Finally, $\mathbf{g}(Q(\{x_1\})J^*)Q(x_1)N_q(x_1)$ is just the short way to write $\mathbf{g}(Q(\{x_1\})J^*)Q(x_1)J^* \cdot JN_q(x_1)$.

\medskip

For regular potentials,  the matrix-valued function $N_{q}$ solving \eqref{eq27} coincides with the classical solution of \eqref{eq26}. More precisely, if $q = s\,dx$ for some $s\in L^1_{\loc}(\R_+)$, and we define $S$ by $S\,dx = J Q$, then $g(Q(\{x_1\})J^*) = \idm$ for every $x_1\in\R_+$, and \eqref{eq26}, \eqref{eq27} are both equivalent  to the integral equation
$$
N_q(x) = \IDM + \int_{0}^{x}S(x_1)N_q(x_1)\,dx_1, \;\; x \ge 0, \qquad N_q(x) = \IDM \mbox{on } (-\infty, 0),
$$
that can be solved by iterations:
\begin{align}
\label{eq: series representation}
N_q(x) = \IDM + 
\int_{0}^{x} S(x_1)\,dx_1 &+
\int_{0}^{x}S(x_1)\int_{u}^{x_1} S(x_2)\,dx_2\,dx_1
\\
\nonumber
&+ \int_{0}^{x}S(x_1)\int_{0}^{x_1} S(x_2)\int_{0}^{x_2}S(x_3)\,dx_3\,dx_2\,dx_1 + \ldots.
\end{align}
Here, the series converges in the matrix norm and the $n$-th term can be bounded by
\begin{gather}
\label{eq: series bound}
\left\|\int_0^x\ldots\int_0^{x_{n - 1}}S(x_1)S(x_2) \cdots S(x_n)\,dx_n \ldots dx_1\right\|\le \frac{1}{n!}\left(\int_0^x\|S(x_1)\|\,dx_1\right)^n\le \frac{(4\|s\|_{L^1[0,x]})^n}{n!}.
\end{gather}
For a general $q \in \cM$, the solution $N_{q}$ can be approximated by solutions of regular equations. More precisely, if $\phi$ is a nonnegative continuous function supported on $[-1, 1]$, $\|\phi\|_{L^1(\R)} = 1$, and  $\phi_\eps = \eps^{-1}\phi(x/\eps)$, $\eps> 0$, is the corresponding approximate unity, then the solutions of regularized equations 
\begin{equation}\label{eq30}
JN_{q^{(\eps)}}'(x) + Q^{(\eps)}(x)N_{q^{(\eps)}}(x) = 0, \quad N_{q^{(\eps)}}(-1) = \IDM, 
\quad  q^{(\eps)}(x) = \int_{\R}\phi_\eps(y-x)q(y),
\end{equation}
converge pointwise on $\R$ to $N_{q}$, see \cite{Persson} (we use initial condition at  $x= -1$ instead of $x=0$ because $\supp q^{(\eps)} \subset [-1, \infty)$ for small $\eps > 0$, one can also use any other point in $(-\infty, 0)$ for the initial condition). In particular, for every $q \in \cM$ we have $\det N_{q} = 1$ on $\R$ and the multiplicative property 
 \begin{equation}\label{eq29}
N_q(w_2) = N_q(w_2, w_1)N_q(w_1), 
\qquad -\infty < w_1 < w_2 < \infty,
\end{equation}
holds, where $N_{q}(w_2, w_1)$ is the value at $w_2 - w_1$ of the solution of \eqref{eq26} for the potential $x \mapsto \chi_{\R_+}(x)q(x + w_1)$, $\chi_{\R_+}$ being the indicator function of $\R_+$. Indeed, it is enough to use the multiplicative property for regularized solutions $N_{q^{(\eps)}}$ and take the limit. 

\medskip

For regular potentials $q$ (i.e., for $q = s\,dx$ with some $s\in L^1_{\loc}(\R_+)$) the corresponding Dirac operators are related to canonical systems as follows. One need to take the solution of \eqref{eq26} and define the Hamiltonian 
\begin{equation}\label{eq28}
\Hh_{q}(x) = N_q^*(x)N_q(x), \qquad x \ge 0.
\end{equation}
Then, the Dirac operator $\Di_{q}$ and the operator of canonical system $\Di_{\Hh_{q}}$ are unitary equivalent, see, e.g., \cite{Romanov} or Section 2.4 in \cite{B2020}. As we will see in a moment, the same relation holds for $q \in \cM$.

\medskip

At first, we need some notation. By ${\rm BV}_{\loc}^{r}(\R_+)$ we will denote the set of all continuous from the right functions $X: \R \to \C^2$ such that their coordinate functions have a finite variation on any interval $[0, L]$, $L> 0$, and the restriction of $X$ to the set $(-\infty, 0)$ is a constant vector in $\C^2$. This constant vector will be denoted by $X(0-)$. This agrees with the standard notation
$$
X(0-) = \lim_{x \to 0,\; x < 0} X(x),
$$
because much stronger property $X(x) = X(0-)$ is assumed for all $x < 0$ if $X \in {\rm BV}_{\loc}^{r}(\R_+)$. Recall that for $\alpha \in [0, 2\pi)$ we denoted by $e_{\alpha}$ the vector 
$$
e_{\alpha} = 
\begin{pmatrix}
\sin \alpha \\ \cos \alpha
\end{pmatrix}.
$$
For $q \in \cM$ and $X \in  {\rm BV}_{\loc}^{r}(\R_+)$, we will say that $\Di_q X \in L^2(\R_+, \C^2)$ if there exists a function $Y \in L^2(\R_+, \C^2)$ solving 
equation 
$$
JX'(x) + QX(x) = Y(x), \qquad x \ge 0,
$$
in the sense of Persson, i.e., such that for all $x \in \R$ we have
\begin{equation}\label{eq111}
JX(x) + \int_{(-\infty, x]}g(Q(\{x_1\})J^*)Q(x_1)X(x_1) = JX(0-) + \int_{(-\infty, x]}Y(x_1)\,dx_1,
\end{equation}
where we extended $Y$ by zero to $(-\infty, 0)$. Here, $Q$ is defined by \eqref{eq45}, in particular, $QX$ is a vector-valued measure.  Let us define the domain of the Dirac $\Di_{q}$ on $\R_+$ corresponding to the boundary condition $\alpha \in [0, 2\pi)$ by 
$$
\dom \Di_q = \Bigl\{X \!\in\!  {\rm BV}_{\loc}^{r}(\R_+) \cap L^2(\R_+, \C^2)\!:  \left\langle X(0-), e_{\alpha}\right\rangle_{\C^2} = 0,\; \Di_q X \!\in\! L^2(\R_+, \C^2) \mbox{ in\,the\,sense\,of\,\eqref{eq111}}
\Bigr\}.
$$
For $X \in \dom \Di_q$, we define 
$$
\Di_q: X \mapsto Y
$$ 
for the unique (up to values on sets of Lebesgue measure zero) function $Y \in L^2(\R_+, \C^2)$ in \eqref{eq111}. With this definition, $\Di_{q}$ sends functions on $\R$ into functions in $L^2(\R_+, \C^2)$. Since we are interested in $D_{q}$ as a densely defined self-adjoint operator on $L^2(\R_+, \C^2)$, an additional step is needed to place $\dom \Di_q$ into $L^2(\R_+, \C^2)$. For this, we note that the values of $X \in \dom \Di_q$ on $\R$ as well as the the function $Y = \Di_q X$ in \eqref{eq111} depend solely on the restriction of $X$ to $\R_+$. Indeed, we only need to check that for $X \in \dom D_q$ the value $X(0-)$ is determined by the restriction of $X$ to $\R_+$. For this, we substitute $x =0$ into \eqref{eq111} and get
\begin{align*}
JX(0-) - JX(0) 
&= g(Q(\{0\})J^*)Q(\{0\})X(0) \\
&= \left[\IDM +  \sum_{n=2}^{\infty}\frac{(Q(\{0\})J^*)^{n-1}}{n!}\right]Q(\{0\})X(0)\\
&= J\left[\sum_{n=1}^{\infty}\frac{(J^*Q(\{0\}))^{n}}{n!}\right]X(0) = Je^{J^* Q(\{0\})}X(0) - JX(0).
\end{align*}
Thus, the value 
\begin{equation}\label{eq1130}
X(0-) = e^{J^* Q(\{0\})}X(0)
\end{equation} 
is determined by the restriction of $X$ onto $\R_+$ and we can consider $\dom \Di_{q}$ as the subset of $L^2(\R_+, \C^2)$. A similar argument applies to any point $x \in \R$. In fact, we have
\begin{equation}\label{eq113}
X(x-) = e^{J^* Q(\{x\})}X(x), \qquad x \in \R,
\end{equation} 
while $Q(\{x\}) \neq 0$ for at most countable set of points $x \in \R$.   Similar argument applied to \eqref{eq27}  gives 
\begin{equation}\label{eq1131}
		N_q(x-) = e^{J^* Q(\{x\})}N_q(x), \qquad x \in \R.
\end{equation}
Let us now recall the definition of the operator of canonical system. We will deal with Hamiltonians $\Hh$ of {\it rank two} almost everywhere on $\R_+$, for the general case see \cite{Romanov} or \cite{Remlingb} (the latter book considers linear relations instead of operators to cover the most general situation). At first, we denote by ${\rm AC}_{\loc}(\R_+)$ the set of all functions $X: \R_+ \to \C^2$ such that their coordinate functions are absolutely continuous on any compact subset of $\R_+$. The Hilbert space $L^2(\Hh)$ is defined by
$$
L^2(\Hh) = \Bigl\{X: \R_+ \to \C^2: \;\; \|X\|^{2}_{L^2(\Hh)} = \int_{\R_+}\left\langle\Hh(x)X(x), X(x)\right\rangle_{\C^2} < \infty \Bigr\}.
$$
For $\tilde X \in  {\rm AC}_{\loc}(\R_+)$, we will say that $\Di_{\Hh} \tilde X \in L^2(\Hh)$ if there exists a function $\tilde Y \in L^2(\Hh)$ solving the equation
$$
J\tilde X'(x) = \Hh(x)\tilde Y(x), \qquad x\ge 0, 
$$
or, in other words, such that
\begin{equation}\label{eq112}
J\tilde X(x) = J\tilde X(0) + \int_{0}^{x}\Hh \tilde Y(x_1)\,dx_1, \qquad x \in \R_+.
\end{equation}
The domain of the canonical system operator $\Di_{\Hh}$ on $\R_+$ corresponding to the boundary condition $\alpha \in [0, 2\pi)$ is defined by  
$$
\dom \Di_{\Hh} = \Bigl\{\tilde X \in  {\rm AC}_{\loc}(\R_+) \cap L^2(\Hh): \bigl\langle \tilde X(0), e_{\alpha}\bigr\rangle_{\C^2}  = 0,\; \Di_{\Hh} \tilde X \in L^2(\Hh) \mbox{ in the sense of \eqref{eq112}}
\Bigr\}.
$$
For $\tilde X \in \dom \Di_{\Hh}$, we define 
$$
\Di_{\Hh}: \tilde X \mapsto \tilde Y
$$ 
for the unique (up to values on sets of Lebesgue measure zero) function $\tilde Y \in L^2(\Hh)$ in \eqref{eq111}. 

\medskip

\begin{Prop}\label{P5}
	Let $q \in \cM$, define $\Hh_q$ by \eqref{eq28}. 
	Then the Dirac operator $\Di_{q}$ in $L^2(\R_+, \C^2)$ defined on the domain $\dom \Di_{q}$ for the boundary condition $\alpha \in [0, 2\pi)$ is unitarily equivalent to the operator $\Di_{\Hh_q}$ in $L^2(\Hh_q)$ defined on the domain $\dom \Di_{\Hh_q}$ for the same boundary condition $\alpha$. In particular,  $\Di_{q}$ is a densely defined self-adjoint operator on $L^2(\R_+, \C^2)$. The unitary equivalence is given by the multiplication operator $U: X \mapsto N_{q}^{-1}X$ from $L^2(\R_+, \C^2)$ to $L^2(\Hh_q)$.
\end{Prop}
\beginpf Recall that $\Hh_q = N_q^* N_q$. In particular, we have  $\det\Hh_q = \det (N_q)^2 = 1$ everywhere on~$\R_+$. The same relation $\Hh_q = N_q^* N_q$ tells us that $X \in L^2(\R_+, \C^2)$ if and only if $UX \in L^2(\Hh)$, and, moreover, $\|X\|_{L^2(\R_+, \C^2)} = \|UX\|_{L^2(\Hh)}$. Thus, the multiplication operator $U: X \mapsto N_{q}^{-1}X$ is unitary from $L^2(\R_+, \C^2)$ to $L^2(\Hh_q)$. For functions $X \in {\rm BV}^{r}_{\loc}(\R_+)$, $\tilde X \in {\rm AC}_{\loc}(\R_+)$, related by $\tilde X = N_{q}^{-1}X$ we have $\bigl\langle \tilde X(0), e_{\alpha}\bigr\rangle_{\C^2} = 0$ if and only if $\bigl\langle N_{q}^{-1}(0) X(0), e_{\alpha}\bigr\rangle_{\C^2} = 0$, i.e., $\bigl\langle N_{q}^{-1}(0) e^{-J^*Q(\{0\})} X(0-), e_{\alpha}\bigr\rangle_{\C^2} = 0$, where $X(0-)$ is defined by \eqref{eq1130}. However, \eqref{eq1131} gives $N_{q}(0) = e^{-J^*Q(\{0\})}N_{q}(0-) = e^{-J^*Q(\{0\})}$. Thus, $N_{q}^{-1}(0) e^{-J^*Q(\{0\})}= I$ 
and we have $\bigl\langle \tilde X(0), e_{\alpha}\bigr\rangle_{\C^2} = 0$ if and only if $\bigl\langle X(0-), e_{\alpha}\bigr\rangle_{\C^2} = 0$. Now, the result is a consequence of the following two lemmas. \qed \qed
\begin{Lem} \label{Zeng1}
If a function $X \in {\rm BV}^{r}_{\loc}(\R_+)$ solves \eqref{eq111} for some $Y \in L^2(\R_+, \C^2)$, then $\tilde X = UX$ solves \eqref{eq112} with $\tilde Y = UY$. In particular, $\tilde Y \in L^2(\Hh)$ and $\tilde X \in {\rm AC}_{\loc}(\R_+)$. 
\end{Lem}
\begin{Lem} \label{Zeng2}
	If a function $\tilde X \in {\rm AC}_{\loc}(\R_+)$ solves \eqref{eq112} for some $\tilde Y \in L^2(\Hh)$, then $X = U^{-1} \tilde X$ solves \eqref{eq111} with $Y = U^{-1}\tilde Y$. In particular, $Y \in L^2(\R_+, \C^2)$ and $X \in {\rm BV}^{r}_{\loc}(\R_+)$.
\end{Lem}
\noindent{\bf Proof of Lemma \ref{Zeng1}.} The result is Claim 3.3 in \cite{Zeng23}. We repeat the proof for the reader's convenience. Integral equation \eqref{eq27} shows that $N_q$ is of locally bounded variation, i.e., its columns belong to the space ${\rm BV}^{r}_{\loc}(\R_+)$. Hence $\tilde X\in {\rm BV}^{r}_{\loc}(\R_+)$ as a product of ${\rm BV}^{r}_{\loc}(\R_+)$ functions. Moreover, \eqref{eq113} and \eqref{eq1131} imply
\begin{gather*}
	\tilde X(x-) = N_q(x-)^{-1}X(x-) = N_q(x)^{-1}e^{-J^* Q(\{x\})}e^{J^* Q(\{x\})}X(x) = N_q(x)^{-1}X(x) = \tilde X(x),
\end{gather*}
i.e., $\tilde X$ is continuous. For simplicity, we write below $\mathbf{G}(x_1)$ instead of $\textbf{g}(Q(\{x_1\})J^*)$. Equality \eqref{eq111} rewrites as
\begin{equation}\label{eq1111}
JN_q(x) \tilde X(x) + \int_{(-\infty, x]}\mathbf{G}(x_1)Q(x_1)N_q(x_1)\tilde X(x_1) = J\tilde X(0) + \int_{(-\infty, x]}Y(x_1)\,dx_1.
\end{equation}
Denote $M(x) = \int_{(-\infty, x]}\textbf{G}(x_1)Q(x_1)N_q(x_1)$.  For $2\times 2$ matrix-valued functions $\Theta_{1}$, $\Theta_{2}$ whose columns belong to ${\rm BV}^{r}_{\loc}(\R_+)$, the formula of integration by parts (see assterion $(iv)$ in Theorem (21.67) in \cite{Hewitt}) reads as
$$
\int_{(-\infty, x]}(d\Theta_{1}(x_1)) \Theta_{2}(x_1-) = \Theta_1(x)\Theta_2(x) - \Theta_1(-\infty)\Theta_2(-\infty) - \int_{(-\infty, x]}\Theta_{1}(x_1)\,d\Theta_{2}(x_1),
$$
where $d\Theta_{1}$, $d\Theta_{2}$ denote the representing measures of $\Theta_{1}$, $\Theta_{2}$ as functions of bounded variations. Clearly, if $\Theta_2$ is continuous, we can replace $\Theta_{2}(x_1-)$ by $\Theta_{2}(x_1)$. Applying this formula to $\Theta_{1} = M$ and to the continuous matrix-valued function $\Theta_{2} = (\tilde X, 0)$ with columns $\tilde X$, $\left(\begin{smallmatrix} 0\\0 \end{smallmatrix}\right)$, and considering the first column in the resulting expressions, we get
$$
	\int_{(-\infty, x]}\mathbf{G}(x_1)Q(x_1)N_q(x_1)\tilde X(x_1) 
	= M(x)\tilde X(x) - \int_{(-\infty, x]}M(x_1)\,d\tilde X(x_1).
$$
Equality \eqref{eq27} gives $M(x) = J - JN_q(x)$, $x \in \R$, hence
\begin{align*}
\int_{(-\infty, x]}\mathbf{G}(x_1)Q(x_1)N_q(x_1)\tilde X(x_1) 
&= (J - JN_q(x))\tilde X(x) - \int_{(-\infty, x]}(J - JN_q(x_1))\,d\tilde X(x_1)	\\
&= \int_{(-\infty, x]}JN_q(x_1)\,d\tilde X(x_1) - JN_q(x)\tilde X(x) + J\tilde X(-\infty).
\end{align*}
Combining this with \eqref{eq1111}, we get (recall that $X$ is constant and $N_q = I$ on $(-\infty, 0)$, so $\tilde X(-\infty) = \tilde X(0-) = \tilde X(0)$ by continuity)
\begin{gather*}
\int_{(-\infty, x]}JN_q(x_1)\,d\tilde X(x_1) = \int_{(-\infty, x]}Y(x_1)\,dx_1.
\end{gather*}
We see that the measures $JN_q(x_1)\,d\tilde X(x_1)$ and $Y(x_1)\,dx_1$ coincide, hence $\tilde X$ is absolutely continuous and
\begin{align*}
J\tilde X(x) - J\tilde X(0) &= \int_0^x J\,d\tilde X(x_1) = \int_0^x JN_q(x_1)^{-1}J^{-1}JN_q(x_1)\,d\tilde X(x_1) \\
&= \int_0^x JN_q(x_1)^{-1}J^{-1} Y(x_1)\,dx_1 = \int_0^x JN_q(x_1)^{-1}J^{-1} N_q(x_1) \tilde Y(x_1)\,dx_1\\
&= \int_0^x \Hh(x_1) \tilde Y(x_1)\,dx_1,
\end{align*}
because $JA^{-1}J^{-1} = A^*$ for every real $2\times 2$ matrix $A$ with unit determinant. The lemma follows. \qed

\medskip

Knowing that operators $\Di_q$ and $\Di_{\Hh_q}$ are unitary equivalent, it is of no surprise that they have the same spectral measures and Weyl functions. 
Recall that the Weyl function of $\Di_{q}$ for $q \in L^2(\R_+)$ is defined by
$$
m_{q}(z) = \lim_{x \to +\infty} \frac{n_{22}(x,z)}{n_{21}(x,z)}, \qquad \Im z > 0, 
$$
where $N = \left(\begin{smallmatrix}
n_{11} & n_{12} \\
n_{21} & n_{22} 
\end{smallmatrix}\right)$ 
is the solution of the Dirac system $JN'(x,z) + Q(x)N(x, z) = zN(x, z)$, $N(0, z) = \idm$, $z \in \C$. Let us use the same definition for potentials $q \in \cM$, with the interpretation of solution $N$ in the sense of Persson: 
\begin{equation}\label{eq27bis}
JN(x, z) + \int_{(-\infty, x]}\mathbf{g}((Q(\{x_1\}) - z)J^*)(Q(x_1)-z)N(x_1, z)= J,
\qquad x \in \R.
\end{equation}
Note that $N_q(x) = N(x, 0)$, $x \in \R_+$, for the function $N_q$ solving \eqref{eq27}. A variant of Lemma \ref{Zeng1} for $N$ in place of $X$ implies that $N(x, z)  = N_q(x) M(x, z)$, where $M$ is the fundamental matrix solution for the Hamiltonian $\Hh_q = N_q^* N_q$. In particular, for $\omega=\infty$, $\tilde \omega(x) \doteq \sigma_1 N_q(x)^{T} \sigma_1 \omega$ and each $\Im z > 0$ we have
$$
m_{q}(z) = \lim_{x \to +\infty} \sigma_1 N(x,z)^{T} \sigma_1 \omega = \lim_{x \to +\infty} \sigma_1 M(x,z)^{T} \sigma_1 \tilde \omega(x) = m_{\Hh_q}(z), 
$$
because $\tilde \omega(x) \in  \R \cup \infty$ for every $x \in \R_+$ (according to the Weyl circles analysis, for every $\Hh\in \Hbbs$, $z \in \C_+$, the sets $\sigma_1 M(x,z)^{T} \sigma_1(\C_+ \cup \R \cup \infty)$ shrink to the singleton set $\{m_{\Hh}(z)\}$ as $x \to \infty$, see, e.g.,\ Section 8 in \cite{Romanov}).  Thus, for every $q \in \cM$ we have 
\begin{equation}\label{eq114}
m_q = m_{\Hh_q}, \qquad \Hh_q = N_q^* N_q,
\end{equation}
as in the standard theory for $q \in L^1_{\loc}(\R_+)$. The Schur function $f_q$ is then defined via \eqref{sf}. Clearly, we have $f_q = f_{\Hh_q}$. Our final remark concerning the general theory of Dirac operators with measures is the spectral theorem for these operators. For potentials $q \in L^1_{\loc}(\R_+)$, it can be found in Sections 7, 14 of \cite{Den06}.
\begin{Prop}\label{P8}
Let $q \in \cM$. Define the main spectral measure $\mu_q$ of the Dirac operator $\Di_{q}$ on $\R_+$ to be the representing measure of the harmonic function 
$\Im m_q = \Re(\tfrac{1+f_q}{1-f_q})$, i.e., by 
\begin{equation}\label{eq81}
\frac{1-|f_{q}(z)|^2}{|1 - f_{q}(z)|^2} = \frac{1}{\pi}\int_{\R} \frac{\Im z}{|x - z|^2}\, d\mu_{q}(x), \qquad z \in \C_+,   
\end{equation}
where $f_q$ is the Schur function of $\Di_q$. Then the operator $\Di_q$ on $L^2(\R_+, \C^2)$ is unitarily equivalent to the multiplication operator by the independent variable in $L^2(\mu_q)$. The operator $V: L^2(\R_+, \C^2) \to L^2(\mu_q)$ of unitary equivalence is densely defined by
$$
V: X \mapsto \frac{1}{\sqrt{\pi}}\int_{\R_+}\bigl\langle X(x), N_1(x, \ov{z}) \bigr\rangle_{\C^2}\,dx, \qquad N_1(x, z) = N(x, z)\oz,
$$ 
on smooth functions $X: \R_+\to \C$ with compact support.
\end{Prop}
\beginpf In view of \eqref{eq114}, the result is a direct consequence of Proposition \ref{P5} and the spectral theorem for operators of canonical systems (see Section 8 in \cite{Romanov}, or Chapter 3 in \cite{Remlingb}). \qed

\medskip

In what follows we will deal with a subclass of $\cM$ -- potentials $q$ supported on lattices. To simplify notation, we define 
$$
\cM_{\ell} = \{q \in \cM: \; \supp q \subset \ell\Z_+\}.
$$
Take $q\in \cM_{\ell}$ and	let $Q$ be the matrix potential associated with $q$ via \eqref{eq45}. Define also the constant $2\times 2$ matrices 
$Q[\ell k] = Q(\{\ell k\})$. We have
\begin{equation}\label{eq33}
N_{q}(x)=e^{JQ[\ell n]}e^{JQ[\ell(n-1)]} \cdot \ldots \cdot e^{JQ[0]}, \qquad n = \lfloor x/\ell\rfloor, \qquad x\ge 0.
\end{equation}
Here $\lfloor y \rfloor$ is the integer part of $y$ (the maximal integer $k$ such that $k \le y$).

\medskip		

\subsection{Auxiliary lemmas}
Recall that $\D$ denotes the open unit disk $\{z \in \C:\;|z|<1\}$. Take a number $\ell > 0$ and consider the conformal map  $e^{2i\ell z}: \Omega_{\C_+} \to \Omega_{\D}$ from  $\Omega_{\C_+} = \{z \in \C_+: |\Re z|<\pi/2\ell\}$ to  $\Omega_{\D} = \{z \in \D: z \notin (-1,0]\}$. We will denote by $\omega_\ell: \Omega_{\D} \to \Omega_{\C_+}$ the inverse conformal map.

\medskip
\begin{Lem}\label{l10}
Suppose that the Schur function $g\in S(\C_+)$  satisfies $g(z + \pi/\ell) = g(z)$ for some $\ell > 0$ and all $z \in \C_+$, i.e., $g$ is $\pi/\ell$--periodic. Then $G = g(\omega_\ell(\cdot))\colon \Omega_{\D}\to \C$  extends analytically to the whole open unit disk~$\D$. In particular, 
\begin{itemize}
	\item[$(a)$] there exists the finite limit $g(\infty) = \lim_{y \to +\infty}g(iy)$\textup{;}
	\item[$(b)$] if $g(\infty) = 0$, then $g = e^{2i\ell z}f$, where $f$ is also $\pi/\ell$--periodic and $f\in S(\C_+)$.
\end{itemize}
Thus, one step of the Schur's algorithm $e^{2i\ell z}g_1 = \frac{g - g(\infty)}{1 - \ov{g(\infty)}g}$ for a $\pi/\ell$--periodic Schur function $g$ in $\C_+$ corresponds to the one step of the Schur's algorithm $zG_1 = \frac{G - G(0)}{1 - \ov{G(0)}G}$ for its counterpart $G$ in $\D$.
\end{Lem}
\beginpf Since $g(z + \pi/\ell) = g(z)$ for all $z \in \C_+$, the function $G = g(\omega_\ell(\cdot))$ extends to a continuous function on $\D \setminus\{0\}$. Therefore, $G$ is analytic in $\D \setminus\{0\}$. The point $0$ is a removable singularity of $G$ because $G$ is bounded. 
From here we see  that there exists the limit 
$$
G(0) = \lim_{\eps\to 0}G(\eps) = \lim_{\eps\to 0}g(\omega_\ell(\eps)) = \lim_{y\to+\infty}g(iy) = g(\infty).
$$
To prove $(b)$, notice that $g(\infty) = 0$ implies $G(0) =0$. Then by Schwarz lemma we have $G(\lambda) = \lambda F(\lambda)$, $\lambda \in \D$, for some $F\in S(\D)$. In particular, $g(z) = \omega_\ell^{-1}(z)F(\omega_\ell^{-1}(z)) = e^{2 i \ell z} f(z)$ for $f\in S(\C_+)$.
\qed

\medskip

\begin{Lem}\label{l38}
Let $A = \left(\begin{smallmatrix}
	a & b\\b &-a
\end{smallmatrix}\right)$ be a zero-trace symmetric matrix. Then we have
\begin{gather*}
	\exp(A) = \cosh(\lambda)\IDM + \frac{\sinh(\lambda)}{\lambda}A, \quad \lambda = \sqrt{a^2 + b^2}.
\end{gather*}
\end{Lem}
\beginpf
We have $A^2 = \left(\begin{smallmatrix}
a & b\\b &-a
\end{smallmatrix}\right)\left(\begin{smallmatrix}
a & b\\b &-a
\end{smallmatrix}\right) = \lambda^2 \idm$. Therefore
\begin{gather*}
\exp(A) = \sum_{k= 0}^{\infty}\frac{A^k}{k!}  = \sum_{k= 0}^{\infty}\frac{A^{2k}}{(2k)!} + \sum_{k= 0}^{\infty}\frac{A^{2k + 1}}{(2k + 1)!} 
= \sum_{k= 0}^{\infty}\frac{\lambda^{2k}}{(2k)!}\IDM + \sum_{k= 0}^{\infty}\frac{\lambda^{2k}}{(2k + 1)!}A
\end{gather*}
The series in the right hand-side are Taylor series of $\cosh(\lambda)$ and $\sinh(\lambda)/\lambda$. The lemma follows.
\qed

\subsection{Proof of Theorem \ref{t2}}\label{s32}
Given a Hamiltonian $\Hh$ on $\R_+$, and $0\le u<v\le\infty$, let us denote by $M(v, u, z)$ the value at $v-u$ of the fundamental matrix solution for the Hamiltonian $\Hh(u+\cdot)$ on $\R_+$.

\begin{Lem}\label{pl1}
If $q \in \cM_{\ell}$ and $f_q$ is the Schur function of $\Di_q$, then $f_q(z + \pi/\ell)=f_q(z)$ for all $z\in \C_+$.
\end{Lem}
\beginpf Let $\Hh_q = N_q^* N_q$ for the matrix-valued function $N_q$ in \eqref{eq33}, and let $M_q$ be the fundamental solution of the corresponding canonical system \eqref{cs}. The Hamiltonian $\Hh_q$ is a constant rank two Hamiltonian on $[\ell k, \ell (k + 1))$ for every $k\in \Z_+$. By Lemma \ref{l16}, the function $M_k(z) = M_q(\ell (k + 1), \ell k, z)$ has the form $M_k(z) = A_k^{-1}E_{\ell z} A_k$ for $A_k = \sqrt{\Hh_q(\ell k)}$ and the function $E_{\ell z}$ defined in Lemma \ref{l110}. In particular, $M_k$ is $\pi/\ell$--antiperiodic for every $k\in \Z_+$, i.e., it satisfies $M_k(z + \pi/\ell) = -M_q(z)$ for all $z\in \C_+$.  For each $n$ the chain rule for solutions of ODE gives 
\begin{equation}\notag
M_q(\ell n, z)  = M_{n-1}(z) M_{n-2}(z) \cdots M_0(z).
\end{equation}
Therefore, $M_q(2\ell n, z)$ is $\pi/\ell$-periodic for every $n\in \Z_+$. Now periodicity of $m_{\Hh_q}$ follows from the definition \eqref{wf} of Weyl function:
$$
m_{\Hh_q}(z) \doteq \lim_{x \to \infty}\sigma_1 M_q(x,z)^{T} \sigma_1 \omega
\doteq
\lim\limits_{n \to \infty,\; n \in \Z}\sigma_1 M_q(2\ell n,z)^{T} \sigma_1 \omega, \qquad 
z \in \C_+.
$$
To conclude, recall that $f_q$ is related to $m_q = m_{\Hh_q}$ via \eqref{sf}. \qed 

\begin{Lem}\label{lem: induction base}
Let $q \in \cM_{\ell}$ for some $\ell > 0$, and let $f_q$ be the Schur function of $\Di_q$. Then $f_q(\infty)$ exists and satisfies
$$
f_q(\infty) = \varkappa^{-1}(q(\{0\})),
$$ 
where $\varkappa$ is defined in \eqref{eq44}. Moreover, the Schur function of $\Di_{q_\ell}$, $q_\ell = q(\ell + \cdot)$, satisfies 
\begin{gather}\label{eq251}
e^{2i \ell z}f_{q_{\ell}}(z) = \frac{f_{q}(z)-f_q(\infty)}{1-\ov{f_q(\infty)}f_{q}(z)}, \qquad z \in \C_+.
\end{gather}
\end{Lem}
\beginpf Recall that $f_q$, $f_{q_\ell}$ coincide with Schur functions of Hamiltonians $\Hh_q = N_q^*N_q$, $\Hh_{q_\ell} = N_{q_\ell}^*N_{q_\ell}$, where
\begin{align*}
N_{q}(x)
&= e^{JQ[\ell n]}e^{JQ[\ell(n-1)]} \cdots e^{JQ[0]}, \hspace{12,5mm} n = \lfloor x/\ell\rfloor, \qquad x\ge 0,\\
N_{q_\ell}(x)
&= e^{JQ[\ell (n + 1)]}e^{JQ[\ell(n-1)]} \cdots e^{JQ[1]}, \qquad n = \lfloor x/\ell\rfloor, \qquad x\ge 0.
\end{align*}
In particular, we have
\begin{gather*}
\Hh_q(x) = 
\begin{cases}
	N_q(0)^*N_q(0), & x \in [0, \ell),\\
	N_q(0)^*\Hh_{q_\ell}(x-\ell)N_q(0), & x \ge \ell.
\end{cases}
\end{gather*}
Lemma \ref{l16}, Lemma \ref{l110} and relation \eqref{eq: LR matrices} give
\begin{gather*}
m_q\doteq \sigma_1 N_q(0)^{*}\sigma_1 E_{lz}m_{q_\ell},\qquad f_q\doteq R\sigma_1 N_q(0)^{*}\sigma_1 E_{lz}Lf_{q_\ell}.
\end{gather*}
Denote $\zeta = q(\{0\})$. We have $JQ[0] = \left(\begin{smallmatrix}-\Re \zeta & \Im \zeta \\ \Im \zeta & \Re \zeta\end{smallmatrix}\right)$. Lemma \ref{l38} gives
\begin{gather*}
N_q(0)^* = N_q(0) = e^{JQ[0]} = \cosh(|\zeta|) \cdot \IDM + \frac{\sinh(|\zeta|)}{|\zeta|} \cdot J Q[0].
\end{gather*}
Furthermore, the straightforward calculation shows
\begin{align*}
RE_{\ell z}L 
&= 
\begin{pmatrix}
1 & -i\\
1 & i
\end{pmatrix}
\begin{pmatrix}
\cos \ell z & \sin \ell z\\
-\sin \ell z & \cos \ell z
\end{pmatrix}
\begin{pmatrix}
i & i\\
-1 & 1
\end{pmatrix}
= 2i
\begin{pmatrix}
e^{i\ell z} & 0\\
0 & e^{-i\ell z}
\end{pmatrix},
\\
R \sigma_1 J Q[0]\sigma_1 E_{\ell z}L 
&= 
\begin{pmatrix}
1 & -i\\
1 & i
\end{pmatrix}
\begin{pmatrix}
\Re \zeta & \Im \zeta\\ \Im \zeta & -\Re \zeta
\end{pmatrix}
\begin{pmatrix}
\cos \ell z & \sin \ell z\\
-\sin \ell z & \cos \ell z
\end{pmatrix}
\begin{pmatrix}
i & i\\
-1 & 1
\end{pmatrix},\\
&=2i
\begin{pmatrix}
0 & \ov{\zeta}e^{-i\ell z}\\
\zeta e^{i\ell z} & 0
\end{pmatrix}.
\end{align*}
When working with the relation ``$\doteq$'', the constant factor $2i$ can be cancelled, hence
\begin{align*}
f_{q} \doteq  
\begin{pmatrix}
	\cosh(|\zeta|)e^{i\ell z} & \frac{\sinh(|\zeta|)}{|\zeta|}\ov{\zeta}e^{-i\ell z}\\
	\frac{\sinh(|\zeta|)}{|\zeta|}\zeta e^{i\ell z} & \cosh(|\zeta|)e^{-i\ell z} 
\end{pmatrix}    
f_{q_\ell}
\doteq
\begin{pmatrix}
	e^{2i\ell z} & \frac{\tanh(|\zeta|)}{|\zeta|}\ov{\zeta}\\
	\frac{\tanh(|\zeta|)}{|\zeta|}\zeta e^{2i\ell z} & 1 
\end{pmatrix}    
f_{q_\ell}.
\end{align*}
Note that the function $\varkappa^{-1}: \zeta \mapsto \frac{\tanh|\zeta|}{|\zeta|}\ov{\zeta}$ from $\C$ to $\D$ is inverse to the function $\varkappa: s \mapsto \frac{\ov{s}}{2|s|}\log\frac{1+|s|}{1-|s|}$ from $\D$ to $\C$. This follows from the fact that $\tanh^{-1}|s| = \frac{1}{2} \log\frac{1+|s|}{1-|s|}$. Therefore the obtained formula for $f_q$ rewrites as 
\begin{gather*}
f_q = 
\begin{pmatrix}
	e^{2i\ell z} & \varkappa^{-1}(\zeta)\\
	\ov{\varkappa^{-1}(\zeta)}e^{2i\ell z} & 1 
\end{pmatrix}    
f_{q_\ell},\qquad
f_q = \frac{e^{2i\ell z}f_{q_\ell} + \varkappa^{-1}(\zeta)}{\ov{\varkappa^{-1}(\zeta)}e^{2i\ell z}f_{q_\ell} + 1},\qquad 
e^{2i \ell z}f_{q_{\ell}} = \frac{f_{q}-\varkappa^{-1}(\zeta)}{1-\ov{\varkappa^{-1}(\zeta)}f_{q}}.
\end{gather*}
It remains to prove that  $f_q(\infty)=\varkappa^{-1}(\zeta)$. For this we simply write
\begin{gather*}
\lim_{y\to\infty}(f_{q}(iy)-\varkappa^{-1}(\zeta)) = \lim_{y\to\infty} e^{2i \ell \cdot iy}f_{q_{\ell}}(iy)(1-\ov{\varkappa^{-1}(\zeta)}f_{q}(iy)) = 0,
\end{gather*}
where we used the fact that both $f_q$ and $f_{q_\ell}$ are bounded in $\C_+$.
\qed
\medskip

\noindent{\bf Proof of Theorem \ref{t2}.}
Let $q$ be a discrete measure on $\R_+$ such that $\supp q \subset \ell \Z_{+}$ for some $\ell > 0$. 
The Schur function $f_q: \C_+ \to \D$ of $\Di_{q}$ on $\R_+$ coincides with the Schur function of the Hamiltonian $\Hh_{q} = N_{q}^*N_{q}$ on $\R_+$, where $N_q$ is defined in \eqref{eq33}. By Lemma \ref{pl1}, $f_q$ is $\pi/\ell$-periodic in $\C_+$. For $\varkappa$ defined in \eqref{eq44}, Lemma \ref{lem: induction base} states that $q(\{0\}) = \varkappa(f_{q}(\infty))$ and $f_{q_\ell}$ is the first Schur's iterate of $f_{q}$. Then \eqref{eq39} follows by the induction principle.
Since there exists infinitely many Schur iterates of $f_q$ (equivalently, $|f_{q,k}(\infty)| < 1$ for each $k \in \Z_+$), we have $f_q \in S_{\ell, *}(\C_+)$. 
It remains to show that for every $f \in S_{\ell, *}(\C_+)$ there exists $q\in \cM_{\ell}$ such that $f = f_q$. Take $f \in S_{\ell, *}(\C_+)$, find its recurrence coefficients $\{f_k(\infty)\}_{k\in \Z_+}$, and define $q$ by
$$
q(\{k\ell\}) = \varkappa(f_{k}(\infty)), \qquad k\in \Z_+.
$$
The first part of the proof shows that recurrence coefficients of $f$ and $f_q$ coincide. The same is true for the corresponding Schur functions $F = f \circ\omega_\ell$, $F_q = f_q \circ\omega_\ell$ in $\D$ in Lemma \ref{l10} and the standard recurrence coefficients in Schur's algorithm in $\D$ for these functions. It follows that $F = F_q$, see formula $(1.3.43)$ in \cite{Simonbook1}. Thus, $f =f_q$ and the proof is completed. \qed

\subsection{Proof of Corollaries \ref{c3}-\ref{c5}}
Consider a Schur function $F \in S_*(\D)$. Note that $\frac{1+zF}{1-zF}$, $\frac{1+F}{1-F}$ are analytic functions in $\D$ with positive real part. It follows that there are finite nonnegative measures $\sigma$, $\nu$ on the unit circle $\T$ being the boundary values of the corresponding positive harmonic functions in $\D$, i.e.,
\begin{align}
\frac{1-|\lambda F(\lambda)|^2}{|1 - \lambda F(\lambda)|^2} 
&=  \int_{\T} \frac{1-|\lambda|^2}{|1 - \bar \lambda \xi|^2}\,d\sigma(\xi) \label{eq75},\\
\frac{1-|F(\lambda)|^2}{|1 - F(\lambda)|^2}
&=\int_{\T} \frac{1-|\lambda|^2}{|1 - \bar \lambda  \xi|^2}\,d\nu(\xi), \label{eq76}
\end{align}
for all $\lambda \in \D$, see Section I.3 in \cite{Garnett}. Taking $\lambda = 0$, one can see that $\sigma$ is a probability measure on~$\T$. In general, the measure $\nu$ is not probabilistic and there is no simple way of expressing $\nu/\nu(\T)$ in terms of $\sigma$. However, $\nu/\nu(\T)$ coincides with the measure $\sigma_{\alpha,1}$ generated via \eqref{eq75} by the first Schur iterate $F_{\alpha,1}(\lambda) = \frac{1}{\lambda}\frac{F_{\alpha}(\lambda)-F_{\alpha}(0)}{1 - \ov{F_{\alpha}(0)}F_{\alpha}(\lambda)}$ of the Schur function $F_{\alpha} = \alpha F$ for a specific value $\alpha \in \T$. More precisely, we have the following relation.
\begin{Lem}\label{lReduction}
Suppose that $F \in S(\D)$ and $\alpha \in \T$ are such that $\alpha = \frac{1-\ov{F(0)}}{1-F(0)}$, $|F(0)| < 1$. Then
\begin{equation}\label{eqReduction}
	\frac{1-|\lambda F_{1,\alpha}|^2}{|1 - \lambda F_{\alpha,1}|^2} = \frac{1-|F(0)|^2}{|1+\alpha F(0)|^2}\frac{1-|F|^2}{|1 - F|^2}, \qquad \lambda \in \D.
\end{equation}
In particular, we have $\nu/\nu(\T) = \sigma_{\alpha,1}$, i.e.,
\begin{equation}\label{eq76bis}
	\frac{1-|\lambda F_{\alpha,1}(\lambda)|^2}{|1 - \lambda F_{\alpha,1}(\lambda)|^2} 
	=  \frac{1}{\nu(\T)}\int_{\T} \frac{1-|\lambda|^2}{|1 - \bar \lambda \xi|^2}\,d\nu(\xi), \qquad \lambda \in \D.
\end{equation}
\end{Lem}
\beginpf We have $\lambda F_{\alpha,1} = \frac{F_{\alpha}-F_{\alpha}(0)}{1 - \ov{F_{\alpha}(0)}F_{\alpha}} = 
\frac{\alpha F-\alpha F(0)}{1 - \ov{F(0)}F}$. It follows that
$$
1 - |\lambda F_{\alpha,1}|^2 = \frac{(1-|F(0)|^2)(1-|F|^2)}{|1 - \ov{F(0)}F|^2},  
$$
and, taking into account $\ov{F(0)} + \alpha=1+\alpha F(0) $, we get
$$
1 - \lambda F_{\alpha,1} = \frac{1 - \ov{F(0)}F-\alpha F+\alpha F(0)}{1 - \ov{F(0)}F} = \frac{(1+\alpha F(0))(1- F)}{1 - \ov{F(0)}F}.
$$
These two relations yield \eqref{eqReduction}. Integrating \eqref{eqReduction} over $\T$, we see that $\frac{1-|F(0)|^2}{|1+\alpha F(0)|^2} = \frac{1}{\nu(\T)}$. Now, \eqref{eq76bis} follows from \eqref{eq76}. \qed

\medskip

Let a Schur function $F \in S_*(\D)$ and a probability measure $\sigma$ on $\T$ be related as in \eqref{eq75}. There is a strong connection between properties of orthogonal polynomials in $L^2(\sigma)$ and the properties of $F$. Let $\{\Phi_{k}\}_{k \ge 0}$ be the family of polynomials such that $\Phi_{k} = z^k + \ldots + \Phi_{k}(0)$ for each $k \ge 0$ and 
$$
(\Phi_{k_1}, \Phi_{k_2})_{L^2(\sigma)} = 0 \; \mbox{ for all } \; 0 \le k_1 < k_2 < \infty.
$$
Let also $\Phi_{k}^{*}(z) = z^{k}\ov{\Phi_{k}(1/\bar z)}$ denote the so-called reflected orthogonal polynomials. It is well-known (see, e.g., Theorem 1.5.2 in \cite{Simonbook1}) that 
\begin{equation}\label{rr}
\Phi_{k+1}(z) = z\Phi_{k}(z) - \ov{\alpha_k}\Phi^*_{k}(z), \qquad k \ge 0,
\end{equation}
for all $z \in \C$ and some sequence $\{\alpha_{k}\}_{k \ge 0} \subset \D$, which is called the sequence of recurrence coefficients of~$\sigma$. Moreover, any $\{\alpha_{k}\}_{k \ge 0} \subset \D$ arises uniquely as the sequence of recurrence coefficients of some probability measure $\sigma$ supported on an infinite subset of $\T$, see Section 1.7 in \cite{Simonbook1}. The following theorem is due to Geronimus \cite{Ger44}, its modern exposition can be found in Section 3.1 of \cite{Simonbook1}.
\begin{Thm}[Geronimus' theorem] \label{tG}
Assume that a Schur function $F \in S_*(\D)$ and a probability measure $\sigma$ on $\T$ are related by \eqref{eq75}. Let $\{F_k(0)\}_{k \ge 0}$ be the recurrence coefficients of $F$, see \eqref{sa}, and let $\{\alpha_{k}\}_{k \ge 0} \subset \D$ be the recurrence coefficients of $\sigma$ in \eqref{rr}. Then $F_{k}(0) = \alpha_k$ for all $k \ge 0$.
\end{Thm}
Geronimus' theorem allows to formulate many results of the theory of orthogonal polynomials in the language of Schur functions. Let us state several landmark theorems of the theory in this form. In all of them, $\sigma = w_\sigma\,dm_{\T} + \sigma_{s}$ (respectively, $\nu = w_\nu\,dm_{\T} + \nu_{s}$) is the Radon-Nikodym decomposition of $\sigma$ (respectively, $\nu$) into the absolutely continuous and singular parts. Here and below $m_{\T}$ denotes the Lebesgue measure on the unit circle $\T$ normalized so that $m_{\T}(\T) = 1$. For $F$, $\nu$ related by \eqref{eq76} it is known that 
\begin{equation}\label{eq77}
w_{\nu}(\xi) = \frac{1-|F(\xi)|^2}{|1 - F(\xi)|^2} \mbox{ for Lebesgue almost all } \xi \in \T,
\end{equation}
and 
\begin{equation}\label{eq78}
\nu = *\lim_{r \to 1}\frac{1-|F(r\xi)|^2}{|1 - F(r\xi)|^2}\,d m_{\T}(\xi),
\end{equation}
where $*\lim$ stands for the *-weak convergence of measures. For the proof of relations \eqref{eq77}, \eqref{eq78}, see Theorem I.5.3, Theorem I.3.1 in \cite{Garnett}. Similar relations hold for the measure $\sigma$ (for $\sigma$, one need to replace $F$ by $\xi F$ in \eqref{eq77}, \eqref{eq78}).  In the next three results we assume that a Schur function $F \in S_*(\D)$ and measures $\sigma$, $\nu$ on $\T$ are related by \eqref{eq75}, \eqref{eq76} and  $\{F_k(0)\}_{k \ge 0} = \{\alpha_{k}\}_{k \ge 0}$ are the recurrence coefficients of $F$, $\sigma$.
\begin{Thm}[Rakhmanov's theorem] \label{tR}
If $|F| < 1$ almost everywhere on $\T$, then $\lim_{k \to \infty}F_k(0) = 0$.
\end{Thm}
\beginpf By \eqref{eq77} if $|F| < 1$ almost everywhere on $\T$, then $w_{\sigma} > 0$ almost everywhere on $\T$. By the classical Rakhmanov's theorem \cite{R77}, see also Section 9.1 in \cite{Simonbook1}, this implies $\lim_{k \ge 0} \alpha_k = 0$. It remains to use Geronimus theorem. \qed 

\medskip

Below, $W^1(\T)$ stands for the Wiener algebra of all absolutely convergent Fourier series on $\T$.

\begin{Thm}[Baxter's theorem] \label{tB}
The following assertions are equivalent: 
\begin{itemize}
	\item[$(1)$] $F \in S_*(\D) \cap W^1(\T)$, $\max_{\xi \in \T} |F(\xi)| < 1$\textup{;}
	\item[$(2)$] $\sum_{k \ge 0} |F_k(0)| < \infty$\textup{;}
	\item[$(3)$] $\sigma = w_\sigma \, dm_{\T}$, where $w_\sigma \in W^1(\T)$, $\inf_{\T}w_\sigma > 0$\textup{;}
	\item[$(4)$] $\nu = w_\nu \, dm_{\T}$, where $w_\nu \in W^1(\T)$, $\inf_{\T}w_\nu > 0$.
\end{itemize}
\end{Thm}
\beginpf Equivalence $(2) \Leftrightarrow (3)$ is the original Baxter's theorem \cite{Baxter} (modulo Geronimus' theorem \ref{tG}). Equivalence $(1) \Leftrightarrow (2) \& (3)$ was proved by Golinskii in \cite{Gol97}. So, it suffices to prove the equivalence $(2) \Leftrightarrow (4)$. For this we take $\alpha \in \T$, set $F_{\alpha} = \alpha F$ and note that $(2)$ is equivalent to 
$\sum_{k \ge 0} |F_{\alpha, k}(0)| < \infty$ because $F_{\alpha, k}(0) = \alpha F_{k}(0)$ for all $k \ge 0$. Assumption $\sum_{k \ge 1} |F_{\alpha, k}(0)| < \infty$ is equivalent to the assumption $(2)$ for $F_{\alpha, 1}$, the first Schur iterate of $F_{\alpha}$. In turn, the latter is equivalent to the fact that the probability measure $\sigma_{\alpha,1}$ corresponding to $F_{\alpha, 1}$ satisfies assumption $(3)$ with $\sigma$ replaced by $\sigma_{\alpha, 1}$. Now, choosing $\alpha = \frac{1-\ov{F(0)}}{1-F(0)}$ we see from Lemma \ref{lReduction} that $\sigma_{\alpha, 1} = \nu/\nu(\T)$, and the result follows.  \qed 

\medskip

Next theorem deals with the Sobolev space $H_{1/2}(\T)$ -- the set of functions 
$$
H_{1/2}(\T) = \biggl\{g \in L^1(\T): \; \sum_{k \in \Z} |k||\hat g(k)|^2 < \infty\biggr\},
$$
where $\hat g(k)$ is the $k$-th Fourier coefficient of a function on the unit circle $\T$, see \eqref{eq80}.  
\begin{Thm}[Szeg\H{o}-Golinskii-Ibragimov theorem]\label{tSGI} 
The following assertions are equivalent: 
\begin{itemize}
	\item[$(1)$] $\sum_{k \ge 0} k|F_k(0)|^2 < \infty$\textup{;}
	\item[$(2)$] $\sigma = w_\sigma \, dm_{\T}$, where $w_{\sigma} = e^{\phi}$ for some $\phi \in H_{1/2}(\T)$\textup{;}
	\item[$(3)$] $\nu = w_\nu \, dm_{\T}$, where $w_{\nu} = e^{\psi}$ for some $\psi \in H_{1/2}(\T)$.
\end{itemize}
\end{Thm}
\beginpf The proof is similar to the proof of Theorem \ref{tB}. The equivalence $(1) \Leftrightarrow (2)$ is the original Szeg\H{o}-Golinskii-Ibragimov theorem \cite{GI71} (modulo Geronimus theorem \ref{tG}). Let us show that $(3)$ is equivalent to $(1)$. For this we take $\alpha \in \T$, set $F_{\alpha} = \alpha F$ and note that $(1)$ is equivalent to 
$\sum_{k \ge 0} k|F_{\alpha, k}(0)|^2 < \infty$ because $F_{\alpha, k}(0) = \alpha F_{k}(0)$ for all $k \ge 0$. Assumption $\sum_{k \ge 0} k|F_{\alpha, k}(0)|^2 < \infty$ is equivalent to the assumption $(1)$ for $F_{\alpha, 1}$, the first Schur iterate of $F_{\alpha}$. In turn, the latter is equivalent to the fact that the probability measure $\sigma_{\alpha,1}$ corresponding to $F_{\alpha, 1}$ satisfies assumption $(2)$ with $\sigma$ replaced by $\sigma_{\alpha, 1}$. Now, choosing $\alpha = \frac{1-\ov{F(0)}}{1-F(0)}$ we see from Lemma \ref{lReduction} that $\sigma_{\alpha, 1} = \nu/\nu(\T)$, and the result follows.  \qed 

\medskip

Let $q \in \cM$, and let $f_q$ is the Schur function of $\Di_q$. Recall that the main spectral measure $\mu_q$ of the Dirac operator $\Di_{q}$ on $\R_+$ is defined by 
\begin{equation}\label{eq81bis}
\frac{1-|f_{q}(z)|^2}{|1 - f_{q}(z)|^2} = \frac{1}{\pi}\int_{\R} \frac{\Im z}{|x - z|^2}\, d\mu_{q}(x), \qquad z \in \C_+.   
\end{equation}
Writing $\mu_q = w_{q}\,dx+\mu_{q,s}$ for $\mu_{q,s} \bot dx$, we obtain
\begin{equation}\label{eq85}
w_{\sigma}(x) = \frac{1-|f_q(x)|^2}{|1 - f_q(x)|^2} \; \mbox{ for Lebesgue almost all } \; x \in \R,
\end{equation}
and
\begin{equation}\label{eq83}
\mu_q = *\lim_{\eps \downarrow 0}\frac{1-|f_q(x + i\eps)|^2}{|1-f_q(x+i \eps)|^2}\,dx,
\end{equation}
where the limit is understood in the $*$-weak sense. These relations are analogous to \eqref{eq77}, \eqref{eq78} and have essentially the same proofs.

\medskip

Given $q \in \cM_{\ell}$, we know that $f_q \in S_{\ell, *}(\C_+)$, see Theorem \ref{t2}. By Lemma \ref{l10}, there exists $F_{q} \in S_*(\D)$  such that $F_{q}(e^{2i\ell z}) = f_{q}(z)$, $z \in \C_+$. Below we prove Corollaries \ref{c3}-\ref{c5} based on Theorem \ref{t2} and the following proposition. 

\begin{Prop}\label{p-reduction}
Let $q\in \cM_{\ell}$ and $f_q \in S_{\ell, *}(\C_+)$ be the measure-valued potential and  its Schur function. Also let $F_{q} \in S_*(\D)$ be such that $F_{q}(e^{2i\ell z}) = f_{q}(z)$, $z \in \C_+$. Then the spectral measure $\mu_q$ of $\Di_{q}$ is periodic with period $\pi/\ell$ and satisfies 
\begin{equation}\label{eq89}
	\frac{1-|F_{q}(\lambda)|^2}{|1 - F_{q}(\lambda)|^2} = \frac{\ell}{\pi}\int_{[-\frac{\pi}{2\ell},\frac{\pi}{2\ell})} \frac{1-|\lambda|^2}{|1 - \bar \lambda e^{2i\ell x}|^2}\, d\mu_{q}(x), \qquad  \lambda \in \D.
\end{equation}
In particular, if $\nu_{q}$ is the measure on $\T$ generated by $F_{q}$ via \eqref{eq76}, then 
$\frac{\ell}{\pi}\mu_{q}(E) = \nu_{q}(E^*)$, where $E^* = \{e^{2i\ell x}, \; x \in E\}$ for every Borel set $E \subset [-\frac{\pi}{2\ell},\frac{\pi}{2\ell})$. 
\end{Prop}
\beginpf Set $w_\eps(x) = \frac{1-|f_q(x + i\eps)|^2}{|1-f_q(x+i \eps)|^2}$, $x \in \R$, and take arbitrary $x_0 \in \R$. By \eqref{eq83}, for the proof of \eqref{eq89} it suffices to show that for every $\lambda \in \D$ we have
\begin{equation}\label{eq84}
\frac{1-|F_{q}(\lambda)|^2}{|1 - F_{q}(\lambda)|^2} = \lim_{\eps \to 0}\frac{\ell}{\pi}\int_{[x_0-\frac{\pi}{2\ell},x_0+\frac{\pi}{2\ell})} \frac{1-|\lambda|^2}{|1 - \bar \lambda e^{2i\ell x}|^2} w_{\eps}(x)\, dx.
\end{equation}
For $\eps> 0$, denote $r_\eps = e^{-2\ell \eps}$. In view of $F_{q}(e^{2i\ell z}) = f_{q}(z)$, $z \in \C_+$, relation \eqref{eq84} can be rewritten in the form  
\begin{align*}
\frac{1-|F_{q}(\lambda)|^2}{|1 - F_{q}(\lambda)|^2} 
&= \lim_{\eps \to 0}\frac{\ell}{\pi}\int_{[x_0-\frac{\pi}{2\ell},x_0+\frac{\pi}{2\ell})} \frac{1-|\lambda|^2}{|1 - \bar \lambda e^{2i\ell x}|^2} \frac{1-|F_q(e^{2i\ell (x + i\eps)})|^2}{|1-F_q(e^{2i\ell (x + i\eps)})|^2}\, dx\\
&= \lim_{\eps \to 0}\int_{\T} \frac{1-|\lambda|^2}{|1 - \bar \lambda \xi|^2} \frac{1-|F_q(r_\eps \xi)|^2}{|1- F_q(r_\eps \xi)|^2}\, dm_{\T}(\xi)\\
&= \lim_{\eps \to 0}\frac{1-|F_{q}(r_\eps \lambda)|^2}{|1 - F_{q}(r_\eps \lambda)|^2},
\end{align*}
which holds by continuity of $F_{q}$ in $\D$. Periodicity of $\mu_q$ follows from the fact that $x_0 \in \R$ in \eqref{eq84} is arbitrary, while the left hand side does not depend on $x_0$. Relation $\frac{\ell}{\pi}\mu_{q}(E) = \nu_{q}(E^*)$ is just the change of variables in the Lebesgue integral. \qed 

\medskip

In the proof of corollaries below we assume that $f_q$ and $F_q$ are related as in Proposition \ref{p-reduction}. In particular, their Schur's iterations satisfy $f_{q,k}(\infty) = F_{q,k}(0)$ for all $k \ge 0$.

\medskip

{\noindent \bf Proof of Corollary \ref{c3}.} Suppose that $\sigma_{ac}(\Di_{q}) = \R$. It follows that 
$|f_{q}| < 1$ almost everywhere on $\R_+$, see \eqref{eq85}. Then, the corresponding Schur function $F_{q}$ in $\D$ is such that $|F_q| < 1$ almost everywhere on $\T$. By Rakhmanov's theorem \ref{tR}, it follows that $F_{q,k}(0) \to 0$ as $k \to +\infty$. Then, Theorem \ref{t2} and the fact that $f_{q,k}(\infty) = F_{q,k}(0)$ imply $q(\{k\}) \to 0$ as $k \to +\infty$. \qed 

\medskip

{\noindent \bf Proof of Corollary \ref{c4}.} By Theorem \ref{t2}, $\sum_{k \ge 0}|q(\{k\})| < \infty$ if and only if $\sum_{k \ge 0} |F_{q, k}(0)| < \infty$. By Baxter's theorem \ref{tB}, this is equivalent to $\nu = w_\nu \, dm_{\T}$, $w_\nu \in W^1(\T)$, $\inf_{\T}w_\nu > 0$, for the measure $\nu$ constructed from $F_{q}$ by \eqref{eq76}. By Proposition \ref{p-reduction}, we have $\mu_q = w_{q}\,dx$, where $w_q(x) = w_{\nu}(e^{2ix})$ is a continuous  $\pi$-periodic function on $\R$ with $\inf_{x\in \R} |w_q (x)| = \inf_{\xi\in \T} |w_\nu(\xi)|  > 0$.  Moreover, for every $k\in \Z$ we have
\begin{gather*}
\frac{1}{\pi}\int_{-\pi/2}^{\pi/2}e^{-2ikx}w_q(x)\,dx = \int_{\T}\ov{\xi}^kw_{\nu}(\xi)\, dm_{\T}. 
\end{gather*}
It follows that assumption $\sum_{k \ge 0}|q(\{k\})| < \infty$ implies $\mu_q = w_{q}\,dx$ for a $\pi$-periodic function $w_q$ such that 
$\inf_{x\in \R} |w_q (x)| > 0$ and
\begin{gather*}
w_{q} = \sum_{k \in \Z} c_k e^{2i k x}, \qquad \sum_{k \in \Z} |c_k| = \sum_{k \in \Z} |\hat{w}_\nu(k)| = \|w_\nu\|_{W^1(\T)} < \infty.
\end{gather*}
The argument is reversible: any measure $\mu_q$ with these properties gives rise to the function $F_q$ such that $\sum_{k \ge 0} |F_{q, k}(0)| < \infty$, i.e., $\sum_{k \ge 0}|q(\{k\})| < \infty$ by Theorem \ref{t2}. \qed

\medskip 

{\noindent \bf Proof of Corollary \ref{c5}.} By Theorem \ref{t2}, assumption $\sum_{k \ge 0}k|q(\{k\})|^2 < \infty$ is equivalent to the assumption $\sum_{k \ge 0}k|F_{q, k}(0)|^2 < \infty$. Then, Szeg\H{o}-Golinskii-Ibragimov theorem \ref{tSGI} says that the last condition is equivalent to $\nu = e^{\psi} \,dm_{\T}$, $\psi \in H_{1/2}(\T)$ for the measure $\nu$ related to $F_{q}$ via \eqref{eq76}. By Proposition \ref{p-reduction}, this can be further equivalently reformulated in the following way: $\mu_q = w_q\,dx$ for a positive $\pi$-periodic function $w_q$ that satisfies $\log w_q(x) = \psi(e^{2ix})$ on $\R$. For these objects, we have
\begin{gather*}
\log w_q(x) = \sum_{k \in \Z} c_k e^{2i k x}, \qquad c_k = \frac{1}{\pi}\int_{-\pi/2}^{\pi/2}e^{-2ikx}\log w_q\,dx = \int_{\T}\ov{\xi}^k\psi(\xi)\, dm_{\T}.
\end{gather*}
In particular, the series $\sum_{k \in \Z} |k||c_k|^2 = \sum_{k \in \Z} |k||\hat{\psi}(k)|^2$ converge or not simultaneously. From here we see that $\sum_{k \ge 0}k|q(\{k\})|^2 < \infty$ if and only if $\mu_q = w_q\,dx$, where $\log w_q = \sum_{k \in \Z}c_k e^{2ikx}$ with $\sum_{k \in \Z}|k||c_k|^2 < \infty$. \qed

\medskip

\section{Spectral continuity estimates for potentials in \texorpdfstring{$L^2(\R_+)$}{L2(R+)}}\label{s4}
Our aim in this section is to prove Theorem \ref{t1}. The strategy of the proof is to approximate potentials $q, \tilde q \in L^2(\R_+)$ by discrete measures supported on $\ell \Z_+$, $\ell \to 0$. By Theorem \ref{t2}, solution of the direct and inverse spectral problems for measures supported on $\ell \Z_+$ relies on a Schur algorithm. This and the continuity estimates for the Schur algorithm that we prove below will be the principal ingredients of the proof of Theorem \ref{t1}.
\subsection{Estimates for Schur's algorithm}\label{s41}
Recall that the Schur's algorithm for functions $F \in S_*(\D)$ is defined by 
\begin{equation}\label{sa-bis}
F_0 = F, \qquad 
zF_{k+1} = \frac{F_k - F_k(0)}{1 - \ov{F_k(0)}F_k}, \qquad k \ge 0.
\end{equation}
It is known that the mapping $F \mapsto \{F_k(0)\}_{k \ge 0}$ is a bijection from $S_*(\D)$ onto the space $\ell^0(\Z_+, \D)$ consisting of all sequences $\{q(k)\}_{k \in \Z_+}$ indexed by nonnegative integers $\Z_+ = \{k \in \Z:\;k \ge 0\}$ such that $|q(k)| < 1$ for every $k \in \Z_+$. Moreover, one can introduce topologies on $S_*(\D)$, $\ell^0(\Z_+, \D)$ so that this mapping became a homeomorphism, see Section 1.3.6 in \cite{Simonbook1}. One version of the Szeg\H{o} theorem (see Section 2.7.8 in \cite{Simonbook1}) states that for every function $F \in S_*(\D)$ we have the following relation between recurrence coefficients $\{F_k(0)\}_{k \ge 0}$ and the logarithmic integral of $F$:
\begin{gather}\label{eq: szego theorem}
\eta(F):=\prod_{k = 0}^{\infty}(1 - |F_k(0)|^2) = \exp\left(\int_{\T}\log(1 - |F|^2)\,dm_{\T}\right).
\end{gather}
Since $|F_k(0)| < 1$ for each $k$ and $|F| \le 1$ almost everywhere on $\T$, the quantities $\prod_{k = 0}^{\infty}(1 - |F_k(0)|^2)$ and $\exp\left(\int_{\T}\log(1 - |F|^2)\,dm_{\T}\right)$ are defined for every $F \in S_*(\D)$ but could be, in general, equal to zero. We will denote by $\szc$ the class of all functions $F \in S_*(\D)$ such that $\eta(F) > 0$. Equivalently, $F \in \szc$ if $F \in S_*(\D)$ and any of the following two equivalent assertions holds:
\begin{itemize}
\item $\sum_{k \ge 0}|F_k(0)|^2 < \infty$, or
\item $\log(1-|F|^2) \in L^1(\T)$.
\end{itemize}
For a function $h \in L^1(\T)$ and $k \in \Z$, let us denote by $\hat h(k)$ its $k$-th Fourier coefficient,
\begin{equation}\label{eq80}
\hat h(k) = \int_{\T}h(\xi)\bar\xi^k\,dm_{\T}(\xi).    
\end{equation}
Take $r > 0$. Denote by $W^1(r\T)$ the linear space of functions $h \in L^1(\T)$ such that $\sum_{k\in \Z}r^k|\hat{h}(k)| < \infty$ equipped with the norm $\|\sum_{k\in \Z}c_kz^k\|_{W^1(r\T)} = \sum_{k\in\Z}r^k|c_k|$. If $h_{1,2} \in L^2(\T)$, then $h_1 h_2 \in L^1(\T)$ and the following multiplicative inequality holds:
\begin{align*}
\|h_1 h_2\|_{W^1(r\T)}
&=\sum_{k \in \Z} r^{k}|(\hat h_1 * \hat h_2)(k)| \\
&\le \sum_{k \in \Z} \sum_{j \in \Z} r^{k} |\hat h_1(j)\hat h_2(k-j)| \\
&= \sum_{j \in \Z}r^j|\hat h_1(j)|\sum_{k \in \Z}r^{k-j}|\hat h_2(k-j)| \\
&=\|h_1\|_{W^1(r\T)} \cdot \|h_2\|_{W^1(r\T)}.
\end{align*}
Consider the space 
$$
B^{1}(r\T) = \Bigl\{F \in S_*(\D): \;\; \sum_{k \ge 0}r^{k}|F_{k}(0)| < \infty\Bigr\}.
$$
For Schur functions $F, G \in \szc$, and $r \in (0, 1)$, we are going to estimate 
\begin{equation}\label{eq181}
\|F-G\|_{W^1(r\T)} = \sum_{k=0}^{\infty}r^k|\hat F(k) - \hat G(k)|
\end{equation}
in terms of 
\begin{gather}
\label{eq: SW metric def}
\rho_{B^1(r\T)}(F,G) := \sum_{k=0}^{\infty}r^k|F_{k}(0) - G_{k}(0)|.
\end{gather}
We set $\eta(F, G) = \min(\eta(F), \eta(G))$, where $\eta(\cdot)$ is defined in \eqref{eq: szego theorem}.
\begin{Thm}\label{t7}
For all $F, G \in \szc$, $r\in (0,1)$ such that $1-r\ge 12\log(\eta(F,G)^{-1})$ we have
\begin{equation}\label{eq18}
	\frac{1}{2}\rho_{B^1(r\T)}(F,G) \le
	\|F - G\|_{W^1(r\T)}
	\le 2\rho_{B^1(r\T)}(F,G).
\end{equation}
\end{Thm}
By definition, we have $B^1(r\T) \subset S_*(\D)$ for all $r \in (0,1)$. The situation changes for $r=1$ and $r > 1$. Schur functions $F$ lying in $B^1(\T)$ were considered by Baxter in his famous work \cite{Baxter}, see Theorem \ref{tB} above. For $r > 1$, Schur functions have exponentially decaying recurrence coefficients. This case was studied by P.\,Nevai and V.\,Totik in \cite{NT89}, see also Section 7 in \cite{Simonbook1}, \cite{Simon06}, \cite{GM06} and references in these works. However, we are not aware of inequalities analogous to \eqref{eq18} for $r \ge 1$. This remains an interesting open direction. We also do not know if \eqref{eq18} holds (possibly with some different constants) if we change $12$ in the statement of Theorem \ref{t7} by arbitrary $\eps> 0$. See Lemma \ref{l2p3} below for the case where $g = 0$.

\medskip

For the proof of Theorem \ref{t7} we need some auxiliary results. Take a function $F \in \szc$. The inequality $x \le -\log(1 - x)$ holds for all $x\in [0,1)$, let us apply it to \eqref{eq: szego theorem}. We have
\begin{align}
\label{eq: rec coeffs norm}
\sum_{k = 0}^\infty|F_k(0)|^2 \le \sum_{k = 0}^\infty-\log(1 - |F_k(0)|^2) &= \log \eta(F)^{-1},	\\
\label{eq: fourier coeffs norm}
\sum_{k = 0}^\infty|\hat{F}(k)|^2
= \int_{\T}|F|^2\,dm_{\T} \le \int_{\T}-\log(1 - |F|^2)\,dm_{\T} &= 
\log \eta(F)^{-1}. 
\end{align}
\begin{Lem}\label{lemma: 1/11 bound}
Assume that $F\in \szc$, $r\in (0,1)$ are such that $1-r\ge T\log\eta(F)^{-1}$ for some $T > 0$. Then $\|F_n\|^2_{W^1(r\T)}\le 1/T$ for every $n\ge 0$.
\end{Lem}
\beginpf From \eqref{eq: fourier coeffs norm} we get
\begin{align*}
\|F_n\|_{W^1(r\T)} = \sum_{k=0}^{\infty}r^k|\hat F_n(k)|&\le \bigg(\sum_{k=0}^{\infty}r^{2k}\cdot \sum_{k=0}^{\infty}|\hat F_n(k)|^2\bigg)^{1/2}
\le \sqrt{\frac{\log \eta(F_n)^{-1}}{1 - r^2}}.
\end{align*}
The definition \eqref{eq: szego theorem} of $\eta(F)$ yields $\eta(F)\le \eta(F_n)$ for each $n \ge 0$, hence
$$
\|F_n\|^2_{W^1(r\T)}\le \frac{\log \eta(F_n)^{-1}}{1 - r^2}\le \frac{\log \eta(F)^{-1}}{1 - r} \le \frac{1}{T}.
\eqno\qed
$$
\begin{Lem}\label{l2p3}
Assume that $F\in \szc$, $r\in (0,1)$ are such that $1-r\ge T\log\eta(F)^{-1}$ for some $T > 0$. Then 
$$
\frac{T - 1}{T}\sum_{k = 0}^\infty|F_k(0)|r^k\le \|F\|_{W^1(r\T)}\le \frac{T + 1}{T}\sum_{k = 0}^\infty|F_k(0)|r^k.
$$
\end{Lem}
\beginpf
One step $F_n \mapsto F_{n+1}$ of the Schur's algorithm \eqref{sa-bis} can be rewritten in the following form:
\begin{gather}
\label{eq: schur alg with R}
zF_{n + 1} = F_n - F_n(0) + z\ov{F_n(0)}F_n F_{n + 1}.
\end{gather}
Then, 
\begin{align*}
z^2F_{n + 1} 
&= 
zF_n - zF_n(0) + z^2\ov{F_n(0)}F_n F_{n + 1}\\
&=
F_{n-1} - F_{n-1}(0) - zF_n(0) + z\ov{F_{n-1}(0)}F_{n-1} F_{n} + z^2\ov{F_n(0)}F_n F_{n + 1}.
\end{align*}
Iterating \eqref{eq: schur alg with R} further, we obtain
\begin{gather*}
z^{n + 1} F_{n + 1} = F - \sum_{k  =0 }^{n}z^k F_k(0) + \sum_{k  =0 }^{n}z^{k + 1}\ov{F_k(0)}F_k F_{k+ 1}.
\end{gather*}
If we send $n\to\infty$, this becomes
\begin{gather}
\label{f in terms of Rf}
F = \sum_{k  =0 }^{\infty}z^k F_k(0) - \sum_{k = 0}^{\infty }z^{k + 1}\ov{F_k(0)}F_k F_{k+ 1}, \qquad z \in \D.
\end{gather}
Lemma \ref{lemma: 1/11 bound} applied for $F_k$ and $F_{k + 1}$ gives
\begin{gather*}
\bigg\|\sum_{k = 0}^{\infty }z^{k + 1}\ov{F_k(0)}F_k F_{k+ 1}\bigg\|_{W^1(r\T)}\le \sum_{k = 0}^\infty r^{k + 1}|F_k(0)|\cdot\|F_k\|_{W^1(r\T)}\cdot \|F_{k + 1}\|_{W^1(r\T)} \le \frac{1}{T}\sum_{k = 0}^\infty r^k|F_k(0)|,
\end{gather*}
where we used the multiplicative property of the norm $\| \cdot \|_{W^1(r\T)}$. The lemma follows from \eqref{f in terms of Rf} and the triangle inequality. \qed

\begin{Lem}\label{lem sum prod norm coeff}
Assume that $F,G\in \szc$, $r\in (0,1)$ are such that $1-r\ge T\log\eta(F,G)^{-1}$ for some $T > 0$. Then we have
$$
\sum_{k = 0}^\infty |G_k(0)|\cdot\|F_{k}\|_{W^1(r\T)}\le \frac{T + 1}{T^2}.
$$
\end{Lem}
\beginpf
From the previous lemma we get
$$
\|F_{k}\|_{W^1(r\T)}\le  \frac{T + 1}{T}\sum_{l = 0}^\infty r^{l}|F_{k + l}(0)|.
$$
We use this,  the Cauchy-Schwarz inequality,  and \eqref{eq: rec coeffs norm} to obtain 
\begin{align}
\sum_{k = 0}^\infty |G_k(0)|\cdot\|F_{k}\|_{W^1(r\T)}&\le \frac{T + 1}{T}\sum_{k = 0}^\infty |G_k(0)|\sum_{l = 0}^\infty r^{l}|F_{k + l}(0)|
= \frac{T + 1}{T}\sum_{l = 0}^\infty r^{l}\sum_{k = 0}^\infty |G_k(0)|\cdot |F_{k + l}(0)| \notag	\\
&\le \frac{T + 1}{T}\sum_{l = 0}^\infty r^{l} \cdot \sqrt{\log\eta(F_l)^{-1} \log\eta(G)^{-1}} \le \frac{(T + 1)\log\eta(F, G)^{-1}}{T(1-r)}.\notag
\end{align}
Since $1-r\ge T\log\eta(F,G)^{-1}$, the lemma is proved. \qed

\medskip

\noindent {\bf Proof of Theorem \ref{t7}.} 
We first prove Theorem  \ref{t7} under the assumption $F_{n + 1} = G_{n + 1} = 0$ for some $n\ge 0$. For simplicity, in some places we write $\|\cdot\|$ instead of $\|\cdot\|_{W^1(r\T)}$. Let us argue by induction. Fix $k \in \Z_+$ and assume that
\begin{equation}
\label{induction claim}
\frac{1}{2}\sum_{l = 0}^{\infty}r^{l}|F_{j + l}(0) - G_{j + l}(0)|\le\|F_j - G_j\|_{W^1(r\T)}\le 2\sum_{l = 0}^{\infty}r^{l}|F_{j + l}(0) - G_{j + l}(0)| 
\end{equation}
for every $j > k$. We want to prove \eqref{induction claim} for $j = k$. Note that \eqref{induction claim} clearly holds for $j \ge n + 1$, therefore, there is no problem with the initial step of induction. Relation \eqref{f in terms of Rf} for $F_k$, $G_k$ gives
\begin{align*}
F_k - G_k =& \bigg(\sum_{l  =0 }^{\infty}z^l F_{k + l}(0) - \sum_{l = 0}^{\infty}z^{l + 1} \ov{F_{k + l}(0)}F_{k + l}F_{k + l + 1}\bigg) -
\\
&- \bigg(\sum_{l  =0 }^{\infty}z^l G_{k + l}(0) - \sum_{l = 0}^{\infty}z^{l + 1} \ov{G_{k + l}(0)}G_{k + l}G_{k + l + 1}\bigg)
\\
=&\sum_{l  =0 }^{\infty}z^l (F_{k + l}(0) - G_{k + l}(0)) - \sum_{l  =0 }^{\infty}z^{l + 1}\big( \ov{F_{k + l}(0)}F_{k + l}F_{k + l + 1} - \ov{G_{k + l}(0)}G_{k + l}G_{k + l + 1}\big).
\end{align*}
For every $s = k+l$, the expression in the last brackets can be rewritten in the form 
\begin{align*}
\ov{F_{s}(0)}F_{s}F_{s + 1} - \ov{G_{s}(0)}G_{s}G_{s + 1} 
=& \,(\ov{F_s(0)} - \ov{G_s(0)})F_sF_{s + 1} + \\
&+ \ov{G_s(0)}F_s(F_{s + 1} - G_{s + 1}) + \ov{G_s(0)}G_{s + 1}(F_s - G_s).
\end{align*}
Therefore, we have
\begin{align}
\nonumber
F_k - G_k &= \sum_{l  =0 }^{\infty}z^l \big(F_{k + l}(0) - G_{k + l}(0)\big) - \sum_{l  =0 }^{\infty}z^{l + 1}F_{k + l}F_{k + l + 1}\big(\ov{F_{k + l}(0)} - \ov{G_{k + l}(0)} \big) -
\\
\label{eq: fk - gk 1}
&-\sum_{l  =0 }^{\infty}z^{l + 1}\ov{G_{k + l}(0)}F_{k + l}\big(F_{k + l + 1} - G_{k + l + 1}\big) - \sum_{l  =0 }^{\infty}z^{l + 1}\ov{G_{k + l}(0)}G_{k + l + 1}\big(F_{k + l} - G_{k + l}\big).
\end{align}
For the second term in the latter sum we use the triangle inequality, the multiplicative property of the norm $\| \cdot \|_{W^1(r\T)}$, and Lemma \ref{lemma: 1/11 bound} to obtain 
\begin{align}
\nonumber
\bigg\|\sum_{l  =0 }^{\infty}z^{l + 1}F_{k + l}F_{k + l + 1}\big(\ov{F_{k + l}(0)} - \ov{G_{k + l}(0)} \big)\bigg\| &\le
\sum_{l  =0 }^{\infty}r^{l + 1}\|F_{k + l}\|\cdot\|F_{k + l + 1}\|\cdot\big|F_{k + l}(0) - G_{k + l}(0) \big|
\\
\label{eq: fk - gk 2}
&\le \frac{1}{T}\sum_{l  =0 }^{\infty}r^l|F_{k + l}(0) - G_{k + l}(0)|.
\end{align}
Let us rewrite the third and the fourth term in \eqref{eq: fk - gk 1} in the form 
\begin{align*}
\sum_{l=0}^{\infty}z^{l + 1}\ov{G_{k + l}(0)}F_{k + l}\big(F_{k + l + 1} - G_{k + l + 1}\big) 
=& 
\sum_{l=1}^{\infty}z^{l}\ov{G_{k + l -1}(0)}F_{k + l - 1}\big(F_{k + l} - G_{k + l}\big),\\
\sum_{l  =0 }^{\infty}z^{l + 1}\ov{G_{k + l}(0)}G_{k + l + 1}\big(F_{k + l} - G_{k + l}\big) 
=&\sum_{l  = 1}^{\infty}z^{l}z\ov{G_{k + l}(0)}G_{k + l + 1}\big(F_{k + l} - G_{k + l}\big)\\
&+ z\ov{G_{k}(0)}G_{k+1}\big(F_{k} - G_{k}\big).
\end{align*}
We see that the sum of these two expressions equals
\begin{gather*}
z \ov{G_{k}(0)}G_{k + 1}(F_{k} - G_{k}) + \sum_{l  =1 }^{\infty}z^{l }\big(\ov{G_{k + l - 1}(0)} F_{k + l - 1} + z\ov{G_{k + l}(0)}G_{k + l + 1}\big)\big(F_{k + l} - G_{k + l}\big).
\end{gather*}
For every $l\ge 1$ we have $k + l > k$, hence $F_{k + l} - G_{k + l}$ can be estimated by the induction assumption \eqref{induction claim}. Together with Lemma \ref{lem sum prod norm coeff} this gives (we use again the multiplicative property of the norm $\| \cdot \|_{W^1(r\T)}$)
\begin{align}
\nonumber
\bigg\|\sum_{l  =1 }^{\infty}z^{l }&(z\ov{G_{k + l}(0)}G_{k + l + 1} + \ov{G_{k + l - 1}(0)} F_{k + l - 1})\big(F_{k + l} - G_{k + l}\big)\bigg\| \le
\\
\nonumber
&\le 2\sum_{l  =1 }^{\infty}r^{l }\big(|G_{k + l}(0)|\cdot\|G_{k + l + 1}\| + |G_{k + l - 1}(0)|\cdot\| F_{k + l - 1}\|\big)\sum_{t = 0}^\infty r^{t}|F_{k + l + t}(0) - G_{k + l + t}(0)|
\\
\nonumber
&\le 2\sum_{l = 1}^{\infty}\big(|G_{k + l}(0)|\cdot\|G_{k + l + 1}\| + |G_{k + l - 1}(0)|\cdot\| F_{k + l - 1}\|\big) \cdot \sup_{\ell \ge 1}  \sum_{t = 0}^\infty r^{t+\ell}|F_{k + \ell + t}(0) - G_{k + \ell + t}(0)|
\\
\label{eq: fk - gk 3}
&\le \frac{4(T + 1)}{T^2} \sum_{s = 1}^\infty r^s|F_{k + s}(0) - G_{k + s}(0)|.
\end{align}
Denote $C(T) = \frac{4(T + 1)}{T^2} + \frac{1}{T} = \frac{5T + 4}{T^2}$.
The substitution of \eqref{eq: fk - gk 2} and \eqref{eq: fk - gk 3} into \eqref{eq: fk - gk 1} implies
\begin{multline}\label{eq: dsbnidsere}
\left(1 - C(T)\right)\sum_{s = 0}^\infty r^s|F_{k + s}(0) - G_{k + s}(0)|\le \\
\le \|(F_k - G_k) - z\ov{G_{k}(0)}G_{k + 1}(F_{k} - G_{k})\|  \le \\
\le \left(1 + C(T)\right)\sum_{s = 0}^\infty r^s|F_{k + s}(0) - G_{k + s}(0)|.
\end{multline}
By Lemma \ref{lemma: 1/11 bound}, we have $\|\ov{G_k(0)}G_{k + 1}\|\le \|G_k\|\cdot \|G_{k + 1}\| \le 1/T$, hence
\begin{gather}
\label{eq: dsbnidsere2}
\frac{T-1}{T}\|F_k - G_k\|\le \|(F_k - G_k) - z\ov{G_{k}(0)}G_{k + 1}(F_{k} - G_{k}) \|  \le \frac{T+1}{T}\|F_k - G_k\|.
\end{gather}
Therefore, for $T > 1$ from \eqref{eq: dsbnidsere} and \eqref{eq: dsbnidsere2} we see that
\begin{align*}
\frac{T(1-C(T))}{T+1}
\sum_{l = 0}^{\infty}r^{l}|F_{k + l}(0) - G_{k + l}(0)| \le \|F_k - G_k\|_{W^1(r\T)},\\ \frac{T(1+C(T))}{T-1}\sum_{l = 0}^{\infty}r^{l}|F_{k + l}(0) - G_{k + l}(0)| \ge \|F_k - G_k\|_{W^1(r\T)}.
\end{align*}
To complete the induction step, it remains to check that for $T = 12$ we have
$$
\frac{1}{2} \le \frac{T(1-C(T))}{T+1},
\qquad 
\frac{T(1+C(T))}{T-1} \le 2, 
$$
where $C(T) = \frac{5T + 4}{T^2}$. This is indeed the case. 

\medskip

To finish the proof, we need to get rid of the assumption $F_{n + 1} = G_{n + 1} = 0$. For this, we take arbitrary Schur functions $F, G \in \szc$, fix $n \in \Z_+$, and consider the functions $F^{(n)}, G^{(n)} \in \szc$ such that 
$$
F^{(n)}_k(0)
=
\begin{cases}
F_k(0), & k \le n,\\
0, & k > n,
\end{cases}
\qquad 
G^{(n)}_k(0)
=
\begin{cases}
G_k(0), & k \le n,\\
0, & k > n.
\end{cases}
$$
From \eqref{eq: szego theorem} we know that $\eta(F^{(n)}, G^{(n)})\ge \eta(F, G)$, hence the previous part of the proof works for $F^{(n)}$ and $G^{(n)}$. It gives
\begin{gather*}
\frac{1}{2}\|F^{(n)} - G^{(n)}\|_{W^1(r\T)} \le
\rho_{B^1(r\T)}(F^{(n)}, G^{(n)})
\le
2\|F^{(n)}-G^{(n)}\|_{W^1(r\T)}.
\end{gather*}
As $n\to\infty$, we have
$$
\rho_{B^1(r\T)}(F^{(n)}, G^{(n)}) \to \rho_{B^1(r\T)}(F, G),\qquad 
\|F^{(n)} - G^{(n)}\|_{W^1(r\T)}\to \|F - G\|_{W^1(r\T)}.
$$
Indeed, the first convergence immediately follows from the definition \eqref{eq: SW metric def} of the metric $\rho_{B^1(r\T)}$; the second convergence holds because $\hat{F}(k) = \widehat{F^{(n)}}(k)$, $\hat{G}(k) = \widehat{G^{(n)}}(k)$ for all $n \ge k$, see Theorem 1.5.5 in \cite{Simonbook1}.
\qed

\medskip

\subsection{Approximation of \texorpdfstring{$L^2(\R_+)$}{L2(R+)}-potentials by discrete measures}
In this short section we prove some auxiliary lemmas related to the approximation of functions in $L^2(\R_+)$ by measures supported on $\ell \Z_+$ as $\ell \to 0$. These lemmas will be used in the proof of Theorem~\ref{t1}. Let $q \in L^2(\R_+)$, $\ell > 0$. From now on and till the end of the paper, we denote
\begin{equation} \label{eq: qr def}
q_{\ell} = \sum_{k=0}^{\infty}q[\ell k]\delta_{\ell k}, \qquad q[\ell k] = \int_{\ell k}^{\ell(k+1)}q(x)\,dx,
\qquad
Q[\ell k] = \begin{pmatrix}\Im q[\ell k] & \Re q[\ell k] \\ \Re q[\ell k] & -\Im q[\ell k]\end{pmatrix},
\end{equation}
where $\delta_{\ell k}$ is the point mass measure concentrated at the point $\ell k$. This notation should not be confused with the notation $q_{\ell} = q(\ell+\cdot)$ used in previous sections. 
\begin{Lem}\label{l36}
For $q \in L^2(\R_+)$, denote by  $\Hh_{q}$, $\Hh_{q_{\ell}}$ the Hamiltonians on $\R_+$ corresponding to $q$, $q_{\ell}$, respectively, see \eqref{eq28}. Then $\Hh_{q_{\ell}} \to \Hh_{q}$ as $\ell \to 0$ uniformly on compact subsets of $\R_+$. 
\end{Lem}
\beginpf Fix an arbitrary $L > 0$. Since $\Hh_q = N_q^* N_q$ and $\Hh_{q_\ell} = N_{q_\ell}^* N_{q_\ell}$, it suffices to show that $\|N_{q_{\ell}} - N_{q}\|$ tends to $0$ uniformly on $[0, L]$ as $\ell \to 0$. Take $x\in [0, L]$ and set $k_* = \lfloor x/\ell\rfloor$. We have 
$N_{q}(x) = N_q(x,\ell k_*) N_q(\ell k_*, \ell(k_*-1)) \cdot \ldots \cdot  N_q(2\ell, \ell) N_{q}(\ell)$, see \eqref{eq29}. For every $\ell > 0$, $k\in \Z_+$, let $\Delta_{\ell, k} = [\ell k, \ell k + \ell)$. Then,
\begin{gather*}
\int_{\Delta_{\ell, k}}|q(x)|\,dx \le \sqrt{\ell}\left(\int_{\Delta_{\ell, k}}|q(x)|^2\,dx\right)^{1/2} = \sqrt{\ell}\|q\|_{L^2(\Delta_{\ell, k})}.
\end{gather*}
Series representation \eqref{eq: series representation} and inequality \eqref{eq: series bound} give
\begin{align}
\label{tmp: 5:14}
N_q(\ell(k+1), \ell k) &= \IDM + 
\int_{\Delta_{\ell, k}}JQ(x_1)\,dx_1 + O\left(\ell \|q\|^2_{L^2(\Delta_{\ell, k})}\right),\\
\label{tmp: 5:15}
\|N_q(\ell(k+1), \ell k)\| &= 1 + O\left(\|q\|_{L^1(\Delta_{\ell, k})}\right),
\end{align}
where $O(\cdot)$ is uniform with respect to $\ell \in [0,1]$ and $k \in \Z_+$. By the definition \eqref{eq33} of $N_{q_{\ell}}$, we have $N_{q_{\ell}}(x) = e^{JQ[\ell k_*]} e^{JQ[\ell(k_* - 1)]} \cdot \ldots e^{JQ[\ell]}e^{JQ[0]}$. Relation \eqref{tmp: 5:14} gives
\begin{align}
\label{tmp: 5.16}
e^{JQ[\ell k]} &= \IDM + JQ[\ell k] + O(\|Q[\ell k]\|^2) = N_q(\ell(k+1), \ell k) + O\left(\ell \|q\|^2_{L^2(\Delta_{\ell, k})}\right).
\end{align}
As before, we have 
\begin{align}
\label{tmp: 5.17}
\|e^{JQ[\ell k]}\| &= 1 + O\left(\|q\|_{L^1(\Delta_{\ell, k})}\right).
\end{align}
Let us use the following telescopic sum relation
$$
a_1 a_2 \cdots a_{k_*} - b_1 b_2 \cdots b_{k_*} = \sum_{k=1}^{k_*} a_1 \cdots a_{k-1}(a_k - b_k) b_{k+1} \cdots b_{k_*},
$$
where the empty products (appearing for $k=1$ and $k=k_*$ in the r.h.s.) are understood as $1$. 
We obtain
$$
\|N_q(\ell k_*) - N_{q_{\ell}}(x)\| \le C_{q, x} C_{q_{\ell}, x}\sum_{k = 1}^{k_*}\|N_q(\ell k, \ell (k-1)) - e^{JQ[\ell k]}\|,
$$
for
$$
C_{q_{\ell}, x} = \max_{1 \le j \le k_*} \prod_{j = 1}^{k - 1}\|N_q(\ell j, \ell (j-1))\|, \qquad 
C_{q, x} = \max_{1 \le j \le k_*} \prod_{j = k + 1}^{k_*}\|e^{JQ[\ell j]}\|.
$$
The estimates \eqref{tmp: 5:15}, \eqref{tmp: 5.17} tell us that $C_{q_{\ell}, x}$, $C_{q, x}$ are uniformly bounded in $\ell  \in [0, 1]$, $x \in [0, L]$. Relation \eqref{tmp: 5.16} implies
$$	
\sum_{k = 0}^{k_*}\|N_q(\ell(k+1), \ell k) - e^{JQ[\ell k]}\| = O\left(\ell \|q\|^2_{L^2([0,x])}\right),
$$
which tends to zero as $\ell \to 0$ (the constant in $O(\cdot)$ is uniform in $\ell  \in [0,1]$, $x\in [0,L]$). Finally, \eqref{tmp: 5:15} implies that $\|N(x) - N_q(\ell k_*)\| \to 0$ as $\ell \to 0$ uniformly in $x\in [0,L]$, and the proof is completed. \qed 

\medskip 

\begin{Lem}\label{lem: first discr}
For every $q\in L^2(\R_+)$ we have 
\begin{equation}\label{eq91}
	\lim_{\ell \to 0}\frac{1}{\ell }\sum_{k\ge 0} |q[\ell k]|^2 = \|q\|_{L^2(\R_+)}^2.  
\end{equation}
\end{Lem}
\beginpf
Recall that $\Delta_{\ell, k} = [\ell k, \ell k+\ell)$. The Cauchy-Schwarz inequality gives
\begin{gather}\label{eq: lem1 ineq}
\frac{1}{\ell}\sum_{k\ge 0} |q[\ell k]|^2 = \frac{1}{\ell}\sum_{k\ge 0}\bigg|\int_{\Delta_{\ell, k}}q(x)\,dx\bigg|^2 \le \frac{1}{\ell}\sum_{k\ge 0} \ell\cdot \int_{\Delta_{\ell, k}}|q(x)|^2\,dx = \|q\|_{L^2(\R_+)}^2.
\end{gather}
In particular, both sides of \eqref{eq91} with $\lim$ replaced by $\limsup$ or $\liminf$ depend continuously on $q$ in $L^2(\R_+)$-norm. Hence, it suffices to prove \eqref{eq91} only on a dense subset of $L^2(\R_+)$. Let $C_0^\infty(\R_+)$ be the set of infinitely smooth functions on $\R_+$ with compact support. Take some $q\in C_0^\infty(\R_+)$ and let $R> 0$ be such that $\supp q\subset [0, R]$. Then 
\begin{gather*}
\frac{1}{\ell}|q[\ell k]|^2 - \int_{\Delta_{\ell, k}}|q(x)|^2\,dx = \int_{\Delta_{\ell, k}}q(x)\bigg(\frac{1}{\ell}\int_{\Delta_{\ell, k}}\ov{q(y)}\,dy - \ov{q(x)}\bigg)\,dx = O(\ell^2), \quad \ell \to 0.
\end{gather*}
We also have $q[\ell k] = 0$ for $k > R/ \ell$, hence
\begin{align*}
\frac{1}{\ell}\sum_{k\ge 0} |q[\ell k]|^2 = \sum_{k = 0}^{[R/\ell]}\frac{1}{\ell}|q[\ell k]|^2 &= \sum_{k = 0}^{[R/\ell]}\bigg(\int_{\Delta_{\ell, k}}|q(x)|^2\,dx + O(\ell^2)\bigg) 
\\
&= \int_{0}^{\infty}|q(x)|^2\,dx + O(\ell^2 \cdot R/\ell) = \|q\|_{L^2}^2 + O(R\ell),
\end{align*}
which tends to zero with $\ell \to 0$ because $R$ is fixed. \qed

\begin{Lem}\label{final lem 6.2}
For $q\in L^2(\R_+)$ and $A > 0$ we have
\begin{gather}
	\label{eq: lem 62}
	\lim_{\ell \to 0}\sum_{k\ge 0} e^{-A\ell k}|q[\ell k]| = \int_{\R_+} e^{-Ax}|q(x)|dx.
\end{gather}
\end{Lem}
\beginpf
The argument is similar to the proof of the previous lemma. It is easy to check using Cauchy-Schwarz inequality and \eqref{eq: lem1 ineq} that both sides of \eqref{eq: lem 62} with $\lim$ replaced by $\limsup$ or $\liminf$ depend on $q$ continuously in $L^2(\R_+)$-norm. Therefore, we need to prove \eqref{eq: lem 62} only for $q\in C_0^\infty(\R_+)$. As $\ell\to 0$, we have the following estimates
\begin{gather*}
|q[\ell k]| = \ell |q(\ell k)| + O(\ell^2),\qquad \int_{\Delta_{\ell, k}} |q(x)|\,dx = \ell |q(\ell k)| + O(\ell^2),\\
\int_{\Delta_{\ell, k}} (e^{-A\ell k} - e^{-Ax})\,dx = O(e^{-A\ell k}\ell^2). 
\end{gather*}
Therefore, 
\begin{align*}
e^{-A\ell k}|q[\ell k]| - \int_{\Delta_{\ell, k}} e^{-Ax}|q(x)|dx &= e^{-A\ell k}\bigg(|q[\ell k]| - \int_{\Delta_{\ell, k}} |q(x)|\,dx\bigg)
\\
&+ \int_{\Delta_{\ell, k}} (e^{-A\ell k} - e^{-Ax})|q(x)|dx = O(e^{-A\ell k}\ell^2),\quad \ell \to 0.
\end{align*}
Summing up over all $k\ge 0$, we get
\begin{gather*}
\bigg|\sum_{k\ge 0} e^{-A\ell k}|q[\ell k]| - \int_{\R_+} e^{-Ax}|q(x)|dx\bigg| = O\bigg(\sum_{k \ge 0}e^{-A\ell k}\ell^2\bigg) = O\bigg(\frac{\ell^2}{1 - e^{-A\ell}}\bigg) = O(\ell),
\end{gather*}
which completes the proof. \qed

\medskip 

\subsection{Approximation of \texorpdfstring{$S_2(\C_+)$}{S2(C+)}-functions by periodic Schur functions} Take $q \in L^2(\R_+)$, $\ell > 0$, and define $q_{\ell}\in \cM_{\ell}$ by \eqref{eq: qr def}. Let $f_{q_{\ell}}$  be the Schur function of $\Di_{q_{\ell}}$.
By Lemma \ref{pl1}, this function is periodic with period $\pi/\ell$. From Lemma \ref{l10} we know that there exists a Schur function $F_{q_{\ell}}$ in $\D$ that satisfies 
\begin{gather}
\label{eq: fqr and Fqr}
f_{q_{\ell}}(z) = F_{q_{\ell}}(e^{2i\ell z}),\qquad f_{q_{\ell}}(\infty) = F_{q_{\ell}}(0),\qquad z\in \C_+.
\end{gather}
It is natural to expect that $f_{q_\ell}$ in some sense approximates $f_{q}$ as $\ell \to 0$. This is indeed the case, we study this approximation below. For $\ell > 0$, we let $g_{q_{\ell}}$ be such that $g_{q_{\ell}}(x) = f_{q_{\ell}}(x)$ for $|x| \le \pi/(2\ell)$ and $g_{q_{\ell}}(x) = 0$ for $|x| > \pi/(2\ell)$. Note that $g_{q_{\ell}}$ is not analytic and defined only on $\R$. We are going to prove the convergence of $g_{q_\ell}$ to $f_q$ in the following metric space $X$ of measurable functions on $\R$,
\begin{gather}
\label{eq: spacex}
X = \{r\in L^\infty(\R)\colon \|r\|_{L^{\infty}(\R)}\le 1,\;\; \log(1 - |r|^2)\in L^1(\R)\}.
\end{gather}
For every $r, r_1, r_2\in X$ define
\begin{gather}
\label{eq: metric in X}
\rho_{X}^2(r_1, r_2) = \int_{\R}-\log\bigg(1 - \Big|\frac{r_1 - r_2}{1 - \ov{r_1}r_2}\Big|^2\bigg)\,dx = \int_{\R}-\log\bigg(\frac{(1 - |r_1|^2)(1 - |r_2|^2)}{|1 - \ov{r_1}r_2|^2}\bigg)\,dx,
\\ 
E(r_1, r_2) =  \int_{\R}-\log\left|1 - \ov{r_1}r_2 \right|\,dx, \qquad E(r) = E(r, r). \label{eq99}
\end{gather}
Formula \eqref{eq: metric in X} can be written in the following form:
\begin{gather}
\label{eq: dist in terms of E}
\rho_{X}^2(r_1, r_2) = E(r_1) + E(r_2) - 2E(r_1, r_2).
\end{gather}
Theorem 1.6 in \cite{Sylvester} states that $\rho_X$ is a metric on $X$. In particular, it satisfies the triangle inequality
\begin{gather}
\label{eq: rho X triangle ineqaulity}
\rho_X(r_1, r_3)\le \rho_X(r_1, r_2) + \rho_X(r_2, r_3), \qquad r_1, r_2, r_3\in X.
\end{gather}
This follows from the fact that the function $d\colon \D\times \D\to \R$, 
\begin{gather*}
d^2(z,w) = -\log\bigg(1 - \left|\frac{z - w}{1 - \ov{z}w}\right|^2\bigg)
\end{gather*}
is a metric on $\D$. Indeed, by Cauchy-Schwarz inequality, we have 
\begin{align*}
\rho_X(r_1, r_3) &= \int_\R d^2(r_1, r_3)\, dx = \int_\R \big(d(r_1, r_2) + d(r_2, r_3)\big)^2\, dx \\
&\le \int_\R d^2(r_1, r_2)\, dx + \int_\R d^2(r_2, r_3)\, dx + 2\sqrt{\int_\R d^2(r_1, r_2)\, dx\int_\R d^2(r_2, r_3)\, dx} \\
&= \rho_X^2(r_1, r_2) + \rho_X^2(r_2, r_3) + 2\rho_X(r_1, r_2) \rho_X(r_2, r_3) = \big(\rho_X(r_1, r_2) + \rho_X(r_2, r_3) \big)^2.  
\end{align*}
One can see that the space $S_{2}(\C_+)$ in \eqref{eq42} is a closed subspace of $X$ with the induced norm. To prove convergence of $g_{q_\ell}$ to $f_q$ in $X$, we will need the following lemma.
\begin{Lem}\label{lem: eta in limit}
For every $q\in L^2(\R_+)$ we have
\begin{gather}
	\label{eq: in lemma 6.3}
	\frac{1}{\pi}\lim_{\ell \to 0}E(g_{q_{\ell}}) =\lim_{\ell \to 0} \frac{1}{\ell}\log(\eta(F_{q_{\ell}})) = - \|q\|^2_{L^2(\R_+)}.
\end{gather}
\end{Lem}
\beginpf
We have $f_{q_{\ell}}(x) = F_{q_{\ell}}(e^{2i\ell x})$ a.\,e.\,on $[-\pi/2\ell, \pi/2\ell]$, hence
\begin{align*}
E(g_{q_{\ell}}) &= \int_{-\pi/2\ell}^{\pi/2\ell}\log(1-|f_{q_{\ell}}(x)|^2)\,dx = \int_{-\pi/2\ell}^{\pi/2\ell}\log(1-|F_{q_{\ell}}(e^{2i\ell x})|^2)\,dx\\
&= \frac{1}{2\ell}\int_{-\pi}^{\pi}\log(1-|F_{q_{\ell}}(e^{iy})|^2)\,dy = \frac{\pi}{\ell}\log \eta(F_{q_{\ell}}),
\end{align*}    
where the factor $2\pi$ appears in the last inequality because of the normalization of the measure $m_{\T}$ used in \eqref{eq: szego theorem}. 
Therefore it suffices to prove only the second equality in \eqref{eq: in lemma 6.3}.
From the Szeg\H{o} theorem, formula \eqref{eq: szego theorem}, and \eqref{eq: fqr and Fqr} we know that
\begin{gather}
\label{eq: szego theorem applied to Fqr}
\log(\eta(F_{q_{\ell}})) = \sum_{k\ge 0}\log (1 - |F_{q_{\ell},k}(0)|^2) = \sum_{k\ge 0}\log (1 - |f_{q_{\ell},k}(\infty)|^2).
\end{gather}
Theorem \ref{t2} states that $\varkappa(f_{q_{\ell},k}(\infty))=q[\ell k]$, where $\varkappa$ is defined in \eqref{eq44}. As $w \to 0$, we have 
$$
\varkappa(w) - \ov{w} = o(|w|), \qquad	|\varkappa(w)|^2 - |w|^2 = o\left(|w|^2 + |\varkappa(w)|^2\right).
$$
Therefore from \eqref{eq: lem1 ineq} and Lemma \ref{lem: first discr} we get
\begin{gather}
\nonumber
\sum_{k\ge 0}\left||f_{q_{\ell},k}(\infty)|^2 - |q[\ell k]|^2\right| = \sum_{k\ge 0}o\left(|q[\ell k]|^2\right) = o\left(\ell\right),\qquad \ell \to 0,
\\
\label{eq: sum fqrk(infty)}
\lim_{\ell \to 0}\frac{1}{\ell}\sum_{k\ge 0}|f_{q_{\ell},k}(\infty)|^2 = \|q\|^2_{L^2(\R_+)}.
\end{gather}
In particular, we get $\sup_k |f_{q_{\ell},k}(\infty)|^2 = O(\ell)$ uniformly for small $\ell$. It follows that  
\begin{gather*}
\left||f_{q_{\ell},k}(\infty)|^2 - \log (1 - |f_{q_{\ell},k}(\infty)|^2)\right| = O\left(\left|f_{q_{\ell},k}(\infty)\right|^4\right)  = O\left(\ell\left|f_{q_{\ell},k}(\infty)\right|^2\right),\qquad \ell \to 0,
\\
\sum_{k\ge 0}\left||f_{q_{\ell},k}(\infty)|^2 - \log (1 - |f_{q_{\ell},k}(\infty)|^2)\right| = O\bigg(\ell\sum_{k\ge 0}\left|f_{q_{\ell},k}(\infty)\right|^2\bigg) = O(\ell^2),\qquad \ell \to 0.
\end{gather*}
The latter combined with \eqref{eq: szego theorem applied to Fqr} and \eqref{eq: sum fqrk(infty)} concludes the proof. \qed

\medskip

\begin{Lem}\label{lem: convergence in X}
The functions $g_{q_\ell}$ converge to $f_q$ in $X$ as $\ell\to 0$.
\end{Lem}

The proof of this lemma uses two function-theoretic results whose proofs we postpone until the next section. Firstly, if $r_n, r\in S_2(\C_+)$ are such that $r_n \to r$ uniformly on compact subsets in $\C$, then 
\begin{equation}\label{eq95}
\lim\limits_{n\to\infty}E(r_n, \phi) = E(r, \phi) \; \mbox{ for every } \; \phi\in L^1(\R):\;\; \|\phi\|_{L^\infty(\R)} < 1.
\end{equation}
This is Lemma \ref{lem: conv on compacts imply conv in S2}. Secondly, we have $\rho_X(r_n,r) \to 0$ for  $r_n, r\in X$ if and only if
\begin{equation}\label{eq96}
\lim\limits_{n\to\infty}E(r_n) = E(r) \; \mbox{ and \eqref{eq95} holds.}
\end{equation}
This is assertion $(d)$ of Lemma \ref{lem: conv in S2 and conv in measure v2}.

\medskip

\noindent{\bf Proof of Lemma \ref{lem: convergence in X}.} By Lemma \ref{lem: eta in limit} and the isometric relation \eqref{eqSW} in Theorem \ref{tSW}, we have
\begin{align*}
\lim_{\ell\to 0}E(g_{q_\ell}) &=\lim_{\ell \to 0}\int_{-\pi/2\ell}^{\pi/2\ell}-\log(1 - |f_{q_{\ell}}(x)|^2)\,dx 
\\
&= \pi \|q\|^2_{L^2(\R_+)} =\int_\R-\log(1 - |f_q|^2)\,dx = E(f_q).
\end{align*}
So, we only need to prove that $E(g_{q_\ell}, \phi) \to E(f_q, \phi)$ for every function $\phi\in L^1(\R)$ satisfying $\|\phi\|_{L^\infty(\R)} < 1$. Fix such a function $\phi$. Lemma \ref{l36}, Krein -- de Branges Theorem \ref{KdB} and  \eqref{eq95} imply
\begin{gather*}
\lim_{\ell \to 0}\int_\R \log|1 - \ov{f_{q_{\ell}}}\phi|\,dx= \int_\R \log |1 - \ov{f_q}\phi |\,dx.
\end{gather*}
Furthermore, we have 
\begin{align*}
\bigg|\int_\R \log|1 - \ov{f_{q_{\ell}}}\phi|\,dx - \int_\R \log|1 - \ov{g_{q_{\ell}}}\phi|\,dx\bigg| &= \bigg|\int_{|x| > \pi/2\ell} \log|1 - \ov{f_{q_{\ell}}}\phi|\,dx\bigg| 
\le C \int_{|x| > \pi/2\ell}|\phi(x)|\,dx,
\end{align*}
for a constant  $C$ depending only on $\|\phi\|_{L^{\infty}(\R)}$. The r.h.s.\ in the last formula tends to $0$ as $\ell\to 0$. Thus, we get
$\lim_{\ell \to 0} E(g_{q_\ell}, \phi) = \lim_{\ell \to 0} E(f_{q_\ell}, \phi) = \lim_{\ell \to 0} E(f_{q}, \phi)$, and the proof is concluded by the application of \eqref{eq96}. \qed

\medskip

\begin{Lem}\label{final lem 6.11}
For $q, \tilde q\in L^2(\R_+)$ and $A > 0$, we have
\begin{gather*}
\lim_{\ell \to 0}\sum_{k\ge 0}|\hat F_{q_{\ell}}(k) - \hat F_{\tilde q_{\ell}}(k)|e^{-A\ell k}	
= \frac{1}{\sqrt{2\pi}}\int_{\R_+}e^{-A\xi/2}|\hat f_{q} - \hat f_{\tilde q}|\,d\xi.
\end{gather*}

\end{Lem}
\beginpf
First of all, let us relate $\hat{F}_{q_{\ell}}$ and $\hat{g}_{q_{\ell}}$. For every $\ell > 0$ and $k\in \Z_+$, we have 
\begin{align*}
\hat{g}_{q_{\ell}}(2\ell k) &= \frac{1}{\sqrt{2\pi}}\int_{\R}{g}_{q_{\ell}}(x)e^{-i2\ell kx}\,dx = \frac{1}{\sqrt{2\pi}}\int_{-\pi/2\ell}^{\pi/2\ell}{f}_{q_{\ell}}(x)e^{-i2\ell kx}\,dx\\
&= \frac{1}{\sqrt{2\pi}}\int_{-\pi/2\ell}^{\pi/2\ell}{F}_{q_{\ell}}(e^{2i\ell x})e^{-i2\ell kx}\,dx 
= \frac{\sqrt{2\pi}}{2\ell}\cdot \frac{1}{2\pi}\int_{-\pi}^{\pi}{F}_{q_{\ell}}(e^{iy})e^{-iky}\,dy  = \frac{\sqrt{2\pi}}{2\ell}\hat{F}_{q_{\ell}}(k).
\end{align*}
It follows that 
\begin{gather}
\label{eq: fourier coeffs F in terms of g}
\sum_{k\ge 0}|\hat F_{q_{\ell}}(k) - \hat F_{\tilde q_{\ell}}(k)|e^{-A\ell k} = \frac{1}{\sqrt{2\pi}}\sum_{k\ge 0}2\ell|\hat{g}_{q_{\ell}}(2\ell k) - \hat{g}_{\tilde q_{\ell}}(2\ell k)|e^{-A\ell k}.
\end{gather}
Take arbitrary $r_1, r_2\in X$ and observe that the inequality $-\log(1 - x)\ge x$ for $x \in [0,1)$ implies
\begin{gather*}
\rho^2_{X}(r_1, r_2)  = \int_{\R}-\log\bigg(1 - \Big|\frac{r_1 - r_2}{1 - \ov{r_1}r_2}\Big|^2\bigg)\,dx\ge \int_{\R} \Big|\frac{r_1 - r_2}{1 - \ov{r_1}r_2}\Big|^2\,dx \ge \frac{1}{4}\int_{\R} \left|r_1 - r_2\right|^2\,dx.
\end{gather*}
It follows that for every $r_1, r_2\in X$ we have $2\rho_{X}(r_1, r_2)\ge \|r_1 - r_2\|_{L^2(\R)}$.  Then Lemma \ref{lem: convergence in X} gives 
$$
\label{eq: gr to g in L2}
\lim_{\ell \to 0}\|g_{q_{\ell}} - f_q\|_{L^2(\R)} = 0, \qquad \lim_{\ell \to 0}\| g_{\tilde q_{\ell}} - f_{\tilde q}\|_{L^2(\R)} = 0.
$$
The latter, the Cauchy-Schwarz inequality and the isometric property of the Fourier transform  imply 
\begin{gather*}
\lim_{\ell \to 0}\int_{\R_+}|\hat{g}_{q_{\ell}}(\xi) - \hat{g}_{\tilde q_{\ell}}(\xi)|e^{-A\xi/2}\,d\xi =  \int_{\R_+}|\hat f_{q}(\xi) - \hat f_{\tilde q}(\xi)|e^{-A\xi/2}\,d\xi.
\end{gather*}
Thus, the claim of the lemma will follow from
\begin{gather}
\label{eq: darboux sum}
\lim_{\ell \to 0}\bigg|\int_{\R_+}|\hat{g}_{q_{\ell}}(\xi) - \hat{g}_{\tilde q_{\ell}}(\xi)|e^{-A\xi/2}\,d\xi - \sum_{k\ge 0}2r|\hat{g}_{q_{\ell}}(2\ell k) - \hat{g}_{\tilde q_{\ell}}(2\ell k)|e^{-A\ell k}\bigg| = 0.
\end{gather}
Denote  $h = \hat{g}_{q_{\ell}} - \hat{g}_{\tilde q_{\ell}}$. Then the value under the limit in  \eqref{eq: darboux sum} equals
\begin{gather}
\label{eq: darboux sum v2}
\bigg|\int_{\R_+}|h(\xi)|e^{-A\xi/2}\,d\xi - \sum_{k\ge 0}2r|h(2\ell k)|e^{-A\ell k}\bigg|  \le\sum_{k\ge 0}\int_{2\ell k}^{2\ell(k + 1)}\left|h(\xi)e^{-A\xi/2} - h(2\ell k)e^{-A\ell k}\right|\,d\xi.
\end{gather}
For every $\xi\in \R$ we have $|h(\xi)|\le \|g_{q_{\ell}} - g_{\tilde q_{\ell}}\|_{L^1(\R)}$. Moreover, there is an absolute constant $C > 0$ such that uniformly for all small $\ell > 0$, all $k\in\Z_+$ and $\xi\in[2\ell k, 2\ell(k+ 1)]$ we have
\begin{gather*}
|h(\xi) - h(2\ell k)|\le C|\xi - 2\ell k|\cdot\|g_{q_{\ell}} - g_{\tilde q_{\ell}}\|_{L^1(\R)},\qquad |e^{-A\xi/2} - e^{-A\ell k}|\le C|\xi - 2\ell k|e^{-A\ell k}.
\end{gather*}
It follows that 
\begin{gather*}
\left|h(\xi)e^{-A\xi/2} - h(2\ell k)e^{-A\ell k}\right| \le 2Ce^{-A\ell k}|\xi - 2\ell k|\cdot\|g_{q_{\ell}} - g_{\tilde q_{\ell}}\|_{L^1(\R)}.
\end{gather*}
Recall that $g_{q_{\ell}}, g_{\tilde q_{\ell}}$ are supported on the interval of length $\pi/\ell$ and that $\|g_{q_{\ell}}\|_{L^2(\R)}, \|g_{\tilde q_{\ell}}\|_{L^2(\R)}$ are uniformly bounded in $\ell$ by \eqref{eq: gr to g in L2}. We get
\begin{gather*}
\|g_{q_{\ell}} - g_{\tilde q_{\ell}}\|_{L^1(\R)}\le \sqrt{\pi/\ell}\|g_{q_{\ell}} - g_{\tilde q_{\ell}}\|_{L^2(\R)}
\le \sqrt{\pi/\ell}\big(\|g_{q_{\ell}}\|_{L^2(\R)} + \|g_{\tilde q_{\ell}}\|_{L^2(\R)}\big) = O(1/\sqrt{\ell}), \qquad \ell \to 0.
\end{gather*}
Therefore, the $k$-th term in the right-hand side of \eqref{eq: darboux sum v2} is $O(\ell^{3/2}e^{-A\ell k})$ as $\ell \to 0$ and the total sum can be estimated by
\begin{gather*}
O\bigg(\sum_{k\ge 0}\ell^{3/2}e^{-A\ell k}\bigg) = O\bigg(\frac{\ell^{3/2}}{1 - e^{-A\ell}}\bigg) = O(\ell^{1/2}),\qquad \ell \to 0.
\end{gather*}
In particular, it tends to $0$ as $\ell \to 0$, as was required in \eqref{eq: darboux sum}. The proof is concluded.
\qed\medskip

\subsection{Proof of Theorem~\ref{t1}} Recall that we need to prove 
\begin{gather}
\label{eq: theorem 2.1 expand}
\sqrt{\frac{\pi}{2}}\int_{\R_+}|q-\tilde q|e^{-2Ax}\,d x\le \int_{\R_+}|\hat f_q - \hat f_{\tilde q}|e^{-A\xi}\,d\xi \le 2\sqrt{2\pi}\int_{\R_+}|q-\tilde q|e^{-2Ax}\,dx,
\end{gather}
where $q, \tilde q\in L^2(\R_+)$ and $A\in \R$ satisfies $A \ge 12\max\big(\|q\|_{L^2(\R_+)}^2,  \|\tilde q\|^2_{L^2(\R_+)}\big)$.  
Without loss of generality, we can assume that 
\begin{gather}
\label{eq: strict assumption on A}
A > 12\max\big(\|q\|_{L^2(\R_+)}^2,  \|\tilde q\|^2_{L^2(\R_+)}\big),
\end{gather}
the claim for $A = 12\max\big(\|q\|_{L^2(\R_+)}^2,  \|\tilde q\|^2_{L^2(\R_+)}\big)$ will follow by a limiting argument (it worth be mentioned that the constant $12$ is not optimal). Fix $\eps > 0$. Let us show that the inequality
\begin{gather}
\label{eq: prev lemma 8.6}
\frac{1}{2(1 + \eps)}\sum_{k \ge 0}e^{-A\ell k}|q[\ell k]-\tilde q[\ell k]| 
\le 
\sum_{k\ge 0}e^{-A\ell k}|\hat F_{q_{\ell}}(k) - \hat F_{\tilde q_{\ell}}(k)|
\le
\frac{2}{1 - \eps}\sum_{k \ge 0}e^{-A\ell k}|q[\ell k]-\tilde q[\ell k]|
\end{gather}
holds for all sufficiently small $\ell$. Then \eqref{eq: theorem 2.1 expand} will follow if we take $\ell \to 0$ in \eqref{eq: prev lemma 8.6}, apply Lemma \ref{final lem 6.2} and Lemma \ref{final lem 6.11}, and send $\eps\to 0$. Thus, we can now focus on \eqref{eq: prev lemma 8.6}. Let $r = e^{-A\ell}$ and recall the definition \eqref{eq181} of the metric in $W^1(r\T)$. We have 
\begin{gather}
\label{eq: tmp 2}
\sum_{k\ge 0}e^{-A\ell k}|\hat F_{q_{\ell}}(k) - \hat F_{\tilde q_{\ell}}(k)| = \sum_{k\ge 0}r^k|\hat F_{q_{\ell}}(k) - \hat F_{\tilde q_{\ell}}(k)| =\|F_{q_{\ell}}- F_{\tilde q_{\ell}}\|_{W^1(r\T)}. 
\end{gather}
Note that $1 - r = 1 - e^{-A\ell} = A\ell + o(\ell)$ as $\ell\to 0$. Lemma \ref{lem: eta in limit} shows that
\begin{gather*}
\max\big(\|q\|_{L^2}^2, \|\tilde q\|^2_{L^2}\big)\ell= \max\big(\log\eta(F_{q_{\ell}})^{-1}, \log\eta(F_{\tilde q_{\ell}})^{-1}\big) + o(\ell),\quad \ell\to 0.
\end{gather*}
Then assumption \eqref{eq: strict assumption on A} for small $\ell$ implies
\begin{gather*}
1 - r \ge 12\max\big(\log\eta(F_{q_{\ell}})^{-1}, \log\eta(F_{\tilde q_{\ell}})^{-1}\big).
\end{gather*}
Therefore, Theorem \ref{t7} applies to Schur functions $F_{q_{\ell}}, F_{\tilde q_{\ell}}$ on the circle of radius $r =  e^{-A\ell}$ if $\ell$ is small enough. It gives
\begin{gather}
\label{eq: tmp 3}
\frac{1}{2}\rho_{B^1(r\T)}(F_{q_{\ell}}, F_{\tilde q_{\ell}}) \le
\|F_{q_{\ell}}- F_{\tilde q_{\ell}}\|_{W^1(r\T)}
\le 2\rho_{B^1(r\T)}(F_{q_{\ell}}, F_{\tilde q_{\ell}}).
\end{gather}
According to the definition \eqref{eq: SW metric def} of metric $\rho_{B^1(r\T)}$ and \eqref{eq: fqr and Fqr}, we have
\begin{gather}
\label{eq: norm on hT}
\rho_{B^1(r\T)}(F_{q_{\ell}}, F_{\tilde q_{\ell}})  = \sum_{k\ge 0}r^k|F_{q_{\ell},k}(0) -  F_{\tilde q_{\ell},k}(0)|. 
\end{gather}
Equality \eqref{eq: fqr and Fqr} and Theorem \ref{t2} state 
\begin{gather*}
q[\ell k] = \varkappa(F_{q_{\ell},k}(0)),\qquad \tilde
q[\ell k] = \varkappa(F_{\tilde q_{\ell},k}(0)).
\end{gather*}
The straightforward calculation shows that for $u,v\in \C$ small enough we have
\begin{gather}
\label{eq: kappa property}
|(\varkappa(u) - \varkappa(v)) + (\ov{u} - \ov{v})|\le \eps|u - v|.
\end{gather}
From \eqref{eq: sum fqrk(infty)} and \eqref{eq: fqr and Fqr} we know that  $\lim\limits_{\ell \to 0}F_{q_{\ell},k}(0) = \lim\limits_{\ell \to 0}F_{\tilde q_{\ell},k}(0) = 0$ hence \eqref{eq: kappa property} applies. It gives
\begin{align*}
\big|(q[\ell k]-\tilde q[\ell k]) + (\ov{F_{q_{\ell},k}(0)} - \ov{F_{\tilde q_{\ell}, k}(0)})\big|\le \eps\big|F_{q_{\ell},k}(0) - F_{\tilde q_{\ell}, k}(0)\big|,&
\\
(1 - \eps)\big|F_{q_{\ell},k}(0) - F_{\tilde q_{\ell}, k}(0)\big|\le \big|q[\ell k]-\tilde q[\ell k]\big| \le (1 + \eps)\big|F_{q_{\ell},k}(0) - F_{\tilde q_{\ell}, k}(0)\big|.&
\end{align*}
It follows from \eqref{eq: norm on hT} that  
\begin{gather}
\label{eq: tmp 1}
\frac{1}{1 +\eps}\sum_{k \ge 0}e^{-A\ell k}|q[\ell k]-\tilde q[\ell k]|\le \rho_{B^1(r\T)}(F_{q_{\ell}}, F_{\tilde q_{\ell}}) \le \frac{1}{1 -\eps}\sum_{k \ge 0}e^{-A\ell k}|q[\ell k]-\tilde q[\ell k]|.
\end{gather}
To establish \eqref{eq: prev lemma 8.6}, we substitute \eqref{eq: tmp 1} and \eqref{eq: tmp 2} into \eqref{eq: tmp 3}.\qed

\medskip

\section{Sylvester-Winebrenner theorem}\label{section: tSWD}
In the first part of this section we study the metric space $X$ and its subspace $S_{2}(\C_+)$ defined in \eqref{eq: spacex} and \eqref{eq: spaceS}, respectively. Then we use this analysis to prove that  the mapping $\F: q \mapsto f_{q}$ is a homeomorphism from $L^2(\R_+)$ onto $S_{2}(\C_+)$. The fact that $\F$ is a bijection  from $L^2(\R_+)$ to $S_{2}(\C_+)$ was proved in Theorem 11.4 in \cite{Den06}. Therefore, only continuity properties of this mapping need to be investigated. We check continuity of $\F$, $\F^{-1}$ in Section \ref{s52} and prove the lack of uniform continuity of $\F$, $\F^{-1}$ in Sections \ref{s53}, \ref{s54}, correspondingly. This last part implies Theorem \ref{t5}.
\subsection{Properties of the Sylvester-Winebrenner metric} 
Recall that the quantities $E(r_1, r_2)$, $E(r)$, for $r_1, r_2, r \in X$ are defined in \eqref{eq99}. Our aim here is to prove the following result.
\begin{Thm}\label{thm: convergence}
Let $r_n, r\in S_{2}(\C_+)$. Then $r_n \to r$ in $S_{2}(\C_+)$ if and only if $r_n \to r$ uniformly on compact subsets of $\C_+$ and $E(r_n) \to E(r)$.
\end{Thm}
The observation below was already used in the proof of Lemma \ref{final lem 6.11}, but we repeat it here for convenience.
\begin{Lem}\label{lem: conv in X implies conv in L2} 
For every $r_1, r_2\in X$ we have $2\rho_{X}(r_1, r_2)\ge \|r_1 - r_2\|_{L^2(\R)}$. 
\end{Lem}
\beginpf
Let us apply the inequality $-\log(1 - x)\ge x$ in formula \eqref{eq: metric in X} defining the metric $\rho_{X}$. We obtain
$$
\rho^2_{X}(r_1, r_2)  = \int_{\R}-\log\bigg(1 - \Big|\frac{r_1 - r_2}{1 - \ov{r_1}r_2}\Big|^2\bigg)\,dx\ge \int_{\R} \Big|\frac{r_1 - r_2}{1 - \ov{r_1}r_2}\Big|^2\,dx \ge \frac{1}{4}\int_{\R} \left|r_1 - r_2\right|^2\,dx.
\eqno\qed
$$
\begin{Lem}\label{lem: convergence mixed energy}
Let $r_n,r\in X$ be such that $r_n$ are bounded in $X$ and $r_n\to r$ as $n\to\infty$ in the Lebesgue measure on $\R$. Then $\lim_{n\to\infty}\limits  E(r_n, \phi)= E(r, \phi)$ for every $\phi\in X$.
\end{Lem}
\beginpf Considering $r_2 = 0$ in Lemma \ref{lem: conv in X implies conv in L2}, we see that $\{r_n\}_{n \ge 0}$ is a uniformly bounded sequence in~$L^2(\R)$. Since $r_n\to r$ in the Lebesgue measure on $\R$, we then have $r_n\to r$ weakly in $L^2(\R)$. Moreover, for each $k \ge 1$, $\{r_n^k\}_{n \ge 0}$ is also a uniformly bounded sequence in $L^2(\R)$ and $r_n^k \to r^k$ in the Lebesgue measure on $\R$. Thus, we have $r_n^k\to r^k$ weakly in $L^2(\R)$ for every $k \ge 1$.
Let us take $\phi \in X$ and represent $E(r_n, \phi)$ in the form
$$
E(r_n, \phi) = \Re\int_{\R}-\log(1 - \ov{r_n}\phi)\,dx = \Re\sum_{k = 1}^{N}\int_{\R}\frac{(\ov{r_n}\phi)^k}{k}\,dx + \Re\int_{\R}\sum_{k = N+1}^{\infty}\frac{(\ov{r_n}\phi)^k}{k}\,dx.
$$
Here and below $\log z$ denotes the analytic branch of the logarithm in $\C\setminus(-\infty, 0]$ such that $\log 1 = 0$. We are going to show that 
$$
\lim_{N\to \infty}\sup_n\left|\Re\int_{\R}\sum_{k = N+1}^{\infty}\frac{(\ov{r_n}\phi)^k}{k}\,dx\right| = 0.
$$
For this we write
\begin{align*}
\left|\int_{\R}\sum_{k = N+1}^{\infty}\frac{(\ov{r_n}\phi)^k}{k}\,dx\right|^2 
&\le \left(\int_{\R}|r_n\phi|^{N}\sum_{k = 1}^{\infty}\frac{|r_n\phi|^k}{k}\,dx\right)^2\\
&\le \left(\int_{\R} |r_n\phi|^N\sqrt{\log(1 - |r_n|^2)\log(1 - |\phi|^2)}\,dx\right)^2\\
&\le \int_{\R} \log(1 - |r_n|^2)\,dx \int_{\R} |\phi|^{2N}\log(1 - |\phi|^2)\,dx,
\end{align*} 
where we used Cauchy-Schwarz inequality (for sums and for integrals) and the fact that $|r_n| \le 1$ on $\R$. Since $\{r_n\}_{n \ge 0}$ is bounded in $X$, we have $\sup_n\int_{\R} |\log(1 - |r_n|^2)|\,dx = \sup_n\rho^2_X(0, r_n)< \infty$. Moreover, $\int_{\R} |\phi|^{2N}\log(1 - |\phi|^2)\,dx \to 0$ by the Lebesgue dominated convergence theorem with the majorant $|\log (1 - |\phi|^2)|$. Thus, we only need to check that 
$$
\lim_{N \to \infty}\lim_{n \to \infty}\Re\sum_{k = 1}^{N}\int_{\R}\frac{(\ov{r_n}\phi)^k}{k}\,dx = E(r, \phi).
$$
We have $\phi \in L^2(\R)$ by Lemma \ref{lem: conv in X implies conv in L2} for $r_1 = \phi$, $r_2 =0$. In view of $\|\phi\|_{L^\infty(\R)} \le 1$, this implies $\phi^k \in L^2(\R)$ for every $k\ge 1$. Then for each fixed $N \ge 1$ we obtain
$$
\lim_{n \to \infty}\Re\sum_{k = 1}^{N}\int_{\R}\frac{(\ov{r_n}\phi)^k}{k}\,dx = \Re\int_{\R}\sum_{k = 1}^{N}\frac{(\ov{r}\phi)^k}{k}\,dx
$$
from the weak convergence $r_n^k\to r^k$ in $L^2(\R)$. Moreover,  
$$
\left|\sum_{k = 1}^{N}\frac{(\ov{r}\phi)^k}{k}\right| \le \sum_{k = 1}^{N}\frac{|r\phi|^k}{k} \le \sqrt{\log(1 - |r|^2)\log(1 - |\phi|^2)},
$$ 
where the r.h.s.\ belongs to $L^1(\R)$ because $r, \phi \in X$. Thus, by the Lebesgue dominated convergence theorem we have 
$$
\lim_{N \to \infty}\Re\int_{\R}\sum_{k = 1}^{N}\frac{(\ov{r}\phi)^k}{k}\,dx =\Re\int_{\R}\sum_{k = 1}^{\infty}\frac{(\ov{r}\phi)^k}{k} = \int_{\R}-\log|1 - \ov{r}\phi|\,dx = E(r, \phi).
\eqno\qed
$$
\begin{Lem}\label{lem: conv in S2 and conv in measure v2}
Let $r_n,r\in X$. The following assertions are equivalent:
\begin{enumerate}
	\item[(a)] $r_n$ converges to $r$ in $X$\textup{;}
	\item[(b)] $\lim\limits_{n\to\infty}E(r_n) = E(r)$ and $r_n$ converges to $r$ in Lebesgue measure on $\R$\textup{;}
	\item[(c)] $\lim\limits_{n\to\infty}E(r_n) = E(r)$ and $\lim\limits_{n\to\infty}E(r_n, \phi) = E(r, \phi)$ for every $\phi\in X$.
	\item[(d)] $\lim\limits_{n\to\infty}E(r_n) = E(r)$ and $\lim\limits_{n\to\infty}E(r_n, \phi) = E(r, \phi)$ for every $\phi\in L^1(\R)$ with $\|\phi\|_{L^\infty(\R)} < 1$\textup{;}
	\item[(e)] $\lim\limits_{n\to\infty}E(r_n) = E(r)$ and $\lim\limits_{n\to\infty}E(r_n, r) = E(r)$\textup{;}
\end{enumerate}
\end{Lem}
\beginpf Since $E(r) = \rho^2_X(0, r)$ for every $r \in X$, we have $(a) \Rightarrow (b)$ by Lemma  \ref{lem: conv in X implies conv in L2}. Then, $(b) \Rightarrow (c)$ by Lemma \ref{lem: convergence mixed energy}. Clearly, $(c) \Rightarrow (d)$. To show that  $(d) \Rightarrow (a)$, we fix $\eps > 0$, take $r \in X$ and find $\phi\in L^1(\R)\cap L^{\infty}(\R)$ such that $\rho_X(\phi, r) < \eps$ (for instance, one can take $\phi = \chi_{F}r$, where $\chi_{F}$ is the characteristic function of the set $F=\{x: |r(x)| \le 1-1/N\} \cap [-N, N]$ for a sufficiently large $N$). Then, as $n \to \infty$, we have 
$$
\rho_X(r_n, r) \le \rho_X(r_n,\phi) + \rho_X(\phi, r) = \sqrt{E(r_n) + E(\phi) - 2E(r_n, \phi)} + \rho_X(\phi, r) \to 2\rho_X(\phi, r)
$$
by assumption $(d)$ and formula \eqref{eq: dist in terms of E}. Thus, $\rho_X(r_n, r) < 3\eps$ for all $n$ large enough. Since $\eps>0$ is arbitrary, this yields $(a)$. Equivalence of $(a)$ and $(e)$ is a simple consequence of \eqref{eq: dist in terms of E}. \qed

\medskip

\begin{Lem}\label{lem: conv on compacts imply conv in S2}
If $r_n, r\in S_2(\C_+)$ are such that $r_n \to r$ uniformly on compact subsets in $\C_+$, then 
\begin{gather}\label{eq: lem conv log}
	\lim\limits_{n\to\infty}E(r_n, \phi) = E(r, \phi)
\end{gather}
for every $\phi\in L^1(\R)\cap L^\infty(\R)$ with $\|\phi\|_{L^\infty(\R)} < 1$.
\end{Lem}
\beginpf As in the proof of Lemma \ref{lem: convergence mixed energy}, we have
$$
E(r_n, \phi) = \Re\int_{\R}-\log(1 - \ov{r_n}\phi)\,dx = \Re\sum_{k = 1}^{N}\int_{\R}\frac{(\ov{r_n}\phi)^k}{k}\,dx + \Re\int_{\R}\sum_{k = N+1}^{\infty}\frac{(\ov{r_n}\phi)^k}{k}\,dx,
$$
and 
$$
\lim_{N\to \infty}\sup_n\left|\Re\int_{\R}\sum_{k = N+1}^{\infty}\frac{(\ov{r_n}\phi)^k}{k}\,dx\right| \le  
\lim_{N\to \infty}\int_{\R}|\phi|^N\sum_{k = 1}^{\infty}\frac{|\phi|^k}{k}\,dx = \lim_{N\to \infty}\int_{\R}|\phi|^N \log(1-|\phi|)\,dx = 0,
$$
because $\log(1-|\phi|) \in L^{1}(\R)$. Thus, to prove that $E(r_n, \phi) \to E(r, \phi)$ it suffices to check that
\begin{align}
\label{eq: conv of terms in log}
&\lim_{n\to\infty}\Re\sum_{k = 1}^{N}\int_{\R}\frac{(\ov{r_n}\phi)^k}{k}\,dx = \Re\sum_{k = 1}^{N}\int_{\R}\frac{(\ov{r}\phi)^k}{k}\,dx,\\
&\lim_{N\to\infty}\Re\sum_{k = 1}^{N}\int_{\R}\frac{(\ov{r}\phi)^k}{k}\,dx = \int_{\R}\Re\sum_{k = 1}^{\infty}\frac{(\ov{r}\phi)^k}{k}\,dx =E(r, \phi).\label{eq102}
\end{align}
As in the proof of Lemma \ref{lem: convergence mixed energy}, relation \eqref{eq102} is a consequence of the Lebesgue dominated convergence theorem (this time -- with the majorant $|\log(1-|\phi|)|$). So, we can focus on the proof of \eqref{eq: conv of terms in log}. Since $|\bar r_n^k|\le 1$ and $\phi^k \in L^1(\R)$ for all $n$, $k$, it suffices to prove that 
\begin{equation}\label{eq1031}
\lim_{n\to\infty}\int_{\R}\ov{r_n}^{k}\psi\,dx = \int_{\R}\ov{r}^{k}\psi\,dx
\end{equation}
for functions $\psi$ from a dense subset of $L^1(\R)$. By our assumption, $\lim_{n\to\infty} r_n^k(z) = r^k(z)$ for every $k \ge 1$ and $z\in \C_+$. Thus, \eqref{eq1031} holds for $\psi$ from the set of all finite linear combinations of Poisson kernels $P_z: x \mapsto \frac{1}{\pi}\Im \frac{1}{x - z}$, $z\in \C_+$. This set is dense in $L^1(\R)$, which completes the proof. \qed 

\medskip

\noindent {\bf Proof of Theorem \ref{thm: convergence}}. Let $r_n, r\in S_{2}(\C_+)$ be such that $r_n \to r$ in $S_{2}(\C_+)$. Then 
$E(r_n) \to E(r)$ by Lemma \ref{lem: conv in S2 and conv in measure v2}. Moreover, we have $r_n \to r$ in the Hardy space $H^2(\C_+)$ by Lemma 
\ref{lem: conv in X implies conv in L2}.  In particular, $r_n \to r$ uniformly on compact subsets of $\C_+$. Conversely, assume that 
$r_n \to r$ uniformly on compact subsets of $\C_+$ and $E(r_n) \to E(r)$. Then Lemma \ref{lem: conv on compacts imply conv in S2} and Lemma \ref{lem: conv in S2 and conv in measure v2} imply $r_n \to r$ in $S_{2}(\C_+)$. \qed

\medskip

\subsection{Spectral map is a homeomorphism}\label{s52} As we mentioned at the beginning of Section \ref{section: tSWD}, the fact that the mapping $\F: q \mapsto f_{q}$ is a bijection  from $L^2(\R_+)$ to $S_{2}(\C_+)$ is proved in Theorem 11.4 in \cite{Den06}; formula \eqref{eqSW} is contained in Corollary 11.2 in the same paper \cite{Den06}. The proof of continuity in Theorem \ref{tSW} relies on the following known lemma.
\begin{Lem}\label{lem: krein de branges lemma}
Suppose that $\{q_n\}$ is a bounded sequence in $L^2(\R_+)$, and let $\{f_{q_n}\}$ be the corresponding sequence of Schur functions. Then, $\{q_n\}$ converges weakly in $L^2(\R_+)$ if and only if $\{f_{q_n}\}$ converges uniformly on compact subsets of $\C_+$.    
\end{Lem}
\beginpf Suppose that a sequence $\{q_n\}$ converges weakly to $q \in L^2(\R_+)$. Then we have
$$
\lim_{n \to \infty} \max_{0\le y\le L}\left|\int_{0}^{y}q_n(x)\,dx -\int_{0}^{y} q(x)\,dx\right| = 0
$$
for every $L \ge 0$. It follows that the corresponding Weyl functions, $\{m_{q_n}\}$, converge on compact subsets on $\C_+$ to the Weyl function $m_{q}$ of $q$. To prove this fact it suffices to note that the solutions of equations
$$
JN_{q_n}'(x) + Q_n(x)N_{q_n}(x) = 0, \qquad N_{q_n}(0) = \IDM, \qquad Q_n = \begin{pmatrix} \Im q_n & \Re q_n \\ \Re q_n & -\Im q_n \end{pmatrix}
$$
converge to the solution $N_{q}$ of the corresponding equation for $q$ uniformly on $[0, L]$ for each $L \ge 0$ (see \eqref{eq: series representation}, \eqref{eq: series bound}) and use Theorem \ref{KdB} for Hamiltonians $\Hh_{q_n} = N_{q_n}^*N_{q_n}$, $\Hh_{q} = N_{q}^*N_{q}$. The locally uniform convergence $m_{q_n} \to m_{q}$ on $\C_+$ implies the locally uniform convergence $f_{q_n} \to f_{q}$.

\medskip

Conversely, assume that $\{q_n\}$ is a bounded sequence in $L^2(\R_+)$ such that Schur functions $f_{q_n}$ converge on compact subsets of $\C_+$ to some function $f$. Using the weak compactness of closed bounded subsets of $L^2(\R_+)$, we see that from each subsequence $\{q_{n_k}\}$ one can extract another subsequence $\{q_{n_{k_j}}\}$ weakly converging to a function $q \in L^2(\R_+)$. From the first part of the proof we get $f_q = f$. In particular, $q$ is determined uniquely by $f$ (the spectral correspondence $q \mapsto f_q$ is a bijection, see Theorem 11.4 in \cite{Den06}) and the weak limit of $\{q_{n_k}\}$ equals $q$ for all subsequences of the sequence $\{q_{n}\}$. It follows that $\{q_{n}\}$ weakly converges to $q$. \qed

\medskip

\noindent\textbf{Proof of continuity in Theorem \ref{tSW}.}
Let us prove that the map $\F: L^2(\R_+) \to S_2(\C_+)$ is a homeomorphism assuming that this map is a bijection and the sum rule \eqref{eqSW} holds for all $q \in L^2(\R_+)$. Consider some  potentials $q_n, q\in L^2(\R_+)$  and let $f_n, f \in S_{2}(\C_+)$ be the corresponding Schur functions. It is well-known that $q_n \to q$ in $L^2(\R_+)$ if and only if $\|q_n\|_{L^2(\R_+)} \to \|q\|_{L^2(\R_+)}$ and $q_n \to q$ weakly in $L^2(\R_+)$.  In view of \eqref{eqSW}, we have $\|q_n\|_{L^2(\R_+)} \to \|q\|_{L^2(\R_+)}$ if and only if $E(f_n) \to E(f)$.  In particular, the sequence $\{q_n\}$ is bounded in $L^2(\R_+)$ if and only if the sequence $\{f_n\}$ is bounded in $S_2(\C_+)$. If one of these sequences is bounded, we know from Lemma \ref{lem: krein de branges lemma} that $q_n \to q$ weakly in $L^2(\R_+)$ if and only if $f_n \to f$ on compact subsets of $\C_+$. Thus, $q_n \to q$ in $L^2(\R_+)$ if and only if $E(f_n) \to E(f)$ and $f_n \to f$ on compact subsets of $\C_+$. By Theorem \ref{thm: convergence} this is equivalent to the convergence $f_n \to f$ in $S_2(\C_+)$. This proves that the map $\F: L^2(\R_+) \to S_2(\C_+)$ is a homeomorphism. \qed

\medskip

\subsection{Direct map is not uniformly continuous}\label{s53} Our next aim is to prove the following proposition. 
\begin{Prop}\label{p57}
There are potentials $u_n$, $\tilde u_n$ in the unit ball of $L^2(\R_+)$ such that \eqref{eq12} holds. In other words, the direct homeomorphism $\F: S_{2}(\C_+) \to L^2(\R_+)$ in Theorem \ref{tSW} is not uniformly continuous on bounded subsets of $S_{2}(\C_+)$. 
\end{Prop}
\beginpf We construct an explicit example. For $T \ge 1$, $v \ge 0$, denote by $q_{T, v}$ the constant imaginary-valued potential on $[0,T]$ with positive imaginary part that satisfies $\|q_{T, v}\|_{L^2[0, T]} = v$. Extend it by zero to $(T, \infty)$. Also let $\eps(T)$ be some positive number that satisfy $100e^{-\sqrt{T}/6} \le \eps(T) \le T^{-1}$. We are going to show that
\begin{gather}
\label{eq: claim for constant potentials}
\|q_{T, 1} - q_{T, 1 - \eps(T)}\|_{L^2(\R_+)} = \eps(T), \qquad \liminf_{T\to \infty}\rho_{S_{2}}(f_{q_{T, 1}}, f_{_{T, 1 - \eps(T)}}) > 0.
\end{gather}
The claim will then follow. Note that the first relation in \eqref{eq: claim for constant potentials} is a simple consequence of the definition of $q_{T, 1}$, $q_{T, 1 - \eps(T)}$, and we only need to check the second relation. To fix notation, take $p > 0$ and let $q\in L^2(\R_+)$ be the piece-wise constant potential such that $q = ip$ on $[0, T]$, $q = 0$ on $(T,\infty)$. Then $Q(t) = \left(\begin{smallmatrix}	p & 0 \\	0 & -p \end{smallmatrix}\right)$ on $[0, T]$ and the fundamental solution $N(t,z)$ of the corresponding Dirac system $JN' + QN = zN$, $N(0, z)=\idm$, for $0 \le t \le T$, $z \in \C$, is given by 
$$
N(t, z) =\exp(tJ(Q - z)) = \exp\left(t\left(\begin{smallmatrix}	0 & p + z \\	p- z & 0 \end{smallmatrix}\right)\right),
$$
that can be checked by the differentiation with respect to $t$. Then, recall Lemma \eqref{l38},
$$
N(T, z) = 
\begin{pmatrix}
\cosh (T\lambda) & \frac{p + z}{\lambda}\sinh(T\lambda)\\
\frac{\lambda}{p + z}\sinh(T\lambda) & \cosh(T\lambda)
\end{pmatrix}, \qquad \lambda = \sqrt{p^2 - z^2}.
$$
Notice that $N(T,z)$ is an entire function in $z$ and its does not depend on the choice of the root $\lambda = \sqrt{p^2 - z^2}$. Indeed, this follows from the fact that functions $\cosh(T\lambda)$, $\sinh(T\lambda)/\lambda$ depend only on the value $\lambda^2$. Further, for $t > T$ we have $Q(t) = \left(\begin{smallmatrix}	0 & 0 \\	0 & 0 \end{smallmatrix}\right)$, hence
\begin{gather*}
N(t,z) = \exp(-(t - T)zJ)N(T, z), \qquad t \ge T.
\end{gather*}
The fundamental matrix solution of the canonical system with the Hamiltonian $\Hh_q = N(t, 0)^*N(t, 0)$ that corresponds to the Dirac operator with the potential $q$ is given by $M(t,z) = N(t, 0)^{-1}N(t, z)$, see, e.g., Section 2.4 in \cite{B2020}.  Therefore, the Weyl function of this canonical system equals
$$
m_{\Hh_q} \doteq \lim_{t \to \infty}\sigma_1 M(t,z)^{T} \sigma_1 \omega \doteq \lim_{t \to \infty}\sigma_1 N(t,z)^{T} (N(t, 0)^{-1})^T \sigma_1\omega
\doteq \lim_{t \to \infty}\sigma_1 N(t,z)^{T} \sigma_1\tilde\omega(t),
$$
where $\tilde\omega(t) = \sigma_1(N(t, 0)^{-1})^T \sigma_1 \omega$ belongs to $\C_+$ because $\omega\in \C_+$ and $\sigma_1(N(t, 0)^{-1})^T \sigma_1 \in \sltr$. So, we have
\begin{align*}
\lim_{t \to \infty}\sigma_1 N(t,z)^{T} \sigma_1\tilde\omega(t) 
&\doteq \sigma_1 N(T,z)^{T} \sigma_1 \lim_{t \to \infty}\sigma_1 \exp(-(t - T)zJ)^{T}  \sigma_1\tilde\omega(t)  \\
&\doteq \sigma_1 N(T,z)^{T} \sigma_1 i \doteq N(T,z) i,
\end{align*}
where we used the fact that the fractional-linear transformation with the matrix $\sigma_1 \exp(-szJ)^{T}  \sigma_1$ for large $s > 0$ maps $\C_+$ into a small disk with center at $i$. Then,
\begin{align*}
f_q 
&= \frac{m_q - i}{m_q + i} = \frac{m_{\Hh_q} - i}{m_{\Hh_q} + i}  \doteq \begin{pmatrix}1 &-i \\ 1 &i \end{pmatrix} N(T,z) i 
\doteq \begin{pmatrix}1 &-i \\ 1 &i \end{pmatrix} 
\begin{pmatrix} i\cosh(T\lambda) + \frac{p + z}{\lambda}\sinh(T\lambda)\\i\frac{\lambda}{p + z}\sinh(T\lambda) + \cosh(T\lambda) \end{pmatrix} \\
&=
\frac{\left(\frac{p + z}{\lambda} + \frac{\lambda}{p + z}\right)\sinh(T\lambda)}{2i\cosh(T\lambda) + \left(\frac{p + z}{\lambda} - \frac{\lambda}{p + z}\right)\sinh(T\lambda)} = \frac{\frac{2p}{\lambda}\sinh(T\lambda)}{2i\cosh(T\lambda) + \frac{2z}{\lambda}\sinh(T\lambda)}\\
&=\frac{-ip\sinh(T\lambda) }{\lambda\cosh(T\lambda) - iz \sinh(T\lambda)}.
\end{align*}
From now on we will assume that $z = x$ with real $x\in [c_1p, c_2p]$, where $c_1 = \sqrt{3/4}$, $c_2 = \sqrt{8/9}$. In this case $p^2 - x^2 > 0$, hence $\lambda$ can be chosen positive, $\lambda\in [p/3, p/2]$. For brevity, we will write $f = f_q$ and $\eps = \eps(T)$ in the remaining part of the proof. We have
\begin{gather*}
f(x) = \frac{-ip(e^{T\lambda} - e^{-T\lambda})}{\lambda(e^{T\lambda} + e^{-T\lambda}) - ix (e^{T\lambda} - e^{-T\lambda})} = \frac{-ip + ipe^{-2T\lambda}}{(\lambda - ix) + (\lambda + ix)e^{-2T\lambda}} =g(x) + e^{-2T\lambda}h(x),
\\
g(x) = \frac{-ip}{\lambda - ix},\qquad
h(x) = e^{2T\lambda}(f_q(x) - g(x)) = \frac{-2ip\lambda}{(x+i\lambda)^2 - p^2 e^{-2T\lambda}}.
\end{gather*}
Let us denote $q = q_{T, 1}$ and $q_{\eps} = q_{T, 1 - \eps}$, and let $f_{\eps}, g_\eps, h_\eps$ be the functions corresponding to $q_{\eps}$. We have
\begin{align*}
|g(x) - g_{\eps}(x)| 
&= \left|\frac{p}{\lambda - ix} -\frac{p_\eps}{\lambda_\eps - ix} \right| = \frac{|p\lambda_\eps - p_\eps\lambda + ix(p_\eps - p)|}{|(\lambda - ix)(\lambda_\eps - ix)|}\ge \frac{|p\lambda_\eps - p_\eps\lambda|}{pp_{\eps}} 
\\
&= \left|\lambda_\eps /p_\eps - \lambda/p\right| =  \left|\sqrt{1 - {x^2}/{p_\eps^2}} - \sqrt{1 - {x^2}/{p^2}}\right| 
= \frac{{x^2}/{p^2} - {x^2}/{p_\eps^2}}{\sqrt{1 - {x^2}/{p_\eps^2}} + \sqrt{1 - {x^2}/{p^2}}}\\
&\ge \frac{{x^2}/{p_\eps^2} - {x^2}/{p^2}}{2} \ge \frac{x^2}{p^2}(1 - p_{\eps}^2/p^2)\ge c_1^2(1 - p_{\eps}^2/p^2).
\end{align*}
Recall that $p = 1/\sqrt{T}$, $p_\eps = (1-\eps)/\sqrt{T}$. Hence $1 - p_{\eps}^2/p^2 = 1 
- (1 - \eps)^2 \ge \eps$ for small $\eps$ (equivalently, for large $T$). It follows that 
\begin{gather}
\label{eq: g - geps}
|g(x) - g_{\eps}(x)| \ge \frac{c_1^2\eps}{2} = \frac{3\eps}{8},\qquad x\in [c_1p,c_2p].
\end{gather}
Furthermore, we  have $\lambda\ge\lambda_\eps \ge p/4$ for small $\eps$,  hence $e^{-2T\lambda}\le e^{-2T\lambda_\eps}\le e^{-2Tp/4} = e^{-\sqrt{T}/2}$ and for $T$ large enough we have
\begin{gather}
\label{eq: f - g bound}
|f(x) - g(x)|= e^{-2T\lambda}|h(x)|= 
\frac{2e^{-2T\lambda}p\lambda}{|(x+i\lambda)^2 - p^2 e^{-2T\lambda}|}\le\frac{2e^{-2T\lambda}p\lambda}{p^2(1 - e^{-2T\lambda})}\le e^{-\sqrt{T}/2},
\\
\label{eq: feps - geps bound}
|f_{\eps}(x) - g_{\eps}(x)|= e^{-2T\lambda_\eps}|h_\eps(x)| = \frac{2e^{-2T\lambda_\eps}p_\eps\lambda_\eps}{|(x+i\lambda_\eps)^2 - p_\eps^2 e^{-2T\lambda_\eps}|}\le\frac{2e^{-2T\lambda_\eps}p_\eps\lambda_\eps}{p_\eps^2(1 - e^{-2T\lambda_\eps})}\le e^{-\sqrt{T}/2}.
\end{gather}
The triangle inequality implies
\begin{align}
\nonumber
|f(x) - f_{\eps}(x)|&\ge |g(x) - g_\eps(x)| - |f_q(x) - g(x)| - |g_\eps(x) - f_{\eps}(x)|
\\
\label{eq: f - feps}
&\ge \frac{3}{8}\eps - e^{-2T\lambda}|h(x)| - e^{-2T\lambda_\eps}|h_\eps(x)| \ge \frac{3}{8}\eps - 2e^{-\sqrt{T}/2} \ge  \frac{\eps}{4}.
\end{align}
By definition of $\lambda$, we have
\begin{gather*}
|g(x)|^2 = \left|\frac{-ip}{\lambda - ix}\right|^2 = \frac{p^2}{\lambda^2 + x^2} = 1,\qquad \left|\frac{g - g_{\eps}}{1 - \ov{g}g_{\eps}}\right|= \left|\frac{g - g_{\eps}}{\ov{g}g - \ov{g}g_{\eps}}\right| = 1.
\end{gather*}
We are ready to estimate $\rho_{S_{2}}(f, f_{\eps})$. We have
\begin{gather*}
\left|\frac{f - f_{\eps}}{1 - \ov{f}f_{\eps}} - \frac{g - g_{\eps}}{1 - \ov{g}g_{\eps}}\right| = \left|\frac{(f - f_{\eps})(1 - \ov{g}g_{\eps}) - (g - g_{\eps})(1 - \ov{f}f_{\eps})}{(1 - \ov{g}g_{\eps})(1 - \ov{f}f_{\eps})}\right| \le \frac{4(|f - g| + |f_{\eps} - g_{\eps}|)}{\left|(1 - \ov{g}g_{\eps})(1 - \ov{f}f_{\eps})\right|}.
\end{gather*}
For the denominator we write $|1 - \ov{g}g_{\eps}| = |g - g_{\eps}|$, $|1 - \ov{f}f_{\eps}| \ge |f - f_{\eps}|$ and apply \eqref{eq: g - geps}, \eqref{eq: f - feps}; for the numerator we use \eqref{eq: f - g bound} and \eqref{eq: feps - geps bound}. This gives
\begin{gather*}
1- \left|\frac{f - f_{\eps}}{1 - \ov{f}f_{\eps}} \right|\le \left|\frac{f - f_{\eps}}{1 - \ov{f}f_{\eps}} - \frac{g - g_{\eps}}{1 - \ov{g}g_{\eps}}\right|\le \frac{4\cdot 2\cdot e^{-\sqrt{T}/2}}{3\eps/8\cdot \eps/4} < \frac{2^8 e^{-\sqrt{T}/2}}{3\eps^2} \le \frac{e^{-\sqrt{T}/6}}{2}.
\end{gather*}
It follows that
$$
1- \left|\frac{f - f_{\eps}}{1 - \ov{f}f_{\eps}} \right|^2\le e^{-\sqrt{T}/6},\qquad -\log\bigg(1- \left|\frac{f - f_{\eps}}{1 - \ov{f}f_{\eps}} \right|^2\bigg)\ge \frac{\sqrt{T}}{6}.
$$
Finally, we write
\begin{gather*}
\rho_{S_{2}}(f, f_\eps)^2\ge \int_{c_1p}^{c_2p}-\log\bigg(1 - \left|\frac{f - f_{\eps}}{1 - \ov{f}f_{\eps}}\right|^2\bigg)\, dx\ge \frac{(c_2 - c_1)p\sqrt{T}}{6} = \frac{c_2 - c_1}{6} > 0.
\end{gather*}
The proof is concluded. \qed

\medskip

\begin{Rema} The argument can be adapted to prove that $\F$ is not uniformly continuous on any ball of radius $\delta > 0$ with center at $0$. 
\end{Rema}

\subsection{Inverse map is not uniformly continuous}\label{s54}
In this section we prove that the mapping $\F^{-1}\colon S_{2}(\C_+)\to L^2(\R_+)$ from Theorem \ref{tSW} is not uniformly continuous. The proof is based on the non-injectivity of the direct scattering transform of the Dirac operator on $\R$. Let us introduce the required objects, following notation from \cite{BD24}. The results of the present section are adaptations of the results from Section 6 in \cite{BG2024}, where the discrete scattering transform was considered. 
\medskip

\noindent We start with the function $q\in L^2(\R)$ and introduce $q^+,q^-\in L^2(\R_+)$,
\begin{gather*}
q^+(x) = -\ov{q(x)},\qquad q^-(x) = q(-x),\qquad x\ge 0.
\end{gather*}
The functions $f^+$ and $f^-$ are the Schur functions corresponding to $q^+$ and $q^-$ respectively as described in the introduction of the present paper, see \eqref{sf}. From Theorem \ref{tSW} we know that $\log(1 - |f^\pm|^2)\in L^1(\R)$ hence there exist outer in $\C_+$ functions $\fa^\pm$ satisfying $|\fa^\pm|^{-2} = 1 - |f^\pm|^2$, see Theorem 4.4 in \cite{Garnett}; the analytic functions $\fb^\pm$ are defined so that $f^\pm = \fb^\pm / \fa^\pm$ in $\C_+$. Finally, we let
\begin{gather}
\label{eq: refl def}
a = \fa^+\fa^- - \fb^+\fb^-,\qquad b = \fa^-\ov{\fb^+} - \fb^-\ov{\fa^+},\qquad \rc_q = b/a.
\end{gather}
The function $a$ is outer in $\C_+$, and
\begin{gather}
\label{eq: modulus of a b r}
|a|^2 = 1 + |b|^2, \qquad 1-|\rc_q|^2 = 1 - \frac{|b|^2}{|a|^2} = \frac{|a|^2 - |b|^2}{|a|^2} = |a|^{-2}.
\end{gather}
The function $\rc_q$ is called the reflection coefficient of $q$ and the mapping $q\mapsto \rc_q$ acting from $L^2(\R)$ to the space $X$ defined in \eqref{eq: spacex} is called the direct scattering transform. It has some symmetries, see Lemma 2.1 in \cite{BD24}. We will need a part of this lemma that concerns translation of the potential.
\begin{Lem}[Lemma 2.1 in \cite{BD24}]\label{lem: shift reflection coefficient}
Let $q\in L^2(\R)$. We have $\rc_{q_s}(\lambda) = e^{-i\lambda s}\rc_{q}(\lambda)$ for almost every $\lambda\in\R$, where $s\in \R$ and $q_s\colon x\mapsto q(x - s)\in L^2(\R)$ is the $s$-shift of $q$. 
\end{Lem}

The following proposition is also well-known.

\begin{Prop}\label{prop: refl coeff and potent}
For every $q\in L^2(\R)$, we have
\begin{gather*}
	\int_\R |q|^2\,dx = \frac{1}{\pi}\int_\R-\log (1 - |\rc_q|^2)\,dx.
\end{gather*}
\end{Prop}
\beginpf
In the upper half-plane $\C_+$ we can write 
\begin{gather*}
a = \fa^+\fa^- - \fb^+\fb^- = \fa^+\fa^-\left(1 - \frac{\fb^+}{\fa^+}\cdot \frac{\fb^-}{\fa^-}\right) = \fa^+\fa^-\left(1 - f^+f^-\right).
\end{gather*}
Therefore, \eqref{eq: modulus of a b r} implies
\begin{gather*}
\int_\R\log (1 - |\rc_q|^2)\,dx = -\int_\R\log |\fa^+|^2\,dx - \int_\R\log |\fa^-|^2\,dx 	 -\int_\R\log |1 - f^+f^-|^2\,dx.
\end{gather*}
Recall that $\fa^\pm$ are defined so that $|\fa^\pm|^{-2} = 1 - |f^\pm|^2$ a.\,e.\,on $\R$. Theorem \ref{tSW} gives
\begin{gather*}
-\int_\R\log |\fa^+|^2\,dx = \int_\R\log(1 -  |f^+|^2)\,dx = -{\pi}\int_{\R_+}|q^+|^2\,dx,
\\
-\int_\R\log |\fa^-|^2\,dx = \int_\R\log(1 -  |f^-|^2)\,dx = -{\pi}\int_{\R_+}|q^-|^2\,dx,
\\
\int_\R\log |\fa^+|^2\,dx + \int_\R\log |\fa^-|^2\,dx = \pi\int_{\R_+}|q^+|^2\,dx + \pi\int_{\R_+}|q^-|^2\,dx =  \pi\int_{\R}|q|^2\,dx.
\end{gather*}
We have $f^\pm\in L^2(\R)\cap H^\infty(\C_+)$. It follows that $f^\pm\in H^2(\C_+)$, see Corollary 4.3
in \cite{Garnett} and consequently $f^+f^-\in H^1(\C_+)$. We also know that $|f^+f^-|\le 1$ in $\C_+$ hence $\log(1 - f^+f^-)$ is well-defined analytic function in $\C_+$ that also belongs to $H^1(\C_+)$. Therefore we have 
\begin{gather*}
\int_\R\log |1 - f^+f^-|^2\,dx = 2\Re\int_\R\log (1 - f^+f^-)\,dx =0,
\end{gather*}
where we used the fact that the integral of an $H^1(\C_+)$ function over $\R$ is $0$, see Lemma 3.7 in \cite{Garnett}. 
\qed

\begin{Prop}\label{prop: continuity from L^2 to S_2}
If $q_n\to q$ in $L^2(\R)$ as $n\to\infty$ then $\rc_{q_n}\to \rc_{q}$ in $X$ as $n\to\infty$.
\end{Prop}
\beginpf
We have $\|q_n^+ - q^+\|_{L^2(\R_+)}\to 0$ and $\|q_n^- - q^-\|_{L^2(\R_+)}\to 0$ as $n\to\infty$. Then by Theorem \ref{tSW} we get $f_n^+\to f^+$ and $f_n^-\to f^-$ in $S_{2}(\C_+)$ as $n\to\infty$. In particular, $f_n^\pm\to f^\pm$ in measure on $\R$ by Lemma \ref{lem: conv in S2 and conv in measure v2}. Therefore, $\log(1 - |f_n^+|^2)$ converge in measure to $\log(1 - |f^+|^2)$.  Theorem  \ref{tSW} also gives
\begin{gather*}
\lim_{n\to\infty}\|\log(1 - |f_n^+|^2)\|_{L^1(\R)} = \lim_{n\to\infty}{{\pi}}\|q_n^+\|^2_{L^2(\R_+)}= {{\pi}}\|q^+\|^2_{L^2(\R_+)} = \|\log(1 - |f^+|^2)\|_{L^1(\R)}.
\end{gather*}
Then Scheff\'{e}'s lemma, see Lemma 5.10 in \cite{williams1991probability}, implies that $\log(1 - |f_n^+|^2)$ converge  to $\log(1 - |f^+|^2)$ in $L^1(\R)$. Now we write 
\begin{align}
\nonumber
\rc_{q} 
&= \frac{\fa^-\ov{\fb^+} - \fb^-\ov{\fa^+}}{\fa^+ \fa^- - \fb^+\fb^-} = \frac{\ov{\fa^{+}}}{\fa^{+}}\frac{\ov{f^+} - f^{-}}{1 - f^+f^-}\notag = \exp(-2iH(\log|\fa^+|))\frac{\ov{f^+} - f^{-}}{1 - f^+f^-}
\\
\label{eq: r in terms of f}
&=
\exp(iH(\log(1-|f^+|^2)))\frac{\ov{f^+} - f^{-}}{1 - f^+f^-}.
\end{align}
Here $H$ denotes the Hilbert transform and we used the fact that $\fa^+$ is an outer function. Weak continuity of the Hilbert transform $H$ on $L^1(\R)$ (see Section III.2 in \cite{Garnett}) then implies that $\exp(iH(\log(1-|f_n^+|^2)))$ converges in the Lebesgue measure to $\exp(iH(\log(1-|f^+|^2)))$. From here and \eqref{eq: r in terms of f} we see that functions $\rc_{q_n}$ converge to $\rc_{q}$ in Lebesgue measure on $\R$. By Proposition \ref{prop: refl coeff and potent} we also have
\begin{gather*}
\lim_{n\to\infty}E(\rc_{q_n})   = \lim_{n\to\infty}{\pi}\int_\R |q_n|^2\,dx ={\pi}\int_\R |q|^2\,dx = E(\rc_{q}).
\end{gather*}
The convergence in $X$ now follows from  Lemma \ref{lem: conv in S2 and conv in measure v2}. \qed

\medskip

Next proposition is analogous to the noninjectivity of the scattering transform for Jacobi matrices discovered by Volberg and Yuditskii in \cite{VYu02}.

\begin{Prop}[Example 6.1 in \cite{BD24}]\label{prop: not injective}
There exist potentials $q, \tilde q\in L^2(\R)$ such that $q\neq\tilde q$ in $L^2(\R)$ with $\rc_q = \rc_{\tilde q}$ a.e.\,in $\R$.
\end{Prop}
\begin{Lem}\label{lem: change of q to -q}
Let $q\in L^2(\R_+)$ be supported on $(-\infty, 0]$. Then $\rc_q = -f^-$. Moreover, $\rc_q = f_{q^\#}$, where $q^\#\in L^2(\R_+)$ is such that $q^\#(x) = -q^-(x) = q(-x)$ for $x\ge 0$.
\end{Lem}
\beginpf
The potential $q$ is supported on $(-\infty, 0]$ hence $q^+ = 0$ and
\begin{gather*}
f^+ = 0,\quad \fa^+ = 1,\quad \fb^+ = 0, \qquad a = \fa^-,\quad b = -\fb^-,\qquad \rc_{q} = \frac{b}{a} = \frac{-\fb^-}{\fa^-} = -f^-.
\end{gather*}
We also have $f_{q^\#} = - f_{q^-}$, see Section 7 in \cite{Den06}.
The equality $\rc_q = f_{q^\#}$ follows. \qed
\begin{Prop}\label{p514}
There are potentials $q_n$, $\tilde q_n$ in the unit ball of $L^2(\R_+)$ such that \eqref{eq11} holds. In other words, the inverse homeomorphism $\F^{-1}: S_{2}(\C_+) \to L^2(\R_+)$ in Theorem \ref{tSW} is not uniformly continuous on bounded subsets of $S_{2}(\C_+)$. 
\end{Prop}
\beginpf
Let $q$ and $\tilde q$ be the potentials from Proposition \ref{prop: not injective}, denote $\rc = \rc_q = \rc_{\tilde q}$. For any $\delta > 0$, considerations similar to those in Example 6.1 in \cite{BD24} allow to assume without loss of generality that $E(\rc) \le \delta$. In particular, we can assume without loss of generality that potentials $q_n$, $\tilde q_n$ lie in the unit ball of $L^2(\R_+)$. For  $N\in \R$, let 
\begin{gather*}
q_N(x) = 
\begin{cases}
	q(x), &x \le N	,\\
	0, &x  > N,
\end{cases}
\qquad 
q_{N,s}(x) = q_N(x + N),
\\
\tilde q_N(x) = 
\begin{cases}
	\tilde q(x), &x \le N,
	\\
	0, &x  > N,
\end{cases}
\qquad 
\tilde q_{N,s}(x) = \tilde q_N(x + N).
\end{gather*}
We have $q_N\to q$ and $\tilde q_N\to\tilde q$ in $L^2(\R)$ hence for large $N$
\begin{gather}
\label{eq: shifted potentials are separated}
\|q_{N,s}- \tilde q_{N,s}\|_{L^2(\R)} = \|q_{N}- \tilde q_{N}\|_{L^2(\R)}\ge \|q- \tilde q\|_{L^2(\R)}/2.
\end{gather}
By Proposition \ref{prop: continuity from L^2 to S_2}, we have $\rc_{q_N}\to \rc$ and $\rc_{\tilde q_N}\to \rc$ in $X$ as $n\to\infty$. Therefore, for every $\eps > 0$ we can choose large $N$ such that 
\begin{gather*}
\rho_{X}(\rc_{q_N}, \rc)\le \eps/2, \quad \rho_{X}(\rc, \rc_{\tilde q_N})\le \eps/2,\qquad \rho_{X}(\rc_{q_N}, \rc_{\tilde q_N})\le \rho_{X}(\rc_{q_N}, \rc) + \rho_{X}(\rc, \rc_{\tilde q_N})\le \eps.
\end{gather*}
Potentials $q_{N,s}$ and $\tilde q_{N,s}$ are the $N$-shifts of  $q_{N}$ and $\tilde q_{N}$ respectively, hence Lemma \ref{lem: shift reflection coefficient} applies. It gives
\begin{gather*}
\rc_{q_{N,s}}(x) = e^{iNx}\rc_{q_{N}}(x),\qquad \rc_{\tilde q_{N,s}}(x) = e^{iNx}\rc_{\tilde q_{N}}(x),\qquad \rho_{X}(\rc_{q_{N,s}}, \rc_{\tilde q_{N,s}}) = \rho_{X}(\rc_{q_N}, \rc_{\tilde q_N})\le \eps.
\end{gather*}
Both $q_{N,s} $ and $\tilde q_{N,s}$ are supported on $(-\infty, 0]$ hence $\rc_{q_{N,s}} = f_{q_{N,s}^\#}$ and $\rc_{\tilde q_{N,s}} = f_{\tilde q_{N,s}^\#}$ by Lemma \ref{lem: change of q to -q}. Since $\eps >0$ is arbitrary, and
\begin{gather*}
\|q_{N,s}^\# -\tilde q_{N,s}^\#\|_{L^2(\R_+)} = \|q_{N}- \tilde q_{N}\|_{L^2(\R)}\ge \|q- \tilde q\|_{L^2(\R)}/2,
\\
\rho_{S_{2}}(f_{q_{N,s}^\#}, f_{\tilde q_{N,s}^\#}) = \rho_{X}(\rc_{q_{N,s}}, \rc_{\tilde q_{N,s}}) \le \eps,
\end{gather*} 
this concludes the proof. \qed

\medskip

\noindent {\bf Proof of Theorem \ref{t5}.} The theorem is a combination of Proposition \ref{p57} and Proposition \ref{p514}.  \qed

\vspace{1cm}

\phantom{space}

\bibliographystyle{plain} 	
\bibliography{bibfile} 
\end{document}